\documentclass[11pt]{article}

\usepackage[T1]{fontenc}
\usepackage{lmodern}
\usepackage{microtype}
\usepackage{xcolor}
\usepackage{amsmath,amssymb,amsthm,mathtools}
\usepackage{bm}
\usepackage{mathrsfs}
\usepackage{enumitem}
\usepackage{booktabs,tabularx}
\usepackage{geometry}
\usepackage{aliascnt}
\usepackage{hyperref}
\usepackage[capitalize,noabbrev]{cleveref}

\allowdisplaybreaks

\hypersetup{
  pdftitle={Sharp Barycentric Costs in Teichm\"uller Space},
  pdfauthor={Hangyu Li},
  colorlinks=true,
  linkcolor=blue!55!black,
  citecolor=blue!55!black,
  urlcolor=blue!55!black
}

\newtheorem{theorem}{Theorem}[section]
\newaliascnt{proposition}{theorem}
\newtheorem{proposition}[proposition]{Proposition}
\aliascntresetthe{proposition}
\newaliascnt{lemma}{theorem}
\newtheorem{lemma}[lemma]{Lemma}
\aliascntresetthe{lemma}
\newaliascnt{corollary}{theorem}
\newtheorem{corollary}[corollary]{Corollary}
\aliascntresetthe{corollary}
\newaliascnt{claim}{theorem}

\aliascntresetthe{claim}
\theoremstyle{definition}
\newaliascnt{definition}{theorem}
\newtheorem{definition}[definition]{Definition}
\aliascntresetthe{definition}
\newaliascnt{assumption}{theorem}

\aliascntresetthe{assumption}
\newaliascnt{example}{theorem}

\aliascntresetthe{example}
\theoremstyle{remark}
\newaliascnt{remark}{theorem}

\aliascntresetthe{remark}

\crefname{proposition}{Proposition}{Propositions}
\crefname{lemma}{Lemma}{Lemmas}
\crefname{corollary}{Corollary}{Corollaries}
\crefname{definition}{Definition}{Definitions}
\crefname{example}{Example}{Examples}
\crefname{remark}{Remark}{Remarks}

\newcommand{\Bel}{\operatorname{Bel}}
\newcommand{\Teich}{\mathcal T}
\newcommand{\QD}{\mathcal Q}
\newcommand{\Bspace}{\mathcal B}
\newcommand{\Prob}{\mathbb P}
\newcommand{\E}{\mathbb E}
\newcommand{\D}{\mathbb D}
\newcommand{\distH}{d_{\mathbb H}}
\newcommand{\distB}{\mathscr D_2}
\newcommand{\WP}{\mathrm{WP}}
\newcommand{\id}{\operatorname{id}}

\newcommand{\Sym}{\operatorname{Sym}}
\newcommand{\supp}{\operatorname{supp}}
\newcommand{\barH}{\operatorname{bar}_{\mathbb H}}

\newcommand{\cF}{\mathcal F}
\newcommand{\Energy}{\mathscr E}

\newcommand{\cU}{\mathcal U}

\newcommand{\vH}{\mathbf S}

\title{Sharp Barycentric Costs\\in Teichm\"uller Space}
\author{Hangyu Li\\[0.3em]
{\small Department of Mathematics, The Chinese University of Hong Kong}\\
{\small\texttt{hyli@math.cuhk.edu.hk}}}
\date{}

\begin{document}

\maketitle

\begin{abstract}
The affine Beltrami fiber over a closed hyperbolic surface has a quadratic normal response whose least vertical \(L^2\) energy equals the tangent norm dual to the hyperbolic cotangent supremum norm. We establish this calibration by a curvature identity and exact localized responses. Its strictness determines the energy measures of near-minimizers; under an amplitude constraint, it yields a positive excess whose quadratic lower coefficient is explicit when the supporting differential has nondegenerate maxima. Pointwise Poincar\'e barycenters of representatives in one fiber provide a nonlinear realization: equal-weight pairs attain prescribed small outputs exactly, and their creation energy approaches the calibrated norm if and only if the output scale is negligible relative to the squared amplitude bound. Restoration to the full local fiber graph has leading cost four times the squared Weil--Petersson norm. For degenerating surfaces, we identify the additive limit of the calibrated norm across persistent components and collar directions. The induced length metric completes to marked augmented Teichm\"uller space and has distance \(\frac\pi2\sum_i\ell_i+O((\sum_i\ell_i)^2)\) to a fixed pinching stratum. Systole-dependent estimates also give exact nonlinear realizations along degenerating families.
\end{abstract}



\section{Introduction}\label{sec:introduction}

The hyperbolic supremum norm of holomorphic quadratic differentials defines a Finsler cometric on Teichm\"uller space. Its dual tangent norm is known as a member of the \(L^p\) family of cotangent-dual metrics \cite[\S2]{bridgeman2023convergence}. We identify this norm with a least-energy problem arising from the curvature of an affine Teichm\"uller fiber. The normal response to a perturbation tangent to that fiber is quadratic, and its least vertical \(L^2\) energy admits an exact cotangent calibration.

This identification has consequences at two scales. On a fixed surface, the calibration is strict: a near-minimizer concentrates on the maximum set of a supporting quadratic differential. A pointwise amplitude bound prevents arbitrarily fine concentration and forces an excess over the tangent norm. Along a pinching family, the corresponding cotangent unit ball separates into the persistent differentials on the limiting punctured surface and the principal parts of the shrinking collars. The dual norm becomes additive, with an explicit contribution from each pinched curve.

The quadratic response also governs a nonlinear problem. Representatives of the same marked conformal class can be aggregated pointwise by Poincar\'e barycenters of their Beltrami coefficients. The resulting coefficient need not remain in their common fiber. Exact two-point realizations of the prescribed output recover the sharp quadratic energy; the cost of returning the aggregate to the full local fiber graph is instead controlled by the Weil--Petersson norm. The nonlinear estimates retain two \(L^2\) factors under concentration and admit a uniform, systole-dependent version.

\subsection{Main results}

Let \(X\) be a closed hyperbolic surface of genus \(g\ge2\). Write
\(\pi:\Bel(X)\to\Teich(X)\) for the Teichm\"uller projection. With the harmonic complement to \(V=\ker D\pi_0\), the fiber through the origin is locally the graph
\[
\Phi(v)=v+\Psi(v),\qquad \Psi(0)=D\Psi_0=0.
\]
Its affine vertical Hessian is \(\vH_0\), and we set
\(q_X(v)=-\frac12\vH_0(v,v)\). All \(L^2\) norms in what follows use hyperbolic area.

\paragraph{Fiber curvature and exact calibration.}
The starting point is the variational norm
\[
\gamma_X(u)=\inf\{\E\|v\|_2^2:\E v=0,\ \E q_X(v)=u\},
\]
where the infimum is over finite centered laws of vertical directions. The curvature identity and localized point responses show that
\begin{equation}\label{eq:intro-cotangent-calibration}
\gamma_X^*(\phi)=\|\phi\|_{\infty,h}.
\end{equation}
The inequality \(|\langle q_X(v),\phi\rangle|\le
\|\phi\|_{\infty,h}\|v\|_2^2\) is strict for nonzero \(v\) and \(\phi\). Nevertheless, every \(u\) has smooth single-direction responses \(q_X(v)=u\) with energy arbitrarily close to \(\gamma_X(u)\) (\cref{lem:curvature-response,thm:cotangent-calibration,lem:radial-point-response,thm:single-direction-cost}).

\paragraph{Concentration and constrained energy.}
Fix \(u\ne0\) and a supporting differential \(\phi\) with \(\|\phi\|_{\infty,h}=1\). The possible limits of energy measures of near-minimizing laws are precisely the positive measures \(\nu\) with
\begin{equation}\label{eq:intro-contact-measure}
\nu(X)=\gamma_X(u),\qquad
\supp\nu\subset Z_\phi:=\{|\phi|_h=1\},\qquad
u=\int_{Z_\phi}e_X(p,a_\phi(p))\,d\nu(p),
\end{equation}
where \(e_X\) is cotangent point evaluation and
\(a_\phi=\overline\varphi/|\varphi|\). Every such measure is obtained from smooth single directions with exact quadratic response (\cref{thm:energy-measure-description}).

The convexified amplitude-one response body defines \(J_X\). For \(\gamma_X(u)=1\), its excess satisfies
\(\lambda+\kappa_X\lambda^3\le J_X(\lambda u)
\le\lambda+C_{X,u}\lambda^2\). If a supporting \(\phi\) has only nondegenerate maxima, then
\begin{equation}\label{eq:intro-capped-excess}
\liminf_{\lambda\downarrow0}
\frac{J_X(\lambda u)-\lambda}{\lambda^2}
\ge\left(4\pi\sum_{p\in Z_\phi}
\bigl(\det\operatorname{Hess}_h(1-|\phi|_h)(p)\bigr)^{-1/2}\right)^{-1}.
\end{equation}
Nondegenerate global maxima occur on an open dense subset of \(\QD(X)\setminus\{0\}\) for each fixed \(X\). This is a genericity statement about \emph{cotangent differentials}; for a given tangent direction, the estimate applies when one of its supporting differentials has the stated property. The coefficient in \eqref{eq:intro-capped-excess} is an explicit lower bound for the quadratic excess (\cref{thm:capped-excess,thm:nondegenerate-cotangents,prop:explicit-excess-coefficient}).

\paragraph{Exact barycentric realization and restoration.}
Let \(\tau_X\) be the normalized Bers chart and let
\(\mathcal I_X(t;r)\) and \(\mathcal K_X(t;r)\) denote, respectively, the least vertical energy and least integrated Poincar\'e dispersion among centered laws of vertical amplitude at most \(r\) with exact barycentric output \(t\). At each fixed \(X\), for interior feasible \(x\),
\begin{equation}\label{eq:intro-capped-limit}
\frac{\mathcal I_X(\tau_X^{-1}(r^2x);r)}{r^2}
\longrightarrow J_X(x),
\end{equation}
locally uniformly in \(x\). For fixed \(u\ne0\), and for positive \(r_j,\delta_j\to0\), asymptotically sharp creation at exact output \(t_j=\tau_X^{-1}(\delta_j u)\) is possible if and only if \(\delta_j=o(r_j^2)\). On this scale,
\begin{equation}\label{eq:intro-exact-cost}
\frac{\mathcal I_X(t_j;r_j)}{\delta_j}\to\gamma_X(u),
\qquad
\frac{\mathcal K_X(t_j;r_j)}{4\delta_j}\to\gamma_X(u).
\end{equation}
Equal-weight pairs supported on \(\{\Phi(v_j),\Phi(-v_j)\}\) attain the prescribed nonlinear output exactly and approach these infima. For a centered exact-output law, the least excess Fr\'echet energy over the \emph{full local fiber graph} \(\Phi(U_X)\) satisfies
\begin{equation}\label{eq:intro-uniform-restoration}
\Delta_X(\Pi)=4\delta^2\|u\|_{\WP}^2(1+O_X(r^2)).
\end{equation}
The two leading costs therefore arise from different norms, even for the same exact pair (\cref{thm:finite-radius-limit,thm:sharp-exact-cost,cor:sharp-cost-scale,prop:uniform-exact-restoration}). Systole-dependent graph, law and localization estimates produce exact pairs along degenerating families under the sufficient ranges of \cref{thm:uniform-exact-realization,cor:prescribed-degenerating-scales}. The fixed-surface necessity of \(\delta=o(r^2)\) and the sufficient systole-dependent range are recorded separately.

\paragraph{Boundary geometry of the calibrated norm.}
The same norm has a definite limit at the boundary. Suppose the marked surfaces \(X_j\) converge by pinching \(k\) disjoint curves to a finite-area punctured surface \(\Sigma\), with collar differentials normalized by \(b_0^i(\chi_{i,j})=\ell_{i,j}^2\). In the resulting persistent and collar dual coordinates \((u_0,c)\),
\begin{equation}\label{eq:intro-additive-degeneration}
\gamma_{X_j}(u_j)\longrightarrow
\gamma_\Sigma(u_0)+\frac1{4\pi^2}\sum_{i=1}^k|c_i|.
\end{equation}
The length metric \(d_\gamma\) has the same Cauchy sequences as the Weil--Petersson metric and its completion is marked augmented Teichm\"uller space. For a fixed pinching multicurve \(\Gamma\), provided the remaining Fenchel--Nielsen data stay in a compact set away from further degeneration,
\begin{equation}\label{eq:intro-pinching-distance}
d_\gamma(X,D_\Gamma)=\frac\pi2\sum_i\ell_i
+O\!\left(\left(\sum_i\ell_i\right)^2\right).
\end{equation}
With angular twist normalized to have Dehn period \(2\pi\), its coordinate vector has \(\gamma\)-norm asymptotic to \(\ell_i^2/(4\pi)\). The collar limit, the marked completion, and the radial distance law are proved in \cref{thm:additive-cost-degeneration,cor:fn-creation-scales,thm:creation-completion,thm:linear-pinching-distance}.

\subsection{Relation to earlier work}

Bridgeman--Bromberg--Vargas Pallete introduce the hyperbolic \(L^p\) cotangent norms and their dual tangent metrics \cite[\S2]{bridgeman2023convergence}. The norm \(\gamma\) is their \(p=\infty\) tangent metric. The curvature minimization gives its fiber-geometric characterization, while strict calibration and exact localized responses determine the energy measures of near-minimizers and the constrained nonlinear cost. The additive collar formula, marked completion and pinching-distance coefficient are statements about the global geometry of its induced length metric. Their proofs combine the calibrated norm with degenerating cotangent bases and global length-function estimates.

The affine Beltrami model and its infinitesimal quotient are classical \cite{ahlfors1960riemann,fletcher2006quasiconformal,gardiner2000quasiconformal,hubbard2016teichmuller,nag1988complex}. Connections arising from quotient geometry and the Weil--Petersson metric are developed in \cite{fischer1982almost,tromba1986natural,wolpert1986chern}. Miyachi studies second-order tangent and cotangent structures and the infinitesimal duality of the cotangent \(L^1\) norm with the Teichm\"uller norm \cite[Chs.~8--11]{miyachi20241}. Here the fixed affine fiber supplies a quadratic normal response whose vertical energy calibrates the hyperbolic cotangent supremum norm.

The nonlinear estimates retain two \(L^2\) factors under concentration. They use square-integrable Beltrami theory \cite{cui2000integrably,takhtajan2006weil}, Bers variation \cite{guo2003geodesic}, and Beurling resolvents \cite{astala2008elliptic}. Disjoint local fields give a coordinate-square map; its degree provides exact correction \cite{deimling2013nonlinear,arutyunov2016properties}. At the boundary, persistent and collar estimates come from \cite{rupflin2013asymptotics,grosse2019holomorphic}. Length differentials and pants geometry \cite{gardiner1975schiffer,wolpert2008behavior,alessandrini2012behaviour,buser2010geometry} supply global calibrations and marked compactness. Wolpert's stratum-distance expansion gives the comparison with squared Weil--Petersson distance \cite{wolpertgeometry}. Ahlfors--Weill reflection controls the uniform nonlinear range \cite{ahlfors1962uniqueness,nag1987conformally}.

Huang--Ohshika--Pan--Papadopoulos identify the earthquake forward FD-completions and the completions of their symmetrizations with augmented Teichm\"uller space \cite[Theorems~1.21 and 10.27]{huang2024earthquake}. Their argument also uses pants lengths and marked transversals. Their norm assigns lamination length to an infinitesimal earthquake; its distance to the full boundary is asymptotic to \(2\ell_{\mathrm{sys}}\log(1/\ell_{\mathrm{sys}})\) \cite[Theorem~1.29 and Remark~1.30]{huang2024earthquake}. The symmetric norm considered here is dual to the hyperbolic cotangent supremum and has the linear distance \eqref{eq:intro-pinching-distance} to a fixed pinching stratum.

Fr\'echet and Riemannian means supply the pointwise barycenter theory \cite{afsari2011riemannian,frechet1948elements,karcher1977riemannian,sturm2003probability} and local expansions \cite{bhattacharya2003large,pennec2006intrinsic,pennec2019curvature}. At the conformal basepoint, the barycenter agrees with the affine mean through second order, so the leading response is determined by the fiber Hessian; the Poincar\'e model specifies the exact nonlinear aggregate and its loss. Barycentric quasiconformal constructions appear in \cite{douady1986conformally,earle2002barycentric}. Lin studies energy minimization in a boundary-relative Teichm\"uller class and cancellation under dyadic random concatenation \cite{lin2023dyadic}; our constraint prescribes the output of pointwise aggregation of representatives on a closed surface.

\subsection{Organization}

Section~\ref{sec:fiber} derives the curvature identity, exact calibration and point responses. Section~\ref{sec:concentration} treats concentration, bounded amplitude and nondegenerate supporting differentials. Section~\ref{sec:constrained} proves exact realization and restoration. Section~\ref{sec:uniform-scales} establishes uniform exact realization. Section~\ref{sec:cost-degeneration} then derives the additive boundary norm, identifies the marked completion and computes the distance to pinching strata. The appendices contain the nonlinear Bers estimates and the uniform graph and localization proofs.

\section{Fiber curvature and exact calibration}\label{sec:fiber}\label{sec:sharp-cost}

The quadratic normal response of the affine fiber determines its least creation energy. We derive the response from hyperbolic curvature, calibrate it against cotangent point evaluation, and realize each evaluation atom exactly.

\subsection{Beltrami coefficients and normalizations}\label{sec:deterministic}

Fix a closed Riemann surface \(X\) of genus \(g\ge2\), with hyperbolic metric \(ds_X^2=\rho_X^2|dz|^2\) and area \(dA_X=\rho_X^2dx\,dy\). Write
\[
\Bspace(X)=L^\infty(X,\overline K_X\otimes K_X^{-1}),\qquad
\Bel(X)=\{\mu\in\Bspace(X):\|\mu\|_\infty<1\}.
\]
The normalized Beltrami equation defines the holomorphic split submersion \(\pi:\Bel(X)\to\Teich(X)\), with \(\pi(0)=[X]\) and local holomorphic sections \cite{ahlfors1960riemann,gardiner2000quasiconformal,hubbard2016teichmuller,nag1988complex}.

\paragraph{The infinitesimal quotient}
For \(\phi=\varphi dz^2\in\QD(X)\) and \(\mu=\mu(z)d\bar z/dz\), fix the coordinate-invariant pairing
\begin{equation}\label{eq:pairing}
\langle\mu,\phi\rangle=\int_X\mu\varphi\,dx\,dy
=\frac{i}{2}\int_X\mu\varphi\,dz\wedge d\bar z.
\end{equation}
The infinitesimal equivalence theorem identifies
\begin{align}
V&:=\ker D\pi_0,\label{eq:vertical}\\
V&=\QD(X)^\perp
=\{v\in\Bspace(X):\langle v,\phi\rangle=0\text{ for every }\phi\in\QD(X)\}.\label{eq:annihilator}
\end{align}
Thus \(T_{[X]}\Teich(X)\simeq\Bspace(X)/V\), \(T_{[X]}^*\Teich(X)\simeq\QD(X)\), and \(\dim_{\mathbb C}\QD(X)=3g-3\), with induced pairing
\begin{equation}\label{eq:tangent-pairing}
\langle D\pi_0[\mu],\phi\rangle_T=\langle\mu,\phi\rangle.
\end{equation}

\paragraph{Harmonic representatives}\label{subsec:harmonic}

The harmonic representative of \(\phi=\varphi dz^2\) is
\begin{equation}\label{eq:harmonic-beltrami}
h_\phi=\frac{\overline\varphi}{\rho_X^2}\frac{d\bar z}{dz}.
\end{equation}
For \(H=\{h_\phi:\phi\in\QD(X)\}\), nondegeneracy of \eqref{eq:pairing} gives
\begin{equation}\label{eq:splitting}
\Bspace(X)=V\oplus H,\qquad \dim_{\mathbb C}H=3g-3.
\end{equation}
With \(\langle\mu,\nu\rangle_2=\int_X\mu\overline\nu\,dA_X\), the identity \(\langle v,h_\phi\rangle_2=\langle v,\phi\rangle=0\) gives
\begin{equation}\label{eq:orthogonal}
V\perp_{L^2}H.
\end{equation}
Our Weil--Petersson normalization is \cite{gardiner2000quasiconformal,wolf2012weil,wolpert1986chern}
\begin{equation}\label{eq:wp-convention}
\|D\pi_0[h]\|_{\WP}=\|h\|_{L^2(dA_X)},\qquad h\in H.
\end{equation}
The conventions used in the cost and boundary calculations are collected in \cref{tab:normalizations}. Small-noise parameters \(\varepsilon\) are real.

\begin{table}[htbp]
\centering
\small
\caption{Norm and twist conventions. The tensor convention in \cite{grosse2019holomorphic} concerns complex quadratic tensors. Twist derivatives hold length and the other Fenchel--Nielsen coordinates fixed.}\label{tab:normalizations}
\begin{tabularx}{\textwidth}{@{}p{0.25\textwidth}X@{}}
\toprule
Quantity & Convention \\
\midrule
Cotangent norms & \(\|\phi\|_{\WP,*}^2=\int_X|\varphi|^2\rho_X^{-4}\,dA_X\),\quad \(\|\phi\|_{\infty,h}=\max_X|\varphi|\rho_X^{-2}\). \\
Metric variation & \(\dot g=\operatorname{Re}(\upsilon dz^2)\) corresponds to \(h=\overline\upsilon/(2\rho_X^2)\). \\
Complex tensor norm & \(\|\phi\|_{\mathrm{tensor},2}^2=4\|\phi\|_{\WP,*}^2\). \\
Real twist \(s_i\) & A positive Dehn twist adds \(\ell_i\). \\
Angular twist \(\theta_i\) & \(\theta_i=2\pi s_i/\ell_i\); a positive Dehn twist adds \(2\pi\), and \(\partial_{\theta_i}=\ell_i\partial_{s_i}/(2\pi)\). \\
\bottomrule
\end{tabularx}
\end{table}

\subsection{The local fiber and its quadratic response}

Let \(\cF_0=\pi^{-1}([X])\). The second-order response of aggregation is determined by the curvature of this fiber in affine Beltrami space. Its second fundamental form is represented by the vertical Hessian of \(\pi\); the harmonic complement selects a representative of each normal class.

\paragraph{The vertical Teichm\"uller Hessian}

For the infinite-dimensional vertical space \(V\), \(\Sym^2(V)\) denotes the algebraic symmetric square; equivalently, the tensor below may be read as a continuous symmetric complex-bilinear map on \(V\times V\).

For a holomorphic chart \(\tau\) centered at \([X]\), put \(C=\tau\circ\pi\).

\begin{definition}[Vertical Teichm\"uller Hessian]\label{def:vertical-hessian}
For \(v,w\in V\), define
\begin{equation}\label{eq:vertical-hessian}
  \vH_0(v,w)
  :=(D\tau_{[X]})^{-1}
  D^2(\tau\circ\pi)_0[v,w]
  \in T_{[X]}\Teich(X).
\end{equation}
\end{definition}

Here the second derivative is taken in the canonical affine coordinates of \(\Bspace(X)\).  Restriction to vertical directions removes the coordinate-dependent quadratic term in the target chain rule.

The vertical restriction makes \eqref{eq:vertical-hessian} independent of the target chart: in the second-order chain rule, the additional term vanishes because \(D\pi_0v=D\pi_0w=0\).

\paragraph{Local analytic graph of the fiber}

The harmonic splitting gives an isomorphism \(DC_0|_H:H\to\mathbb C^{3g-3}\). The holomorphic implicit-function theorem therefore represents the fiber as a graph over \(V\).

\begin{proposition}[Local analytic fiber graph]\label{prop:fiber-graph}
There are neighborhoods \(0\in U_V\subset V\), \(0\in U_H\subset H\), and a holomorphic map \(\Psi:U_V\to U_H\), with \(\Psi(0)=D\Psi_0=0\), such that
\begin{equation}\label{eq:fiber-graph}
\cF_0\cap(U_V\oplus U_H)
=\{\Phi(v):=v+\Psi(v):v\in U_V\}.
\end{equation}
Writing \(B=D^2\Psi_0:\Sym^2(V)\to H\), on a smaller \(U_V\) one has
\begin{equation}\label{eq:Psi-expansion}
\Psi(v)=\tfrac12B(v,v)+R_\Psi(v),\qquad
\|R_\Psi(v)\|_\infty\le C\|v\|_\infty^3.
\end{equation}
\end{proposition}

\begin{proof}
The holomorphic implicit-function theorem \cite[Appendix~A.5]{hubbard2016teichmuller} applies to \((v,h)\mapsto C(v+h)\), since \(DC_0|_H\) is invertible. Differentiating \(C(v+\Psi(v))=0\) at zero gives \(D\Psi_0=0\); the Banach-holomorphic Taylor expansion on a smaller ball gives \eqref{eq:Psi-expansion}.
\end{proof}

The two second-order objects are related by differentiating the fiber identity.

\begin{proposition}[Quadratic graph identity]\label{prop:B-S}
For all \(v,w\in V\),
\begin{equation}\label{eq:B-S}
  D\pi_0[B(v,w)]
  =-\vH_0(v,w).
\end{equation}
\end{proposition}

\begin{proof}
Twice differentiating the fiber identity gives \(DC_0[B(v,w)]+D^2C_0[v,w]=0\). Apply \((D\tau_{[X]})^{-1}\).
\end{proof}

In the affine coefficient space the quotient-valued second fundamental form is \(\mathrm{II}_0(v,w)=[B(v,w)]\bmod V\). Under the quotient isomorphism \(\overline{D\pi}_0\), \eqref{eq:B-S} reads
\begin{equation}\label{eq:II-S}
\overline{D\pi}_0(\mathrm{II}_0(v,w))=-\vH_0(v,w).
\end{equation}
The harmonic complement selects the representative \(B(v,w)\) of this normal class.

\subsection{Localized responses}\label{sec:nontrivial}

We determine the displacements that the quadratic fiber response can produce by localizing vertical directions on coordinate disks. Compactly supported self-isotopies identify the response with tangent--cotangent pairings. Disjoint supports then make it diagonal in a tangent basis.

\begin{lemma}[Compactly supported isotopy response]\label{lem:local-isotopy-response}
Let \(U\Subset X\) be a coordinate disk and \(W\in C_c^\infty(U,\mathbb C)\). For small real \(t\), the coordinate perturbation \(F_t(z)=z+tW(z)\), extended by the identity, is an orientation-preserving diffeomorphism isotopic to the identity. Extending coefficients by zero, put
\[
v=W_{\bar z}\frac{d\bar z}{dz},\qquad
A_W=W_zW_{\bar z}\frac{d\bar z}{dz}.
\]
Then \(v\in V\) and
\begin{equation}\label{eq:local-response}
\tfrac12\vH_0(v,v)=D\pi_0[A_W].
\end{equation}
\end{lemma}

\begin{proof}
For small \(|t|\), each \(F_{st}\), \(0\le s\le1\), is a local diffeomorphism of degree one, hence a global diffeomorphism; this path gives the isotopy. Its coefficient is
\[
\mu_t=\frac{tW_{\bar z}}{1+tW_z}
=tv-t^2A_W+O_{L^\infty}(t^3).
\]
For \(\phi=\varphi dz^2\), compact support and holomorphicity give \(\langle v,\phi\rangle=-\int_UW\partial_{\bar z}\varphi\,dx\,dy=0\), so \(v\in V\). Since \(\pi(\mu_t)=[X]\), the quadratic term of \(C(\mu_t)=0\) is \(D^2C_0[v,v]=2DC_0[A_W]\), proving \eqref{eq:local-response}.
\end{proof}

For real-valued \(W\), the coefficient \(W_zW_{\bar z}\) in \cref{lem:local-isotopy-response} is the nonnegative density \(|W_z|^2\). We use this observation to localize the tangent--cotangent pairing.

\begin{lemma}[Localized spanning responses]\label{lem:localized-basis}
Let \(n=3g-3\). In any prescribed nonempty open set \(O\subset X\), there are pairwise disjoint coordinate disks \(U_i\Subset O\) and real functions \(W_i\in C_c^\infty(U_i)\) such that
\[
v_i=(W_i)_{\bar z_i}\frac{d\bar z_i}{dz_i}\in V,
\qquad Q_i=\tfrac12\vH_0(v_i,v_i),\qquad 1\le i\le n,
\]
and \(Q_1,\ldots,Q_n\) form a complex basis of \(T_{[X]}\Teich(X)\).
\end{lemma}

\begin{proof}
Choose distinct points \(p_i\in O\) so that evaluation on them is an isomorphism on \(\QD(X)\). They can be chosen inductively: a nonzero differential in the common kernel of earlier evaluations is nonzero somewhere in \(O\) by the identity theorem, and evaluation there lowers the kernel dimension by one. Choose disjoint disks \(U_i\Subset O\), with \(z_i(p_i)=0\), and a basis \(q_j=\varphi_j^{(i)}dz_i^2\) of \(\QD(X)\). The matrix \((\varphi_j^{(i)}(0))\) is invertible.

Fix a nonconstant real \(\chi\in C_c^\infty(\D)\), put \(W_{i,r}=r\chi(z_i/r)\), and extend \(A_{i,r}=|(W_{i,r})_{z_i}|^2d\bar z_i/dz_i\) by zero. Rescaling gives
\begin{equation}\label{eq:localized-pairing-limit}
r^{-2}\langle A_{i,r},q_j\rangle
\longrightarrow c_\chi\varphi_j^{(i)}(0),\qquad
c_\chi=\int_\D|\chi_\zeta|^2\,dA>0.
\end{equation}
Thus the pairing matrix is invertible for small \(r\). Fix such an \(r\) and take \(W_i=W_{i,r}\). By \cref{lem:local-isotopy-response}, \(Q_i=D\pi_0[A_{i,r}]\), so \eqref{eq:tangent-pairing} proves the assertion.
\end{proof}

\begin{theorem}[Full quadratic range]\label{thm:full-range}
The map
\begin{equation}\label{eq:quadratic-response-map}
Q_0:V\to T_{[X]}\Teich(X),\qquad Q_0(v)=\tfrac12\vH_0(v,v),
\end{equation}
is surjective already on the \(n\)-dimensional complex space \(E_0=\operatorname{span}_{\mathbb C}\{v_1,\ldots,v_n\}\) in \cref{lem:localized-basis}. The induced linear map \(\vH_0:\Sym^2(E_0)\to T_{[X]}\Teich(X)\) is also surjective.
\end{theorem}

\begin{proof}
For \(W=\sum_i c_iW_i\), disjoint supports give \(W_zW_{\bar z}=\sum_i c_i^2(W_i)_z(W_i)_{\bar z}\). The isotopy response therefore yields
\begin{equation}\label{eq:diagonal-additivity}
Q_0\!\left(\sum_i c_iv_i\right)=\sum_i c_i^2Q_i.
\end{equation}
Given \(u=\sum_i\alpha_iQ_i\), choose \(c_i^2=\alpha_i\). This proves both surjectivity statements.
\end{proof}

Write \(T_X=T_{[X]}\Teich(X)\), \(n=3g-3\), and
\[
q_X(v)=-\tfrac12\vH_0(v,v)=-Q_0(v),\qquad
V_2=H^{\perp_{L^2}}\subset L^2(X;dA_X).
\]
The sign comes from \eqref{eq:B-S}. The response is first defined on \(V\); the following identity extends it to \(V_2\).

\subsection{The curvature identity}

For \(\phi=\varphi dz^2\in\QD(X)\), set
\begin{equation}\label{eq:cotangent-sup}
\|\phi\|_{\infty,h}=\max_{p\in X}\frac{|\varphi(p)|}{\rho_X(p)^2}.
\end{equation}
This is the hyperbolic \(L^\infty\) cotangent norm in the \(L^p\) family described in \cite[\S2]{bridgeman2023convergence}. In the following formulas write \(\rho=\rho_X\), \(dA=dA_X\).

For a smooth vector field with local complex coefficient \(W\), write
\[
v=\bar\partial W=W_{\bar z}\frac{d\bar z}{dz},
\qquad
NW=W_z+2(\log\rho)_zW.
\]
The scalar \(NW\) is invariant under holomorphic coordinate changes. Define
\[
\|W\|_h^2=\int_X\rho^2|W|^2\,dA.
\]
The harmonic representative description of the tangent space \cite[\S7.7]{hubbard2016teichmuller} identifies the image of \(\bar\partial\) with \(V_2\). More explicitly, the operator \(\bar\partial:H^1(TX)\to V_2\) is an isomorphism. Its kernel consists of holomorphic vector fields, which vanish for \(g\ge2\). Its elliptic range is closed. The adjoint equation is locally \(\partial_z(\rho^2a)=0\), so its kernel is the harmonic Beltrami space \(H\). Thus every \(v\in V_2\) has a unique \(H^1\) potential \(W\).

\begin{lemma}[Covariant response and curvature identity]\label{lem:curvature-response}
For \(v=\bar\partial W\),
\begin{equation}\label{eq:vertical-bochner}
\|NW\|_2^2+\tfrac12\|W\|_h^2=\|v\|_2^2.
\end{equation}
For every \(\phi=\varphi dz^2\in\QD(X)\),
\begin{equation}\label{eq:covariant-response}
\langle q_X(v),\phi\rangle
=-\int_X\varphi
\left((NW)v+\tfrac14\rho^2W^2\right)\,dx\,dy.
\end{equation}
Both formulas extend from smooth fields to \(H^1(TX)\), with the right side of \eqref{eq:covariant-response} defining the extended quadratic response.
\end{lemma}

\begin{proof}
Put \(\sigma=\log\rho\). Curvature \(-1\) gives \(\sigma_{z\bar z}=\rho^2/4\). With the source inner product weighted by \(\rho^4\,dx\,dy\) and both derivative target inner products weighted by \(\rho^2\,dx\,dy\), integration by parts gives
\[
N^*f=-\rho^{-2}f_{\bar z},
\qquad
\bar\partial^*a=-\rho^{-2}(a_z+2\sigma_z a).
\]
Consequently,
\[
N^*NW-\bar\partial^*\bar\partial W
=-2\rho^{-2}\sigma_{z\bar z}W=-\tfrac12W.
\]
Pairing with \(W\) proves \eqref{eq:vertical-bochner}.

For the response formula, use \(F_t(x)=\exp_x(tW_x)\). For small real \(t\), this is a diffeomorphism isotopic to the identity. In a conformal coordinate its expansion is
\[
F_t(z)=z+tW-t^2\sigma_zW^2+O(t^3),
\]
since the geodesic equation is \(\ddot z+2\sigma_z\dot z^2=0\). Its Beltrami coefficient is therefore
\[
\mu_t=tv-t^2\left((NW)v+\sigma_{z\bar z}W^2\right)+O(t^3).
\]
The second-order coefficient is a global Beltrami differential. Since \(\pi(\mu_t)=[X]\), comparison of quadratic terms gives
\[
\tfrac12\vH_0(v,v)
=D\pi_0\left[(NW)v+\tfrac14\rho^2W^2\right].
\]
This proves \eqref{eq:covariant-response} for smooth verticals. To identify it with the original Hessian on all of \(V\), smooth \(v\) and subtract its harmonic projection. This gives smooth \(v_j\in V\), uniformly bounded in \(L^\infty\), converging to \(v\) almost everywhere and in \(L^2\). Their lifts \(v_j^\sharp\), extended by zero outside \(\D\), converge in Euclidean \(L^2(\mathbb C)\) by dominated convergence.

Let \(\mathscr H\) denote the Beurling transform and \(M_j\) multiplication by \(v_j^\sharp\). Choose \(R\ge1\) bounding all \(\|v_j\|_\infty\) and \(\|v\|_\infty\). On the common parameter disk \(|t|<(2R)^{-1}\), the principal solution densities have the Beurling-resolvent representation \cite[\S\S4.5.2, 5.3, 5.7]{astala2008elliptic}
\begin{equation}\label{eq:vertical-resolvent-continuity}
\omega_j(t)=(I-tM_j\mathscr H)^{-1}(t v_j^\sharp),
\qquad \|(I-tM_j\mathscr H)^{-1}\|_{2\to2}\le2.
\end{equation}
The operators \(M_j\) converge strongly on \(L^2(\mathbb C)\). The resolvent identity therefore gives strong convergence of \(\omega_j(t)\) to the density for \(tv^\sharp\). Their Cauchy transforms and derivatives converge locally off \(\overline\D\), so the Bers projections of \(tv_j\) converge to that of \(tv\). Holomorphic dependence and the common Bers bound give locally uniform convergence in \(t\) \cite{ahlfors1960riemann}. Cauchy's formula now yields \(q_X(v_j)\to q_X(v)\).

Elliptic inversion also gives convergence of the potentials in \(H^1\), so the right side of \eqref{eq:covariant-response} has the same limit. Thus the formula holds for the original response on \(V\). Smooth approximation and Cauchy--Schwarz define its extension to \(V_2\).
\end{proof}

\begin{theorem}[\(L^2\)-tame second fundamental form]\label{thm:l2-tame}
There exists a constant \(C_X<\infty\), depending only on the base surface, such that
\begin{equation}\label{eq:l2-tame}
  \|\vH_0(v,w)\|_{\WP}
  \le C_X\|v\|_{L^2(dA_X)}\|w\|_{L^2(dA_X)},
  \qquad v,w\in V.
\end{equation}
Consequently, \(\vH_0\) extends uniquely to a continuous symmetric complex-bilinear map
\begin{equation}\label{eq:l2-extension}
  \vH_0^{(2)}:
  V_2\times V_2
  \longrightarrow T_{[X]}\Teich(X).
\end{equation}
\end{theorem}

\begin{proof}
Elliptic inversion of \(\bar\partial\) and \eqref{eq:vertical-bochner} bound both derivative factors and the potential term in \eqref{eq:covariant-response} by \(C_X\|v\|_2^2\). Cotangent duality and polarization give \eqref{eq:l2-tame}. More explicitly, if the diagonal quadratic form is bounded by \(C_X\|v\|_2^2\), the identity
\[
4b(v,w)=q_X(v+w)-q_X(v-w),\qquad b(v,w)=-\tfrac12\vH_0(v,w),
\]
and rescaling \(v,w\) give a bound by \(C_X\|v\|_2\|w\|_2\) for its symmetric bilinear form. Smooth fields in \(V\) are dense in \(V_2\): subtract the harmonic projection from any smooth \(L^2\) approximation. This proves the unique extension.
\end{proof}

\subsection{The creation norm}

\begin{definition}[Directional quadratic cost]\label{def:directional-cost}
For \(u\in T_X\), set
\begin{equation}\label{eq:directional-cost}
\gamma_X(u):=\inf\left\{\E\|\xi\|_2^2:
\ \E\xi=0,\ \E q_X(\xi)=u\right\}.
\end{equation}
The random variables are strongly measurable and essentially bounded in \(V\). There is no prescribed bound on \(\|\xi\|_\infty\) in this quadratic problem.
\end{definition}

Define the response body
\begin{equation}\label{eq:response-body}
\mathcal A_X:=\{q_X(v):v\in V,\ \|v\|_2\le1\},
\qquad
\mathcal B_X:=\overline{\operatorname{conv}}\,\mathcal A_X.
\end{equation}
The convex hull is taken in the real space underlying \(T_X\). Its gauge is an atomic norm \cite{chandrasekaran2012convex}.

\begin{proposition}[Response norm and cotangent dual]\label{thm:directional-norm}
The function \(\gamma_X\) is a complex norm with unit ball \(\mathcal B_X\), and \(\gamma_X(u)\asymp_X\|u\|_{\WP}\). Its value is unchanged when the infimum in \eqref{eq:directional-cost} is taken over strongly measurable, centered \(V_2\)-valued laws with finite second moment. For \(\lambda\in T_X^*\),
\begin{equation}\label{eq:directional-dual}
\gamma_X^*(\lambda)
=\sup_{\substack{v\in V_2\\\|v\|_2\le1}}|\lambda(q_X(v))|,
\end{equation}
and
\begin{equation}\label{eq:directional-polar}
\gamma_X(u)=\max_{\gamma_X^*(\lambda)\le1}\operatorname{Re}\lambda(u).
\end{equation}
\end{proposition}

\begin{proof}
The bilinear bound makes \(\mathcal A_X\) bounded. Since \(q_X(\sqrt\alpha\,v)=\alpha q_X(v)\), the body is balanced; the localized responses span \(T_X\), so it contains a neighborhood of zero. Its gauge \(g_X\) is a complex norm. For a law of energy \(M>0\), weighting by \(\|\xi\|_2^2/M\) gives
\[
\frac{\E q_X(\xi)}{M}
=\E_{\mathrm{weighted}}q_X\!\left(\frac{\xi}{\|\xi\|_2}\right)
\in\mathcal B_X,
\]
with zero samples omitted. Density of \(V\) in \(V_2\) makes this valid for both classes of laws, so \(g_X(u)\le M\).

If \(t>g_X(u)\), then \(u/t\) lies in the interior of \(\mathcal B_X\), hence in \(\operatorname{conv}\mathcal A_X\) \cite[Theorem~6.3]{rockafellar1997convex}. A finite convex representation
\[
u/t=\sum_k p_kq_X(a_k),\qquad \|a_k\|_2\le1,
\]
gives the centered law assigning probabilities \(p_k/2\) to \(\pm\sqrt t\,a_k\). Its output is \(u\) and its energy is at most \(t\). Letting \(t\downarrow g_X(u)\) proves \(\gamma_X=g_X\) and the Hilbert-space assertion. The zero law treats \(u=0\).

Polarity gives \eqref{eq:directional-polar}. Phase rotation of \(v\) changes \(q_X(v)\) by twice the phase, so the support function is \eqref{eq:directional-dual}. 
\end{proof}

\begin{theorem}[Exact cotangent calibration]\label{thm:cotangent-calibration}
For every \(\phi\in\QD(X)\),
\begin{equation}\label{eq:cotangent-identification}
\gamma_X^*(\phi)=\|\phi\|_{\infty,h}.
\end{equation}
If \(\phi\ne0\) and \(v\in V_2\setminus\{0\}\), then
\begin{equation}\label{eq:strict-calibration}
|\langle q_X(v),\phi\rangle|
<\|\phi\|_{\infty,h}\|v\|_2^2.
\end{equation}
\end{theorem}

\begin{proof}
Set \(C=\|\phi\|_{\infty,h}\), \(M=\|v\|_2^2\), and \(T=\|W\|_h^2\). The response formula gives
\begin{align*}
|\langle q_X(v),\phi\rangle|
&\le C\left(\|NW\|_2\|v\|_2+\tfrac14T\right)\\
&=CM\left(\sqrt{1-x}+\tfrac12x\right),
\qquad x=\frac{T}{2M}\in[0,1].
\end{align*}
For \(v\ne0\), elliptic inversion gives \(T>0\). Since \(\sqrt{1-x}+x/2<1\) for \(x>0\), this proves \eqref{eq:strict-calibration} and the upper bound in \eqref{eq:cotangent-identification}.

For the reverse bound, choose a coordinate centered at \(p\in X\), and a nonconstant real \(\chi\in C_c^\infty(\D)\). Put
\[
W_\varepsilon(z)=\varepsilon\chi(z/\varepsilon),
\qquad v_\varepsilon=(W_\varepsilon)_{\bar z}.
\]
Extend the field by zero. The local isotopy response in \cref{lem:local-isotopy-response} gives
\[
\langle q_X(v_\varepsilon),\phi\rangle
=-\int_X\varphi |(W_\varepsilon)_z|^2\,dx\,dy.
\]
Because \(W_\varepsilon\) is real, \(|(W_\varepsilon)_z|=|(W_\varepsilon)_{\bar z}|\). Rescaling yields
\[
\frac{\langle q_X(v_\varepsilon),\phi\rangle}{\|v_\varepsilon\|_2^2}
\longrightarrow-\frac{\varphi(p)}{\rho(p)^2}.
\]
By \eqref{eq:directional-dual}, taking a point where \(|\varphi|/\rho^2\) is maximal proves the reverse bound. Multiplication of \(v_\varepsilon\) by a phase gives every phase of this limiting response.
\end{proof}

For a point \(p\in X\) and a unit Beltrami coefficient \(a\) at \(p\), define \(e_X(p,a)\in T_X=\QD(X)^*\) by
\[
\langle e_X(p,a),\varphi dz^2\rangle
=a\frac{\varphi(p)}{\rho_X(p)^2}.
\]
The Beltrami transformation rule for \(a\) makes this definition coordinate independent. The atoms are parametrized by the unit circle bundle of the Beltrami line bundle, a compact real three-dimensional manifold.

\begin{corollary}[Geometric response body]\label{cor:geometric-response-body}
The unit ball of \(\gamma_X\) is
\begin{equation}\label{eq:geometric-response-body}
\mathcal B_X=\operatorname{conv}\{e_X(p,a):p\in X,\ |a|=1\}.
\end{equation}
It is not strictly convex. In particular, \(\gamma_X\) is not induced by any real inner product and cannot be a multiple of the Weil--Petersson norm.
\end{corollary}

\begin{proof}
The support function of the compact convex set on the right of \eqref{eq:geometric-response-body} is \(\|\phi\|_{\infty,h}\). The calibration theorem and finite-dimensional polarity identify it with \(\mathcal B_X\).

Every extreme point of the convex hull of a compact set belongs to that set, by Carath\'eodory's theorem \cite[Theorem~17.1]{rockafellar1997convex}. The smooth image \(\{e_X(p,a)\}\) has Hausdorff dimension at most three. If \(\mathcal B_X\) were strictly convex, every boundary point would be extreme. Its boundary would then lie in this image. The boundary of a convex body with interior in \(\mathbb R^{2n}\) is radially bi-Lipschitz equivalent to the unit sphere and has Hausdorff dimension \(2n-1\ge5\). This is impossible. An inner-product unit ball is strictly convex, giving the final assertion.
\end{proof}

\subsection{Exact point responses}

Radial displacement fields give exact point-evaluation responses. Their energy approaches the mass of the corresponding atom, with a quadratic error in the support radius.

\begin{lemma}[Radial point responses]\label{lem:radial-point-response}
Fix \(p\in X\) and a unit Beltrami coefficient \(a\) at \(p\). There are \(\varepsilon_p>0\) and \(C_p<\infty\) such that, for every \(m>0\) and \(0<\varepsilon<\varepsilon_p\), a smooth vertical direction \(v_{m,\varepsilon}\), supported in a coordinate disk about \(p\), satisfies
\begin{align}
q_X(v_{m,\varepsilon})&=m e_X(p,a),\label{eq:radial-exact-response}\\
m\le\|v_{m,\varepsilon}\|_2^2&\le m(1+C_p\varepsilon^2),\label{eq:radial-energy}\\
\|v_{m,\varepsilon}\|_\infty&\le C_p\sqrt m/\varepsilon.\label{eq:radial-amplitude}
\end{align}
The supporting disks shrink to \(p\) as \(\varepsilon\to0\).
\end{lemma}

\begin{proof}
Choose a disk uniformization centered at \(p\), so that its restriction is an embedded coordinate disk with
\[
z(p)=0,\qquad \rho(z)=\frac2{1-|z|^2}.
\]
Fix a nonconstant radial real function \(\chi\in C_c^\infty(\D)\), and set
\[
K_\chi=\int_\D|\chi_z|^2\,dx\,dy
=\int_\D|\chi_{\bar z}|^2\,dx\,dy>0.
\]
Choose \(b\in\mathbb C\) with \(b^2=-ma/(4K_\chi)\), and extend
\(W(z)=b\chi(z/\varepsilon)\) by zero. Put \(v_{m,\varepsilon}=\bar\partial W\). For \(\phi=\varphi dz^2\), the local isotopy formula gives
\[
\langle q_X(v_{m,\varepsilon}),\phi\rangle
=-\frac{b^2}{\varepsilon^2}
\int\varphi(z)|\chi_z(z/\varepsilon)|^2\,dx\,dy.
\]
The weight is radial. The circle mean-value identity for \(\varphi\) makes this integral equal to \(\varepsilon^2K_\chi\varphi(0)\). Since \(\rho(0)^2=4\), the result is \(ma\varphi(0)/4\), proving \eqref{eq:radial-exact-response}.

After rescaling,
\[
\|v_{m,\varepsilon}\|_2^2
=\frac m{K_\chi}\int_\D
\frac{|\chi_{\bar z}(\zeta)|^2}{(1-\varepsilon^2|\zeta|^2)^2}
\,dA_{\mathrm e}(\zeta).
\]
The weight \((1-\varepsilon^2|\zeta|^2)^{-2}\) lies between \(1\) and \(1+C\varepsilon^2\) on a fixed smaller coordinate disk. This proves \eqref{eq:radial-energy}. Differentiating \(W\) gives \eqref{eq:radial-amplitude}.
\end{proof}

\begin{theorem}[Single-direction cost]\label{thm:single-direction-cost}
For every \(u\in T_X\),
\begin{equation}\label{eq:single-direction-cost}
\inf\{\|v\|_2^2:v\in V,\ q_X(v)=u\}=\gamma_X(u).
\end{equation}
For every \(\eta>0\), there is a smooth vertical direction \(v=\bar\partial W\), supported in finitely many disjoint coordinate disks with a nonempty open complement, such that
\[
q_X(v)=u,\qquad \|v\|_2^2\le\gamma_X(u)+\eta.
\]
For \(u\ne0\), the infimum is not attained.
\end{theorem}

\begin{proof}
The symmetric law on \(\{v,-v\}\) gives the lower bound. For \(u\ne0\), write an optimal atomic representation supplied by \cref{cor:geometric-response-body} as
\[
u=\sum_{j=1}^Nm_j e_X(p_j,a_j),\qquad
m_j>0,\quad \sum_jm_j=\gamma_X(u).
\]
Take a normalized supporting differential \(\phi\) for \(u\). Equality in the support inequality forces every atom in this representation to satisfy \(\operatorname{Re}\langle e_X(p_j,a_j),\phi\rangle=1\). At a fixed point its phase is therefore unique, so repeated points can be merged.

Choose disjoint small disks about the resulting distinct points. Apply \cref{lem:radial-point-response} on each disk and sum the fields. Disjoint supports make both the quadratic responses and the energies additive. The output is exactly \(u\), while \eqref{eq:radial-energy} makes the energy tend to \(\sum_jm_j\). The disks may be chosen to leave an open complement. The zero field treats \(u=0\), and strict calibration gives nonattainment when \(u\ne0\).
\end{proof}

Strict calibration also excludes attainment by a centered square-integrable law with nonzero output. Pairing that output with a normalized supporting differential makes
\(\|v\|_2^2-\operatorname{Re}\langle q_X(v),\phi\rangle\)
strictly positive on every nonzero sample, so its expectation cannot vanish.

\begin{corollary}[Localized quadratic cost]\label{cor:cost-convexification}
If the laws in \eqref{eq:directional-cost} are restricted to \(E_0\), their cost is
\begin{equation}\label{eq:localized-sharp-cost}
\gamma_{X,E_0}\!\left(\sum_{i=1}^n\alpha_iQ_i\right)
=\sum_{i=1}^n|\alpha_i|\,\|v_i\|_2^2.
\end{equation}
A symmetric two-point vertical law attains this restricted quadratic cost.
For \(u\ne0\), one has \(\gamma_X(u)<\gamma_{X,E_0}(u)\).
\end{corollary}

\begin{proof}
For \(\xi=\sum_i c_iv_i\), disjoint supports give
\[
\E q_X(\xi)=-\sum_i\E c_i^2Q_i,
\qquad
\E\|\xi\|_2^2=\sum_i\E|c_i|^2\|v_i\|_2^2.
\]
Output \(\sum_i\alpha_iQ_i\) requires \(\E c_i^2=-\alpha_i\), so \(\E|c_i|^2\ge|\alpha_i|\). Choose square roots \(c_i^2=-\alpha_i\) and the law on \(\{\sum_i c_iv_i,-\sum_i c_iv_i\}\) to obtain equality.
For nonzero \(u\), this realizing direction is nonzero. Strict calibration against a normalized supporting cotangent vector gives \(\gamma_X(u)<\|\sum_i c_iv_i\|_2^2=\gamma_{X,E_0}(u)\).
\end{proof}

\section{Concentration and bounded-amplitude creation}\label{sec:concentration}

Strict calibration makes the quadratic infimum unattained. We describe the energy measures of all near-minimizing sequences and then retain the amplitude constraint in a finite-dimensional variational problem. Radial point responses give a matching quadratic excess when a supporting differential has nondegenerate maxima.

\subsection{The calibration deficit}

Fix \(\phi=\varphi dz^2\) with \(\|\phi\|_{\infty,h}=1\), and write
\[
f=\varphi/\rho^2,\qquad c=|f|,\qquad
D_\phi(v)=\|v\|_2^2-\operatorname{Re}\langle q_X(v),\phi\rangle.
\]
Expressions containing \(f\) below use the local coefficient in a conformal coordinate; each displayed integrand is invariant under coordinate change.

\begin{proposition}[Calibration deficit]\label{prop:calibration-deficit}
For \(v=\bar\partial W\in V_2\),
\begin{align}
D_\phi(v)
={}&\frac12\|NW+\overline f\,\overline v\|_2^2
+\frac12\int_X(1-c^2)|v|^2\,dA\notag\\
&+\frac14\int_X\rho^2
\bigl(|W|^2+\operatorname{Re}(fW^2)\bigr)\,dA.\label{eq:calibration-deficit-sum}
\end{align}
All three terms are nonnegative. For a square-integrable \(V_2\)-valued law, put
\[
M=\E\|v\|_2^2,\qquad
D=M-\operatorname{Re}\langle\E q_X(v),\phi\rangle,
\qquad d\nu=\E|v|^2\,dA.
\]
Then
\begin{align}
\int_X(1-c^2)\,d\nu&\le2D,\label{eq:deficit-spatial}\\
\E\|W\|_h^2&\le\sqrt{32MD}.\label{eq:deficit-potential}
\end{align}
Define \(\mathcal E_\phi(\nu)\in T_X\) by
\[
\langle\mathcal E_\phi(\nu),\psi\rangle
=\int_X\overline f\,\frac{\psi}{\rho^2}\,d\nu,
\qquad \psi\in\QD(X),
\]
where \(\psi/\rho^2\) denotes its local coefficient divided by \(\rho^2\). One has
\begin{equation}\label{eq:deficit-vector}
\gamma_X\bigl(\E q_X(v)-\mathcal E_\phi(\nu)\bigr)
\le2\sqrt{2MD}.
\end{equation}
\end{proposition}

\begin{proof}
Expand the first term in \eqref{eq:calibration-deficit-sum} and use
\(\|NW\|_2^2=\|v\|_2^2-\tfrac12\|W\|_h^2\).
The response formula \eqref{eq:covariant-response} gives the identity. The last integrand is nonnegative because \(|f|\le1\). Approximation in \(H^1(TX)\) proves the identity on \(V_2\), and integration proves \eqref{eq:deficit-spatial}.

For an individual nonzero sample, set \(M_v=\|v\|_2^2\), \(T_v=\|W\|_h^2\), and \(x=T_v/(2M_v)\). The proof of \cref{thm:cotangent-calibration} and
\[
1-\sqrt{1-x}-x/2\ge x^2/8,\qquad 0\le x\le1,
\]
give \(D_\phi(v)\ge T_v^2/(32M_v)\). Cauchy--Schwarz after expectation yields \eqref{eq:deficit-potential}; zero samples contribute zero.

Let \(R=NW+\overline f\,\overline v\). For \(\|\psi\|_{\infty,h}\le1\), subtract the definition of \(\mathcal E_\phi(\nu)\) from the covariant response formula. Cauchy--Schwarz bounds the result by
\[
\bigl(\E\|R\|_2^2\bigr)^{1/2}M^{1/2}
+\tfrac14\E\|W\|_h^2
\le\sqrt{2MD}+\tfrac14\sqrt{32MD}.
\]
Taking the cotangent supremum proves \eqref{eq:deficit-vector}.
\end{proof}

\subsection{Energy measures of near-minimizers}

\begin{theorem}[Complete energy-measure description]\label{thm:energy-measure-description}\label{prop:energy-concentration}
Fix \(u\ne0\), and let \(\phi\) satisfy
\[
\|\phi\|_{\infty,h}=1,\qquad
\operatorname{Re}\langle u,\phi\rangle=\gamma_X(u).
\]
Suppose square-integrable \(V_2\)-valued laws have \(\E q_X(v_j)=u\) and
\(\E\|v_j\|_2^2\to\gamma_X(u)\). Every weak limit \(\nu\) of their energy measures satisfies
\begin{equation}\label{eq:contact-energy-measures}
\nu(X)=\gamma_X(u),\qquad
\supp\nu\subset Z_\phi,
\qquad
u=\int_{Z_\phi}e_X\!\left(p,
\frac{\overline{\varphi(p)}}{|\varphi(p)|}\right)\,d\nu(p),
\end{equation}
where \(Z_\phi=\{|\phi|_h=1\}\).

Conversely, every positive measure satisfying \eqref{eq:contact-energy-measures} is the weak limit of \(|v_j|^2dA\) for smooth vertical directions with
\[
q_X(v_j)=u,\qquad \|v_j\|_2^2\longrightarrow\gamma_X(u).
\]
Thus the same set of limiting measures is obtained from arbitrary centered laws and from symmetric two-point laws.
\end{theorem}

\begin{proof}
For the given laws the deficits tend to zero. Compactness of \(X\) and bounded total mass give weakly convergent subsequences. Equations~\eqref{eq:deficit-spatial} and \eqref{eq:deficit-vector} give the support and moment conditions. Total masses converge to \(\gamma_X(u)\). On \(Z_\phi\), \(\overline f=\overline\varphi/|\varphi|\), which identifies the moment with the atom in \eqref{eq:contact-energy-measures}.

For the converse, choose finite positive atomic measures \(\nu_j\) supported on \(Z_\phi\), of mass \(\gamma_X(u)\), converging weakly to \(\nu\). Their atom moments \(u_j\) converge to \(u\). The set \(Z_\phi\) has empty interior: away from zeros, \(\Delta_h\log|\phi|_h=-2\), so \(|\phi|_h\) cannot equal one on an open set. Reserve a coordinate disk \(O\) whose closure avoids \(Z_\phi\).

Apply \cref{lem:radial-point-response} to the finitely many atoms of \(\nu_j\), merging repeated points and choosing pairwise disjoint disks outside \(\overline O\). Their radii can be taken small enough that the sum \(v_j^0\) has output exactly \(u_j\), energy \(\gamma_X(u)+o(1)\), and energy measure tending to \(\nu\). To make the last choice explicit, take a metric inducing weak convergence on bounded measures and require its distance from \(\nu_j\) to be less than \(1/j\).

Choose once and for all a diagonal response basis supported in \(O\), using \cref{lem:localized-basis,thm:full-range}. Coordinate square roots provide \(w_j\) in its span with
\[
q_X(w_j)=u-u_j,\qquad
\|w_j\|_2^2\le C_X\|u-u_j\|_{\WP}\longrightarrow0.
\]
Disjoint supports make \(v_j=v_j^0+w_j\) an exact realizing direction. Its correction has vanishing energy and hence does not change the limiting energy measure. The symmetric laws on \(\{v_j,-v_j\}\) are centered and have the same measures.
\end{proof}

\subsection{The support radius and the concentration gap}

The maximum sets in \cref{prop:energy-concentration} are constrained by the hyperbolic curvature. This gives a quantitative gap for bounded laws and determines the scale on which sharp realization is possible.

\begin{lemma}[Area near a cotangent maximum]\label{lem:cotangent-maximum-area}
Let \(0\ne\phi=\varphi dz^2\in\QD(X)\), normalized by \(\|\phi\|_{\infty,h}=1\), and put
\[
c_\phi=|\varphi|/\rho^2,
\qquad Z_\phi=\{c_\phi=1\}.
\]
The set \(Z_\phi\) has hyperbolic area zero. There are \(D_X,s_X>0\), uniform over all normalized \(\phi\), such that
\begin{equation}\label{eq:cotangent-maximum-area}
\operatorname{Area}\{1-c_\phi<s\}\le D_X\sqrt{s},
\qquad 0<s\le s_X.
\end{equation}
If every maximum is nondegenerate, the area is at most \(D_\phi s\) for \(0<s\le s_\phi\), with constants depending on \(\phi\).
\end{lemma}

\begin{proof}
Use the Laplace operator \(\Delta_h=\rho^{-2}(\partial_x^2+\partial_y^2)\). Away from the zeros of \(\phi\), holomorphicity and curvature \(-1\) give
\[
\Delta_h\log c_\phi=-2.
\]
At a point of \(Z_\phi\), one has \(c_\phi=1\) and \(dc_\phi=0\). Thus \(f=1-c_\phi\) has a positive semidefinite Hessian with \(\Delta_h f=2\). In a suitably rotated conformal coordinate, \(f_{xx}>0\) at that point. Choose a coordinate rectangle on which \(f_{xx}\ge a>0\). For each fixed \(y\), if \(f(x_1,y),f(x_2,y)<s\), strong convexity and \(f\ge0\) imply
\[
0\le f((x_1+x_2)/2,y)
\le s-\frac a8(x_1-x_2)^2.
\]
The sublevel set on that horizontal line therefore has length at most \(\sqrt{8s/a}\). The area density is bounded on the rectangle. A finite collection of these rectangles covers the compact set \(Z_\phi\), and \(f\) has a positive minimum outside their union. Integration gives the area estimate for a fixed \(\phi\). To make it uniform, consider the compact incidence set
\[
\{(\phi,p):\|\phi\|_{\infty,h}=1,\ c_\phi(p)=1\}.
\]
Near this set, \(c_\phi\) depends smoothly on \(\phi\) and \(p\). The rectangles and positive second-derivative bounds above persist in product neighborhoods. Choose a finite product cover. Compactness gives \(s_X>0\) such that every pair with \(1-c_\phi(p)<s_X\) belongs to that cover. For each \(\phi\), sum the area estimates over its active rectangles. Their finite number and uniform bounds give \eqref{eq:cotangent-maximum-area}. Letting \(s\downarrow0\) proves the zero-area assertion.

If all maxima are nondegenerate, they are isolated and hence finite. Near each maximum, \(f\) bounds a positive multiple of the squared coordinate distance. The corresponding sublevel area is \(O(s)\).
\end{proof}

\begin{proposition}[Bounded-amplitude calibration gap]\label{prop:bounded-amplitude-gap}
There is \(\kappa_X>0\) such that, for every normalized nonzero \(\phi\in\QD(X)\), any bounded vertical law with \(\|v\|_{L^\infty(X\times\Omega)}\le r\), \(r>0\), and \(M=\E\|v\|_2^2\) satisfies
\begin{equation}\label{eq:bounded-amplitude-gap}
M-\operatorname{Re}\langle\E q_X(v),\phi\rangle
\ge\kappa_X\frac{M^3}{r^4}.
\end{equation}
In particular,
\begin{equation}\label{eq:uniform-response-capacity}
\gamma_X(\E q_X(v))\le M-\kappa_X\frac{M^3}{r^4}.
\end{equation}
If every maximum of \(|\phi|_h\) is nondegenerate, then, for some \(\kappa_\phi>0\),
\begin{equation}\label{eq:nondegenerate-amplitude-gap}
M-\operatorname{Re}\langle\E q_X(v),\phi\rangle
\ge\kappa_\phi\frac{M^2}{r^2}.
\end{equation}
Centering is not required for these inequalities.
\end{proposition}

\begin{proof}
\Cref{prop:calibration-deficit}, together with \(1-c_\phi^2\ge1-c_\phi\), gives
\begin{equation}\label{eq:weighted-calibration-gap}
M-\operatorname{Re}\langle\E q_X(v),\phi\rangle
\ge\frac12\int_X(1-c_\phi)a\,dA,
\qquad a=\E|v|^2\le r^2.
\end{equation}
The case \(M=0\) is immediate. Otherwise put \(m=M/r^2\le\operatorname{Area}(X)\). Increase \(D_X\) in \cref{lem:cotangent-maximum-area} so that
\(s=(m/(2D_X))^2\le s_X\) for every such \(m\). At most half of the mass of \(a\,dA\) lies in \(\{1-c_\phi<s\}\), since
\[
\int_{\{1-c_\phi<s\}}a\,dA
\le r^2D_X\sqrt{s}=M/2.
\]
Consequently,
\[
\int_X(1-c_\phi)a\,dA\ge sM/2
=\frac{M^3}{8D_X^2r^4}.
\]
This proves \eqref{eq:bounded-amplitude-gap}. Choose a supporting quadratic differential for \(\E q_X(v)\) to obtain \eqref{eq:uniform-response-capacity}; the zero-output case follows from the same inequality with any normalized \(\phi\). Under the nondegeneracy assumption use the area bound \(D_\phi s\), enlarge \(D_\phi\) so that \(s=m/(2D_\phi)\le s_\phi\), and repeat the argument. It gives \eqref{eq:nondegenerate-amplitude-gap}.
\end{proof}

\subsection{The bounded-amplitude response body}

The norm \(\gamma_X\) optimizes over vertical directions without a prescribed amplitude. Retaining that constraint gives a convex response--energy set. Put
\begin{equation}\label{eq:capped-response-body}
\mathscr C_X=
\overline{\operatorname{conv}}
\{(q_X(w),\|w\|_2^2):w\in V,\ \|w\|_\infty\le1\},
\qquad
\mathscr D_X=\operatorname{pr}_{T_X}\mathscr C_X,
\end{equation}
and define the extended-real function
\begin{equation}\label{eq:capped-value}
J_X(x)=\min\{m:(x,m)\in\mathscr C_X\},\qquad x\in\mathscr D_X,
\end{equation}
with \(J_X=+\infty\) outside \(\mathscr D_X\). The convex hull is over the underlying real space.

\begin{proposition}[Bounded-amplitude value function]\label{prop:capped-value-convex}
The set \(\mathscr C_X\) is compact, and \(\mathscr D_X\) contains a neighborhood of zero. The function \(J_X\) is proper, convex and lower semicontinuous, continuous on \(\operatorname{int}\mathscr D_X\), and satisfies
\[
\gamma_X(x)\le J_X(x)\le\operatorname{Area}(X),
\qquad x\in\mathscr D_X.
\]
Its convex conjugate, using the real part of the tangent--cotangent pairing, is
\begin{equation}\label{eq:capped-conjugate}
J_X^*(\phi)=
\sup_{\substack{w\in V\\\|w\|_\infty\le1}}
\left(\operatorname{Re}\langle q_X(w),\phi\rangle-\|w\|_2^2\right).
\end{equation}
\end{proposition}

\begin{proof}
Calibration gives \(\gamma_X(q_X(w))\le\|w\|_2^2\le\operatorname{Area}(X)\). Thus the set in \eqref{eq:capped-response-body} is bounded in a finite-dimensional space, and its closed convex hull is compact. The coordinate-square response basis gives a neighborhood of zero in its output projection. Taking minima on compact fibers proves lower semicontinuity, while convexity follows by combining points of \(\mathscr C_X\). A finite convex function is continuous on the interior of its domain \cite[Theorem~10.1]{rockafellar1997convex}. The norm lower bound passes to closed convex combinations. Finally, maximizing a linear functional on a closed convex hull gives \eqref{eq:capped-conjugate}.
\end{proof}

The closure includes limiting response--energy pairs, which are controlled by compactness in the finite-dimensional output and energy variables.

\begin{lemma}[Finite mixtures with an interior output]\label{lem:interior-mixture}
For \(x\in\operatorname{int}\mathscr D_X\) and \(\eta>0\), there are directions \(w_1,\ldots,w_N\in V\), with \(\|w_i\|_\infty\le1\), and positive weights \(p_i\) such that
\begin{equation}\label{eq:interior-mixture-data}
\sum_i p_i=1,\qquad
\sum_i p_iq_X(w_i)=x,\qquad
\sum_i p_i\|w_i\|_2^2\le J_X(x)+\eta.
\end{equation}
The vectors \(q_X(w_i)\) affinely span \(T_X\).
\end{lemma}

\begin{proof}
Let \(A=\{q_X(w):\|w\|_\infty\le1\}\). Interior points of \(\overline{\operatorname{conv}}A\) belong to the interior of \(\operatorname{conv}A\) \cite[Theorem~6.3]{rockafellar1997convex}. Choose a finite subset of \(A\) whose convex hull contains a neighborhood of \(x\). On a smaller neighborhood, every point is a convex combination with all these weights positive. Their energies are bounded by \(\operatorname{Area}(X)\).

Choose \(\alpha>0\) so small that adding a proportion \(\alpha\) of such a combination changes the energy by less than \(\eta/2\). By the definition of \(\mathscr C_X\), a finite convex combination has output \(x'\) arbitrarily close to \(x\) and energy at most \(J_X(x)+\eta/2\). Take it close enough that
\[
z=\frac{x-(1-\alpha)x'}\alpha
\]
belongs to the preceding smaller neighborhood. Mix it with a positive combination realizing \(z\), in proportions \(1-\alpha\) and \(\alpha\). The output is exactly \(x\), the energy has the required bound, and the added response columns supply the affine span.
\end{proof}

The exact costs \(\mathcal I_X,\mathcal K_X\) are defined in \cref{sec:exact-laws}. Their admissible laws have centered vertical parameters of supremum norm at most \(r\) and a prescribed exact barycentric output.

\begin{theorem}[Finite-radius creation limit]\label{thm:finite-radius-limit}
For every compact \(K\subset\operatorname{int}\mathscr D_X\),
\begin{equation}\label{eq:finite-radius-limit}
\sup_{x\in K}\left|
\frac{\mathcal I_X(\tau_X^{-1}(r^2x);r)}{r^2}-J_X(x)
\right|\longrightarrow0,
\qquad r\downarrow0.
\end{equation}
The same conclusion holds with \(\mathcal K_X/(4r^2)\) in place of \(\mathcal I_X/r^2\). The recovery laws can be chosen as finite mixtures of symmetric two-point parameter laws.
\end{theorem}

\begin{proof}
Let a centered admissible law have exact output \(r^2x\), and put \(w=v/r\), \(m=\E\|w\|_2^2\). Then \(m\le\operatorname{Area}(X)\),
\[
(\E q_X(w),m)\in\mathscr C_X,
\qquad
x=\E q_X(w)+O_{\WP}(rm),
\]
by \cref{prop:uniform-law-expansion}. If \(r_j\to0\) and \(x_j\to x\in K\), compactness of \(\mathscr C_X\) makes every limiting energy \(m\) satisfy \((x,m)\in\mathscr C_X\), hence \(m\ge J_X(x)\). Continuity of \(J_X\) on \(K\) gives the uniform lower bound.

For recovery, fix \(x\) and \(\eta>0\), and choose the data of \cref{lem:interior-mixture}. For weights \(p\) near those data, give weight \(p_i/2\) to each of \(\Phi(rw_i)\) and \(\Phi(-rw_i)\). Denote the normalized Bers output by \(Y(r,p)\). The pointwise Karcher equation, with its invertible derivative in the barycenter variable at zero, shows that \(Y\) is real analytic in \((r,p)\) with values in \(T_X\). This assertion holds in \(L^\infty\) before projection: the finite collection of bounded coefficient fields remains in a common small ball, and the pointwise analytic series has uniform bounds there.

The law is unchanged by \(r\mapsto-r\). Its expansion is therefore
\[
Y(r,p)=r^2\sum_i p_iq_X(w_i)+O(r^4),
\]
uniformly for nearby weights. Thus \(F(r,p)=Y(r,p)/r^2\) extends real analytically to \(r=0\), where \(F(0,p)=\sum_i p_iq_X(w_i)\). The affine-span condition makes its derivative on \(\sum_i\dot p_i=0\) surjective. Restrict to a real \(2n\)-dimensional right-inverse subspace and apply the implicit-function theorem. It gives positive weights \(p(r,x')\), defined for small \(r\) and \(x'\) near \(x\), whose output is exactly \(r^2x'\). Their parameter laws are centered and have radius at most \(r\).

At \((0,x)\) the rescaled energy is at most \(J_X(x)+\eta\). Continuity of the weights and of \(J_X\) bounds it by \(J_X(x')+2\eta\) on a smaller neighborhood for all sufficiently small \(r\). Cover \(K\) by finitely many such neighborhoods. This proves the uniform upper bound; then let \(\eta\downarrow0\). Finally, \eqref{eq:sharp-dispersion-moment} gives
\[
4\mathcal I_X(t;r)\le\mathcal K_X(t;r)
\le4(1+C_Xr^2)\mathcal I_X(t;r),
\]
which proves the dispersion statement.
\end{proof}

\subsection{The excess above the tangent norm}

Write
\begin{equation}\label{eq:capped-pair-value}
J_X^{(2)}(x)=\inf\{\|w\|_2^2:q_X(w)=x,\ w\in V,\ \|w\|_\infty\le1\}.
\end{equation}
This is the single-direction quadratic cost with an amplitude bound. The convexified value satisfies \(J_X\le J_X^{(2)}\).

\begin{theorem}[Small-output excess at bounded amplitude]\label{thm:capped-excess}
Fix \(\gamma_X(u)=1\). There are \(\lambda_0,C_{X,u}>0\) such that, for \(0<\lambda\le\lambda_0\),
\begin{equation}\label{eq:capped-excess-general}
\lambda+\kappa_X\lambda^3
\le J_X(\lambda u)
\le J_X^{(2)}(\lambda u)
\le\lambda+C_{X,u}\lambda^2.
\end{equation}
The upper bound is realized by a smooth direction with supremum norm at most \(1/2\), supported on disjoint disks with a fixed nonempty open complement.

If a normalized supporting differential \(\phi\) for \(u\) has only nondegenerate maxima of \(|\phi|_h\), then, for some \(c_{X,u,\phi}>0\),
\begin{equation}\label{eq:capped-excess-nondegenerate}
c_{X,u,\phi}\lambda^2
\le J_X(\lambda u)-\lambda
\le J_X^{(2)}(\lambda u)-\lambda
\le C_{X,u}\lambda^2.
\end{equation}
\end{theorem}

\begin{proof}
Take a fixed optimal atomic representation
\[
u=\sum_{j=1}^Np_j e_X(p_j',a_j),\qquad
p_j>0,\quad\sum_jp_j=1.
\]
As in \cref{thm:single-direction-cost}, merge repeated points and choose disjoint coordinate disks about the remaining points, leaving a fixed open disk \(O\) outside them. Apply \cref{lem:radial-point-response} with masses \(m_j=\lambda p_j\) and radii
\(\varepsilon_j=L\sqrt{\lambda p_j}\). Choose \(L\) large enough that \eqref{eq:radial-amplitude} bounds every amplitude by \(1/2\), and then take \(\lambda_0\) small enough that all supporting disks remain in the chosen coordinate disks. Their sum \(w_\lambda\) has exact output \(\lambda u\), and
\[
\|w_\lambda\|_2^2
\le\sum_j\lambda p_j(1+C_jL^2\lambda p_j)
\le\lambda+C_{X,u}\lambda^2.
\]

For any finite mixture of directions of amplitude at most one, \cref{prop:bounded-amplitude-gap} gives, against a supporting \(\phi\),
\[
m-\operatorname{Re}\langle x,\phi\rangle\ge\kappa_Xm^3,
\]
and under the nondegeneracy hypothesis it gives the bound \(\kappa_\phi m^2\). These inequalities are preserved under limits of response--energy pairs. Apply them to \((\lambda u,J_X(\lambda u))\in\mathscr C_X\) and use \(J_X(\lambda u)\ge\lambda\). The same lower bounds hold for the single-direction infimum.
\end{proof}

\subsection{Nondegenerate supporting differentials}

Away from its zero divisor, the norm of a holomorphic section of a Hermitian line bundle is critical precisely where its Chern covariant derivative vanishes. Douglas--Shiffman--Zelditch develop the corresponding jet and Hessian description \cite[\S\S2,4]{douglas2004critical}. Applied to \(K_X^2\) with its hyperbolic Hermitian metric, this framework gives the nondegeneracy condition used in the quadratic excess law. We verify its genericity on every closed surface, treating the genus-two Weierstrass points separately, and identify the resulting contact faces.

\begin{theorem}[Nondegenerate cotangent maxima]\label{thm:nondegenerate-cotangents}
For every closed hyperbolic surface \(X\) of genus at least two, the set
\[
\mathcal G_X=\{0\ne\phi\in\QD(X):\text{every global maximum of }|\phi|_h
\text{ is nondegenerate}\}
\]
is open and dense in \(\QD(X)\setminus\{0\}\).
For a normalized \(\phi\in\mathcal G_X\), the maximum set \(Z_\phi\) is finite and the exposed face
\begin{equation}\label{eq:finite-contact-face}
\{u:\gamma_X(u)=\operatorname{Re}\langle u,\phi\rangle=1\}
=\operatorname{conv}\{e_X(p,a_\phi(p)):p\in Z_\phi\},
\qquad a_\phi=\overline\varphi/|\varphi|.
\end{equation}
Every direction in the positive cone over this face has a supporting differential with nondegenerate maxima.

In a disk uniformization centered at a normalized maximum \(p\), write
\(\phi=\varphi(z)dz^2\) and \(\rho(z)=2/(1-|z|^2)\). Then \(\varphi'(0)=0\), and the eigenvalues of the hyperbolic Hessian of \(1-|\phi|_h\) at \(p\) are
\begin{equation}\label{eq:maximum-jet-criterion}
1\pm\frac14\left|\frac{\varphi''(0)}{\varphi(0)}\right|.
\end{equation}
In particular, this maximum is nondegenerate precisely when
\(|\varphi''(0)|<4|\varphi(0)|\).
\end{theorem}

\begin{proof}
We use Serre duality and Riemann--Roch in the form recalled in \cite[Appendices~A.9--A.11]{hubbard2016teichmuller}. First suppose \(g\ge3\). The map from \(H^0(X,K_X^2)\) to its first jets at any \(p\) is onto. Indeed, the obstruction group is
\[
H^1(X,K_X^2(-2p))\simeq H^0(X,K_X^{-1}(2p))^*=0,
\]
since the latter line bundle has degree \(4-2g<0\).
The second-jet evaluation is also onto: its obstruction is
\(H^1(X,K_X^2(-3p))\simeq H^0(X,K_X^{-1}(3p))^*=0\),
since \(\deg K_X^{-1}(3p)=5-2g<0\). Thus \(K_X^2\) satisfies the jet-spanning hypothesis of \cite[Definition~4.1]{douglas2004critical} in this genus range.
On the open set of pairs \((\phi,p)\) with \(\phi(p)\ne0\), consider the section
\(F(\phi,p)=d_p\log|\phi|_h\) of the pulled-back real cotangent bundle of \(X\). Its variation in \(\phi\) along \(\psi\) is
\[
D_\phi F[\psi]=\operatorname{Re}d(\psi/\phi).
\]
First-jet surjectivity lets us set \(\psi(p)=0\) and prescribe its derivative arbitrarily. Hence \(F\) is transverse to the zero section. Parametric transversality, applied on a countable collection of product neighborhoods in this open set, gives a dense residual set of differentials for which every critical point away from the zero divisor is nondegenerate \cite[Ch.~3, Theorem~2.7]{hirsch1976differential}. Global maxima have positive value, so this proves density of \(\mathcal G_X\).

For \(g=2\), the same argument applies except at the six Weierstrass points: the degree-zero bundle \(K_X^{-1}(2p)\) has a section exactly when \(2p\) is canonical. At such a point the hyperelliptic involution has local form \(z\mapsto-z\), preserves the hyperbolic metric, and acts trivially on \(H^0(X,K_X^2)\). The last assertion follows from
\(H^0(X,K_X^2)=\operatorname{Sym}^2H^0(X,K_X)\), both spaces having dimension three. Thus every differential nonzero at that point is critical there. The system \(K_X^2\) has no base points, so choose \(\phi_0(p)\ne0\). Riemann--Roch gives
\[
h^0(K_X^2(-2p))=2,\qquad h^0(K_X^2(-3p))=1.
\]
There is consequently a \(\psi\) vanishing to order exactly two at \(p\). In the family \(\phi_0+t\psi\), \(t\in\mathbb C\), the value at \(p\) is fixed and the complex second derivative of \(\log\phi\) varies freely. The trace of the real Hessian of \(\log|\phi|_h\) is fixed by curvature, while its trace-free part ranges over all trace-free symmetric matrices. Its determinant is therefore a nonzero real-analytic function of the differential on \(\{\phi(p)\ne0\}\). The conditions \(\phi(p)=0\) and degeneracy at these finitely many points can both be avoided on an open dense set. Intersecting with the residual set from transversality on their complement proves density also in genus two.

Nondegenerate maxima are isolated, so \(Z_\phi\) is finite. Choose small disjoint neighborhoods of them on which the Hessian is negative definite. The norm has a strict gap below its maximum outside these neighborhoods. Both properties persist under small perturbations of \(\phi\); the implicit-function theorem continues the unique critical maximum in each neighborhood. Every global maximum of a nearby differential lies among them. This proves openness.

For \eqref{eq:maximum-jet-criterion}, expand at a normalized maximum:
\[
\log|\phi|_h
=\frac12\operatorname{Re}\!\left(\frac{\varphi''(0)}{\varphi(0)}z^2\right)
-2|z|^2+O(|z|^3).
\]
The hyperbolic Hessian is one quarter of the Euclidean Hessian at zero. Negating the displayed quadratic term gives \eqref{eq:maximum-jet-criterion}. Finally, equality against \(\phi\) in the atomic unit ball \eqref{eq:geometric-response-body} forces each atom to lie at a maximum with phase \(a_\phi\). This proves \eqref{eq:finite-contact-face} and the assertion about its positive cone.
\end{proof}

\subsection{A curvature bound for the excess coefficient}

Let \(\phi\in\mathcal G_X\) be normalized and put
\[
g_\phi=\tfrac12(1-|\phi|_h^2),\qquad
H_p=\operatorname{Hess}_h g_\phi(p),\qquad
A_\phi=2\pi\sum_{p\in Z_\phi}(\det H_p)^{-1/2}.
\]
The determinants use hyperbolic orthonormal frames. At a maximum, \(H_p\) equals the Hessian in \eqref{eq:maximum-jet-criterion}.
Define
\begin{equation}\label{eq:bathtub-deficit-function}
B_\phi(m)=\int_0^{1/2}
\bigl(m-\operatorname{Area}\{g_\phi\le s\}\bigr)_+\,ds,
\qquad 0\le m\le\operatorname{Area}(X).
\end{equation}

\begin{proposition}[Explicit lower excess coefficient]\label{prop:explicit-excess-coefficient}
The function \(B_\phi\) is continuous and convex, and
\begin{equation}\label{eq:bathtub-small-mass}
B_\phi(m)=\frac{m^2}{2A_\phi}+o(m^2),\qquad m\downarrow0.
\end{equation}
For every \(x\in\mathscr D_X\),
\begin{equation}\label{eq:capped-bathtub-lower}
J_X(x)-\operatorname{Re}\langle x,\phi\rangle
\ge B_\phi(J_X(x)).
\end{equation}
If \(\gamma_X(u)=\operatorname{Re}\langle u,\phi\rangle=1\), then
\begin{equation}\label{eq:explicit-excess-liminf}
\liminf_{\lambda\downarrow0}
\frac{J_X(\lambda u)-\lambda}{\lambda^2}
\ge\frac{1}{2A_\phi}
=\frac{1}{4\pi\displaystyle\sum_{p\in Z_\phi}(\det H_p)^{-1/2}}.
\end{equation}
The same lower bound holds with \(J_X^{(2)}\) in place of \(J_X\).
\end{proposition}

\begin{proof}
The integral in \eqref{eq:bathtub-deficit-function} is a continuous convex function of \(m\). For a law of amplitude at most one, set \(a=\E|w|^2\), \(m=\int_Xa\,dA\), and \(x=\E q_X(w)\). The deficit identity gives
\[
m-\operatorname{Re}\langle x,\phi\rangle\ge\int_Xg_\phi a\,dA.
\]
Since \(0\le a\le1\), the layer-cake formula yields
\[
\int_Xg_\phi a\,dA
=\int_0^{1/2}\int_{\{g_\phi>s\}}a\,dA\,ds
\ge\int_0^{1/2}\bigl(m-\operatorname{Area}\{g_\phi\le s\}\bigr)_+\,ds.
\]
This inequality holds for every finite mixture. Continuity of \(B_\phi\) passes it to the closed response--energy set \(\mathscr C_X\); taking its minimum-energy point proves \eqref{eq:capped-bathtub-lower}.

In hyperbolic normal coordinates at each \(p\in Z_\phi\),
\[
g_\phi(\xi)=\tfrac12\langle H_p\xi,\xi\rangle+O(|\xi|^3),
\qquad dA=(1+O(|\xi|^2))\,d\xi.
\]
Elliptic rescaling of these finitely many minima gives
\(\operatorname{Area}\{g_\phi\le s\}=A_\phi s+o(s)\).
For each \(\varepsilon>0\), bound this sublevel area between
\((A_\phi-\varepsilon)s\) and \((A_\phi+\varepsilon)s\) for small \(s\). Integration in \eqref{eq:bathtub-deficit-function}, followed by \(\varepsilon\downarrow0\), proves \eqref{eq:bathtub-small-mass}. By \cref{thm:capped-excess},
\(J_X(\lambda u)=\lambda+O(\lambda^2)\).
Substitution in \eqref{eq:capped-bathtub-lower} proves \eqref{eq:explicit-excess-liminf}. The inequality \(J_X^{(2)}\ge J_X\) gives the final assertion.
\end{proof}

The coefficient in \eqref{eq:explicit-excess-liminf} is a lower bound determined by the local curvature of the contact set. The remaining terms in the deficit identity and the vector output constraint also enter the full excess minimization.

\section{Exact realization and the two costs}\label{sec:constrained}

We pass from the quadratic response to exact barycentric output. The local square form supplies a nonlinear correction on disjoint support, and the energy-weighted remainder preserves the leading creation cost. Restoration is then estimated directly from the exact output over the full local fiber graph.

\subsection{Pointwise aggregation}\label{sec:probabilistic}

\paragraph{Random Beltrami fields}

A random Beltrami field is a jointly measurable map \(\mu:X\times\Omega\to\D\), with \(\mu_\omega\in\Bel(X)\) almost surely, on a probability space \((\Omega,\mathscr F,\Prob)\). We assume the uniform ellipticity bound
\begin{equation}\label{eq:uniform-kappa}
  \|\mu\|_{L^\infty(X\times\Omega)}\le \kappa<1.
\end{equation}
Write \(\Pi\) for this probability model. For vector-valued expectations, the vertical parameter \(\xi:\Omega\to V\) is strongly measurable and essentially bounded, so its Bochner expectation belongs to the closed subspace \(V\).

\paragraph{Hyperbolic coefficient barycenter}

The Poincar\'e distance on \(\D\) is
\begin{equation}\label{eq:poincare-distance}
  \distH(a,b)
  =2\operatorname{arctanh}
  \left|
  \frac{a-b}{1-\bar a b}
  \right|.
\end{equation}
Since \((\D,\distH)\) is a Hadamard manifold, every probability measure with finite second moment admits a unique Fr\'echet barycenter \cite[Proposition~4.3]{sturm2003probability}.

\begin{definition}[Ambient hyperbolic representative]\label{def:ambient-barycenter}
For a random Beltrami field satisfying \eqref{eq:uniform-kappa}, define
\begin{equation}\label{eq:beta-def}
  \beta_\Pi(z)
  :=\arg\min_{b\in\D}
  \E\bigl[\distH(b,\mu(z,\cdot))^2\bigr].
\end{equation}
We call \(\beta_\Pi\) the ambient hyperbolic Fr\'echet representative of \(\Pi\).
\end{definition}

\begin{proposition}[Measurability and admissibility]\label{prop:beta-measurable}
Under \eqref{eq:uniform-kappa}, the map \(z\mapsto\beta_\Pi(z)\) is measurable and satisfies
\[
  \|\beta_\Pi\|_\infty\le\kappa<1.
\]
Hence \(\beta_\Pi\in\Bel(X)\).
\end{proposition}

\begin{proof}
Put \(K_\kappa=\{|w|\le\kappa\}\). The pointwise laws of \(\mu(z,\cdot)\) form a measurable map into the probability measures on the compact hyperbolic disk \(K_\kappa\). On this space the \(L^1\)-Wasserstein metric induces weak convergence. The barycenter map is continuous for this metric \cite[Theorem~6.3]{sturm2003probability}, so its composition with the law map is measurable. Since \(K_\kappa\) is hyperbolically convex and contains every pointwise support, it also contains the barycenter \cite[Proposition~6.1]{sturm2003probability}.
\end{proof}

Conformal changes of source coordinate act on coefficients by rotations of \(\D\), which are Poincar\'e isometries. Uniqueness therefore makes \(\beta_\Pi\) a global Beltrami differential. Aggregation uses the affine coefficient model and the hyperbolic area of \(X\).

\paragraph{Integrated Fr\'echet geometry}

Define
\begin{equation}\label{eq:D2-field}
  \distB(\mu,\nu)^2
  :=\int_X
  \distH(\mu(z),\nu(z))^2
  \,dA_X(z).
\end{equation}
Fubini's theorem identifies \(\E\distB(\nu,\mu)^2\) with the integral of the pointwise losses. Hence \(\beta_\Pi\) minimizes the integrated functional:
\begin{equation}\label{eq:ambient-functional}
  \Energy_\Pi(\nu)
  :=\E\distB(\nu,\mu)^2,
  \qquad
  \Energy_\Pi(\beta_\Pi)
  =\inf_\nu\Energy_\Pi(\nu).
\end{equation}

The associated dispersion is
\begin{equation}\label{eq:dispersion}
  \cU_2(\Pi)^2
  :=\Energy_\Pi(\beta_\Pi).
\end{equation}
Its infimum at a prescribed output is the intrinsic creation cost.

\begin{lemma}[Uniform second-moment bounds]\label{lem:barycenter-moment}
There exist \(\rho_0>0\) and \(C<\infty\) such that, for every random coefficient field with \(\|\mu\|_\infty\le\rho\le\rho_0\), its pointwise barycenter \(\beta\) satisfies
\begin{align}
\|\beta\|_2^2&\le C\E\|\mu\|_2^2,\label{eq:beta-second-moment}\\
\|\beta-\E\mu\|_1&\le C\rho\E\|\mu\|_2^2.\label{eq:beta-affine-L1}
\end{align}
Here \(\|\cdot\|_1\) is the hyperbolic \(L^1(X)\) norm.
\end{lemma}

\begin{proof}
The \(\operatorname{CAT}(0)\) variance inequality \cite[Proposition~4.4]{sturm2003probability}, applied pointwise with comparison point zero, gives
\[
\distH(0,\beta)^2\le\E\distH(0,\mu)^2.
\]
On the disk of radius \(\rho_0\), \(2|a|\le\distH(0,a)\le C|a|\). Thus \(|\beta|^2\le C\E|\mu|^2\), proving \eqref{eq:beta-second-moment} after integration.

Put \(\phi_b(a)=(a-b)/(1-\bar b a)\). Uniformly for \(|a|,|b|\le\rho\),
\[
\frac{\operatorname{arctanh}|\phi_b(a)|}{|\phi_b(a)|}\phi_b(a)
=(a-b)+O\bigl(\rho(|a|^2+|b|^2)\bigr).
\]
The Karcher equation \cite{karcher1977riemannian} and the pointwise second-moment bound therefore give
\[
|\beta-\E\mu|\le C\rho\E|\mu|^2.
\]
Integration proves \eqref{eq:beta-affine-L1}.
\end{proof}

\paragraph{Aggregation at the conformal basepoint}\label{sec:ambient}

At zero, the Poincar\'e barycenter agrees with the affine expectation through second order.

\begin{lemma}[Second-order aggregation at the conformal basepoint]\label{lem:second-order-input}
Suppose
\[
  Y_\varepsilon
  =\varepsilon X_1+\varepsilon^2X_2+R_\varepsilon,
\]
where \(X_1,X_2\) are uniformly bounded random fields and
\[
  \|R_\varepsilon\|_{L^\infty(X\times\Omega)}
  \le C|\varepsilon|^3.
\]
If \(b_\varepsilon\) is the pointwise hyperbolic barycenter of \(Y_\varepsilon\), then
\begin{equation}\label{eq:general-zero-barycenter}
  b_\varepsilon
  =\varepsilon\E X_1
  +\varepsilon^2\E X_2
  +O_{L^\infty}(\varepsilon^3).
\end{equation}
\end{lemma}

\begin{proof}
The assumed uniform bounds give \(\|Y_\varepsilon\|_\infty=O(|\varepsilon|)\). Pointwise, the proof of \cref{lem:barycenter-moment} gives
\[
|b_\varepsilon-\E Y_\varepsilon|
\le C\|Y_\varepsilon\|_\infty\E|Y_\varepsilon|^2
=O(|\varepsilon|^3),
\]
uniformly on \(X\). Taking expectations in the assumed expansion proves \eqref{eq:general-zero-barycenter}.
\end{proof}

\begin{proposition}[Quadratic output of pointwise aggregation]\label{prop:aggregation-two-jet}
Let \(\xi\) be a bounded strongly measurable \(V\)-valued law, and put \(m=\E\xi\). Any Beltrami-valued pointwise aggregation rule satisfying
\[
\mathcal M(Y_\varepsilon)=\E Y_\varepsilon+O_{L^\infty}(\varepsilon^3)
\]
for uniformly bounded expansions
\(Y_\varepsilon=\varepsilon X_1+\varepsilon^2X_2+O_{L^\infty}(\varepsilon^3)\)
has, on the fiber samples \(\Phi(\varepsilon\xi)\), the output
\[
\tau\bigl(\pi(\mathcal M(\Phi(\varepsilon\xi)))\bigr)
=\varepsilon^2D\tau_{[X]}\!\left[\E q_X(\xi-m)\right]+O(\varepsilon^3)
\]
in any chart \(\tau\) with \(\tau([X])=0\). Both affine averaging and the Poincar\'e barycenter satisfy this hypothesis.
\end{proposition}
\begin{proof}
The graph expansion gives
\(\mathcal M(\Phi(\varepsilon\xi))=\varepsilon m+\tfrac12\varepsilon^2\E B(\xi,\xi)+O(\varepsilon^3)\).
Taylor expansion of \(\tau\circ\pi\), followed by \eqref{eq:B-S}, leaves
\(\tfrac12D\tau[\vH_0(m,m)-\E\vH_0(\xi,\xi)]\), which equals the claimed centered response. Lemma~\ref{lem:second-order-input} verifies the Poincar\'e case.
\end{proof}

\subsection{Bers coordinates and the exact output}

Fix a cocompact uniformization \(X\simeq\D/\Gamma\), and write \(\D^*=\widehat{\mathbb C}\setminus\overline\D\). Let \(\mathfrak P=\mathfrak b\circ\pi\) be the associated Bers projection, normalized by \(\mathfrak P(0)=0\), into the space of \(\Gamma^*\)-automorphic holomorphic quadratic differentials on \(\D^*\). Here \(\Gamma^*\) is the reflected action on \(\D^*\). We use the Bers norm
\[
  \|\varphi\|_B
  :=\sup_{z\in\D^*}(|z|^2-1)^2|\varphi(z)|.
\]
The norm at infinity is interpreted in the spherical coordinate. The \(L^p\) bounds for the Beurling transform and the Neumann-series solution of the Beltrami equation are classical \cite[Secs.~4.5.2, 5.3, and 5.7]{astala2008elliptic}. We retain the hyperbolic \(L^2\) factors of the compact quotient in the following estimates.

Let \(\mathfrak b\) be the Bers embedding fixed above. Since its target automorphic space has complex dimension \(n=3g-3\), the derivative \(D\mathfrak b_{[X]}\) is a linear isomorphism. Define the normalized chart
\begin{equation}\label{eq:normalized-Bers-chart}
\tau_X(t):=(D\mathfrak b_{[X]})^{-1}[\mathfrak b(t)]
\in T_{[X]}\Teich(X).
\end{equation}
Then \(\tau_X([X])=0\), \((D\tau_X)_{[X]}=\id\), and, on a fixed neighborhood of \([X]\),
\begin{equation}\label{eq:WP-chart-distance}
\|\tau_X(t)\|_{\WP}\asymp_X d_{\WP}([X],t).
\end{equation}
For \(v\) and \(-v\) in the graph domain, set
\begin{equation}\label{eq:two-point-midpoint-map}
\Pi_v:=\tfrac12\delta_{\Phi(v)}+\tfrac12\delta_{\Phi(-v)},
\qquad
\beta_v:=\barH(\Pi_v),
\qquad
G_X(v):=\tau_X(\pi(\beta_v)).
\end{equation}
Thus both samples have class \([X]\), while \(G_X(v)\) records the class produced by their pointwise midpoint.

\subsection{Centered laws and exact creation costs}\label{sec:exact-laws}

Choose a sufficiently small \(r_X>0\) so that \(\Phi\{v:\|v\|_\infty\le r_X\}\) stays in the coefficient balls used in \cref{lem:barycenter-moment,prop:bers-nonlinear-l2}. For \(0<r\le r_X\), consider all strongly measurable \(V\)-valued random parameters satisfying
\begin{equation}\label{eq:centered-law-class}
\E v=0,\qquad \|v\|_{L^\infty(X\times\Omega)}\le r,
\qquad \mu=\Phi(v).
\end{equation}
For their coefficient laws \(\Pi\), define
\begin{align}
\mathcal I_X(t;r)
&:=\inf\{\E\|v\|_2^2:\eqref{eq:centered-law-class},\ \pi(\beta_\Pi)=t\},\label{eq:exact-vertical-cost}\\
\mathcal K_X(t;r)
&:=\inf\{\cU_2(\Pi)^2:\eqref{eq:centered-law-class},\ \pi(\beta_\Pi)=t\}.\label{eq:exact-intrinsic-cost}
\end{align}
The infima range over arbitrary probability spaces and centered laws with the stated bound.

\subsection{Exact local reachability}

\begin{theorem}[Exact local reachability]\label{thm:exact-reach}
Let \(E_0\) be the subspace in \cref{thm:full-range}. On a fixed neighborhood of zero in \(V\), the map \(G_X\) is even and real analytic, with
\begin{equation}\label{eq:midpoint-quartic-normal-form}
G_X(v)=-Q_0(v)+O(\|v\|_\infty^4).
\end{equation}
Restricted to \(E_0\), zero is an isolated zero and its local Brouwer degree, with the complex orientations on domain and target, is \(2^n\).

There exist \(r_X,a_X,C_X>0\) such that, for every \(0<r\le r_X\),
\begin{equation}\label{eq:exact-reach-ball}
\{y:\|y\|_{\WP}\le a_Xr^2\}
\subset G_X\{v\in E_0:\|v\|_\infty<r\}.
\end{equation}
For every such \(y\), a realizing parameter can be chosen with
\begin{equation}\label{eq:exact-reach-size}
G_X(v)=y,
\qquad
\|v\|_\infty\le C_X\|y\|_{\WP}^{1/2},
\qquad
\|v\|_2^2\le C_X\|y\|_{\WP}.
\end{equation}
\end{theorem}

\begin{proof}
The midpoint of two points of the Poincar\'e disk is real analytic near \((0,0)\). Its convergent power series acts uniformly on a fixed \(L^\infty\) ball, so composition with the analytic graph and projection makes \(G_X\) real analytic. Interchanging the two samples gives \(G_X(-v)=G_X(v)\). By \cref{lem:second-order-input,prop:B-S},
\[
\beta_v=\tfrac12B(v,v)+O_{L^\infty}(\|v\|_\infty^3),
\qquad
G_X(v)=-Q_0(v)+O(\|v\|_\infty^3).
\]
Evenness removes the cubic terms in both real Taylor expansions, yielding \eqref{eq:midpoint-quartic-normal-form} and
\begin{equation}\label{eq:symmetric-midpoint-H}
\beta_v=\tfrac12B(v,v)+O_{L^\infty}(\|v\|_\infty^4).
\end{equation}

Write \(v=\sum_i c_iv_i\). Equation~\eqref{eq:diagonal-additivity} gives
\[
q(v):=-Q_0(v)=-\sum_i c_i^2Q_i.
\]
This is the coordinate-square map followed by a complex linear isomorphism. It has only the zero \(v=0\), and its local real degree is \(2^n\): a sufficiently small regular value with all coordinates nonzero has \(2^n\) preimages, each of positive real orientation. The finite-dimensional degree properties used here are recalled in \cite[Ch.~1]{deimling2013nonlinear}; see also \cite{arutyunov2016properties} for even maps with quadratic leading term.

Compactness of the unit sphere of \(E_0\) gives \(a>0\) with
\[
\|q(v)\|_{\WP}\ge a\|v\|_\infty^2.
\]
Choose \(r_X\) so that the remainder in \eqref{eq:midpoint-quartic-normal-form} is at most \((a/4)\|v\|_\infty^2\) for \(\|v\|_\infty\le r_X\). Then zero remains isolated and, for \(\|y\|_{\WP}\le ar^2/4\), the homotopy
\[
q+s(G_X-q)-y,\qquad 0\le s\le1,
\]
has no zero on the boundary of \(D_r=\{v\in E_0:\|v\|_\infty<r\}\). Homotopy invariance and variation of the target value inside this ball give
\[
\deg(G_X-y,D_r,0)=\deg(q,D_r,0)=2^n.
\]
Nonzero degree gives a preimage and proves \eqref{eq:exact-reach-ball} with \(a_X=a/4\). Moreover,
\[
\|G_X(v)\|_{\WP}\ge\tfrac{3a}{4}\|v\|_\infty^2
\]
throughout this neighborhood. Every realizing parameter in \(D_r\) therefore satisfies the first size bound in \eqref{eq:exact-reach-size}. Norm equivalence on \(E_0\) gives the second.
\end{proof}

\subsection{Uniform second-order output of centered laws}

The next estimate keeps the second moment in the remainder. Its constants are independent of the number of samples, their probabilities, and their spatial concentration.

\begin{proposition}[Uniform law expansion]\label{prop:uniform-law-expansion}
After reducing \(r_X\), every centered law in \eqref{eq:centered-law-class} with \(0<r\le r_X\) and \(M=\E\|v\|_2^2\) satisfies
\begin{equation}\label{eq:uniform-law-output}
\tau_X(\pi(\beta_\Pi))
=-\E Q_0(v)+R_\Pi,
\qquad \|R_\Pi\|_{\WP}\le C_XrM.
\end{equation}
Its integrated dispersion satisfies
\begin{equation}\label{eq:sharp-dispersion-moment}
4M\le\cU_2(\Pi)^2\le4(1+C_Xr^2)M.
\end{equation}
\end{proposition}

\begin{proof}
Write \(\mu=\Phi(v)\), \(\beta=\beta_\Pi\), and \(\rho=\|\mu\|_{L^\infty(X\times\Omega)}\le C_Xr\). The cubic fiber estimate gives
\[
\E\mu=\E\Psi(v)=\tfrac12\E B(v,v)+e,
\qquad \|e\|_2\le C_XrM,
\qquad \|\E\mu\|_2\le C_XM.
\]
Orthogonality and the graph estimate give
\begin{equation}\label{eq:law-moment-graph}
\E\|\mu\|_2^2=M+\E\|\Psi(v)\|_2^2\le C_XM,
\qquad \|\E\Psi(v)\|_2\le C_XM.
\end{equation}
The pointwise estimate in the proof of \cref{lem:barycenter-moment} is
\[
|\beta-\E\mu|\le C\rho\,\E|\mu|^2.
\]
Together with \eqref{eq:law-moment-graph}, it yields
\begin{align*}
\|\beta-\E\mu\|_1&\le C_XrM,\\
\|\beta-\E\mu\|_2^2
&\le C\rho^2\int_X(\E|\mu|^2)^2\,dA_X
\le C\rho^4\E\|\mu\|_2^2\le C_Xr^4M.
\end{align*}
Since \(M\le\operatorname{Area}(X)r^2\),
\[
\|\beta\|_2^2\le C_X(M^2+r^4M)\le C_Xr^2M.
\]
The harmonic projection has smooth finite-dimensional range, so it is bounded from \(L^1(X)\) to \(H\). Since the linear Bers projection factors through it,
\begin{equation}\label{eq:linear-Bers-L1}
\|D\mathfrak P_0[\eta]\|_B\le C_X\|\eta\|_1.
\end{equation}
Apply this to the mean expansion, then use \eqref{eq:bers-nonlinear-zero} at \(\beta\). The identity
\[
(D\mathfrak b_{[X]})^{-1}D\mathfrak P_0[\tfrac12B(v,v)]=-Q_0(v)
\]
proves \eqref{eq:uniform-law-output}.

The lower bound in \eqref{eq:sharp-dispersion-moment} follows from \(\distH(a,b)\ge2|a-b|\), the Euclidean variance identity, and \(V\perp H\). For the upper bound, compare the minimizing coefficient \(\beta\) with the affine mean \(\E\mu\). The straight segment between two points of the disk of radius \(\rho\) has Poincar\'e length at most \(2|a-b|/(1-\rho^2)\). Thus
\[
\cU_2(\Pi)^2
\le4(1+C\rho^2)\E\|\mu-\E\mu\|_2^2.
\]
Orthogonality gives
\[
\E\|\mu-\E\mu\|_2^2
=M+\E\|\Psi(v)-\E\Psi(v)\|_2^2
\le M+C_X\E\|v\|_2^4
\le(1+C_Xr^2)M.
\]
This proves the upper bound.
\end{proof}

\subsection{Sharp cost of exact realization}

\begin{theorem}[Sharp cost of exact realization]\label{thm:sharp-exact-cost}
Let \(u\in T_X\), \(0<r_j\le r_X\), \(r_j\to0\), and \(\delta_j>0\) satisfy
\begin{equation}\label{eq:sharp-cost-regime}
\delta_j/r_j^2\longrightarrow0.
\end{equation}
For \(t_j=\tau_X^{-1}(\delta_j u)\),
\begin{equation}\label{eq:sharp-exact-cost}
\frac{\mathcal I_X(t_j;r_j)}{\delta_j}\longrightarrow\gamma_X(u),
\qquad
\frac{\mathcal K_X(t_j;r_j)}{4\delta_j}\longrightarrow\gamma_X(u).
\end{equation}
Both limits are unchanged when the admissible coefficient laws are restricted to the equal-weight pairs \(\Pi_v\) in \eqref{eq:two-point-midpoint-map}. If all vertical samples are required to lie in \(E_0\), the limits equal \(\gamma_{X,E_0}(u)\) in \eqref{eq:localized-sharp-cost}, again with two-point laws.
\end{theorem}

\begin{proof}
The zero law treats \(u=0\). For a centered law with support radius \(r\) and energy \(M\), \cref{prop:uniform-law-expansion} gives
\[
y:=\tau_X(\pi(\beta_\Pi))=\E q_X(v)+R,
\qquad \|R\|_{\WP}\le C_XrM.
\]
The norm bound \(\gamma_X(\E q_X(v))\le M\) implies
\[
\gamma_X(y)\le(1+C_Xr)M,
\qquad
\mathcal I_X(t_j;r_j)\ge\frac{\delta_j\gamma_X(u)}{1+C_Xr_j}.
\]
The dispersion comparison gives \(\mathcal K_X\ge4\mathcal I_X\).

For the upper bound, fix \(\eta>0\). By \cref{thm:single-direction-cost}, choose a smooth compactly supported displacement field whose vertical direction \(v_0\) satisfies
\[
q_X(v_0)=u,
\qquad \|v_0\|_2^2\le\gamma_X(u)+\eta,
\]
and whose support leaves a nonempty open complement. Choose a localized diagonal response subspace \(E\) in that complement. Disjoint supports give
\[
q_X(\sqrt\delta\,v_0+w)=\delta u+q_X(w),
\qquad v_0\perp_{L^2}E.
\]
The square form on \(E\) gives constants \(a,A>0\) with
\[
\|q_X(w)\|_{\WP}\ge a\|w\|_\infty^2,
\qquad \|w\|_2^2\le A\|w\|_\infty^2.
\]
By \eqref{eq:midpoint-quartic-normal-form}, the remainder in
\[
G_X(\sqrt\delta\,v_0+w)-\delta u=q_X(w)+R_\delta(w)
\]
is bounded by \(C\|\sqrt\delta\,v_0+w\|_\infty^4\). Set
\(C_0=C(\|v_0\|_\infty+1)^4\), choose \(K\) with \(aK^2>C_0\), and then take \(\delta\) small enough that \(K\sqrt\delta\le1\). On the closed ball \(\|w\|_\infty\le K\delta\),
\[
\|R_\delta(w)\|_{\WP}\le C_0\delta^2.
\]
The quadratic term strictly dominates this bound on the boundary. The homotopy \(q_X(w)+sR_\delta(w)\), \(0\le s\le1\), preserves the coordinate-square degree \(2^n\). Thus some \(w_\delta\in E\) satisfies
\[
G_X(\sqrt\delta\,v_0+w_\delta)=\delta u,
\qquad \|w_\delta\|_\infty\le K\delta.
\]
Set \(v_j=\sqrt{\delta_j}\,v_0+w_{\delta_j}\). Condition~\eqref{eq:sharp-cost-regime} gives \(\|v_j\|_\infty\le r_j\) for large \(j\), and orthogonality gives
\[
\|v_j\|_2^2\le\delta_j(\gamma_X(u)+\eta)+AK^2\delta_j^2.
\]
The law \(\Pi_{v_j}\) realizes \(t_j\) exactly. Its dispersion is at most \(4(1+C_Xr_j^2)\|v_j\|_2^2\). Let \(j\to\infty\), then \(\eta\downarrow0\), to obtain both upper bounds with two samples.

For laws in \(E_0\), the lower bound uses \(\gamma_{X,E_0}\). The parameters supplied by \cref{thm:exact-reach} satisfy
\[
G_X(v_j)=\delta_j u,\qquad \|v_j\|_\infty=O(\sqrt{\delta_j}),
\qquad q_X(v_j)=\delta_j u+O(\delta_j^2).
\]
Disjoint supports and \eqref{eq:localized-sharp-cost} give
\[
\|v_j\|_2^2=\gamma_{X,E_0}(q_X(v_j))
=\delta_j\gamma_{X,E_0}(u)+O(\delta_j^2).
\]
The support and dispersion bounds give the two restricted limits.
\end{proof}

\begin{corollary}[Necessary and sufficient sharp-cost scale]\label{cor:sharp-cost-scale}
Fix \(u\ne0\) and a supporting quadratic differential \(\phi\) with
\[
\|\phi\|_{\infty,h}=1,
\qquad \operatorname{Re}\langle u,\phi\rangle=\gamma_X(u).
\]
Every sufficiently small centered law of radius \(r\) and energy \(M\), with exact output \(\tau_X(\pi(\beta_\Pi))=\delta u\), satisfies
\begin{equation}\label{eq:nonlinear-capacity-obstruction}
\delta\gamma_X(u)
\le(1+C_Xr)M-\kappa_X\frac{M^3}{r^4}.
\end{equation}
Consequently,
\begin{equation}\label{eq:finite-radius-cost-gap}
\frac{\mathcal I_X(\tau_X^{-1}(\delta u);r)}{\delta}
\ge\frac{\gamma_X(u)}{1+C_Xr}
+\frac{\kappa_X\gamma_X(u)^3}{(1+C_Xr)^4}
\left(\frac{\delta}{r^2}\right)^2.
\end{equation}
Let \(r_j\to0\) and \(\delta_j\to0\) be positive. There exist centered laws of radius at most \(r_j\), realizing \(t_j=\tau_X^{-1}(\delta_j u)\) exactly for all sufficiently large \(j\), with
\begin{equation}\label{eq:sharp-sequence-energy}
\E\|v_j\|_2^2/\delta_j\longrightarrow\gamma_X(u),
\end{equation}
if and only if \(\delta_j/r_j^2\to0\). In that case the laws can be chosen to be equal-weight pairs. The same equivalence holds with \(\cU_2(\Pi_j)^2/(4\delta_j)\to\gamma_X(u)\) in place of \eqref{eq:sharp-sequence-energy}.
\end{corollary}

\begin{proof}
Pair \eqref{eq:uniform-law-output} with \(\phi\) and apply \cref{prop:bounded-amplitude-gap}. Since \(\|\phi\|_{\WP,*}\le\sqrt{\operatorname{Area}(X)}\), the paired nonlinear remainder is at most \(C_XrM\). This proves \eqref{eq:nonlinear-capacity-obstruction}. First discard its last term to obtain \(M\ge\delta\gamma_X(u)/(1+C_Xr)\). Substitute this lower bound into that term and take the infimum to get \eqref{eq:finite-radius-cost-gap}.

For a sequence satisfying \eqref{eq:sharp-sequence-energy}, divide that inequality by \(\delta_j\):
\[
\gamma_X(u)
\le(1+C_Xr_j)\frac{M_j}{\delta_j}
-\kappa_X\left(\frac{M_j}{\delta_j}\right)^3
\left(\frac{\delta_j}{r_j^2}\right)^2.
\]
The first term tends to \(\gamma_X(u)>0\). The last term must tend to zero, which gives the necessary scale. Sufficiency follows from the two-point upper bound in \cref{thm:sharp-exact-cost}: choose an exact pair within \(\delta_j/j\) of the restricted vertical-energy infimum. That infimum divided by \(\delta_j\) tends to \(\gamma_X(u)\), so the chosen pairs satisfy \eqref{eq:sharp-sequence-energy}. Finally, \eqref{eq:sharp-dispersion-moment} gives the sharp dispersion limit from the energy limit. Conversely, a sharp dispersion sequence satisfies \(M_j\le\cU_2(\Pi_j)^2/4\), while the exact-output lower bound gives \(M_j\ge\delta_j\gamma_X(u)/(1+C_Xr_j)\). These two bounds force the sharp energy limit.
\end{proof}

\subsection{A quantitative sharp-cost regime}

The bounded-amplitude excess also gives a rate for exact nonlinear realization. Write \(\mathcal I_X^{(2)}\) and \(\mathcal K_X^{(2)}\) for the corresponding infima restricted to the equal-weight pairs \(\Pi_v\).

\begin{theorem}[Quantitative exact creation]\label{thm:quantitative-exact-creation}
Fix \(\gamma_X(u)=1\). For all sufficiently small \(r,\lambda>0\), put \(\delta=r^2\lambda\) and \(t=\tau_X^{-1}(\delta u)\). Then
\begin{equation}\label{eq:quantitative-exact-upper}
\frac{\mathcal I_X(t;r)}\delta
\le\frac{\mathcal I_X^{(2)}(t;r)}\delta
\le1+C_{X,u}\left(\lambda+\frac{r^2}\lambda\right).
\end{equation}
If a normalized supporting differential for \(u\) has only nondegenerate maxima, then
\begin{equation}\label{eq:quantitative-exact-lower}
\frac{\mathcal I_X(t;r)}\delta
\ge\frac1{1+C_Xr}
+\frac{c_{X,u}\lambda}{(1+C_Xr)^3}.
\end{equation}
Consequently, along \(r_j,\lambda_j\to0\) with \(r_j=o(\lambda_j)\), all four quantities
\[
\frac{\mathcal I_X(t_j;r_j)}{\delta_j}-1,\quad
\frac{\mathcal I_X^{(2)}(t_j;r_j)}{\delta_j}-1,\quad
\frac{\mathcal K_X(t_j;r_j)}{4\delta_j}-1,\quad
\frac{\mathcal K_X^{(2)}(t_j;r_j)}{4\delta_j}-1
\]
are bounded above and below by positive constant multiples of \(\lambda_j\).
\end{theorem}

\begin{proof}
Use the direction \(w_\lambda\) in \cref{thm:capped-excess}, of amplitude at most \(1/2\), and reserve its fixed open complement \(O\). Choose a diagonal response basis in \(O\), with span \(E\), so that
\[
\|q_X(z)\|_{\WP}\ge a\|z\|_\infty^2,
\qquad \|z\|_2^2\le A\|z\|_\infty^2,
\qquad z\in E.
\]
Disjoint supports and \eqref{eq:midpoint-quartic-normal-form} give
\[
G_X(rw_\lambda+z)-r^2\lambda u
=q_X(z)+R_{r,\lambda}(z).
\]
On \(\|z\|_\infty\le Kr^2\), for fixed \(K\) and small \(r\), the total parameter has norm at most \(r\), and
\(\|R_{r,\lambda}(z)\|_{\WP}\le C_Xr^4\), uniformly in \(\lambda\). Choose \(K\) with \(aK^2>C_X\). The quadratic term dominates the remainder on the boundary, so the degree argument of \cref{thm:exact-reach} provides \(z_{r,\lambda}\) with exact output and \(\|z_{r,\lambda}\|_2^2\le AK^2r^4\). Orthogonality of supports gives
\[
\|rw_\lambda+z_{r,\lambda}\|_2^2
\le r^2(\lambda+C_{X,u}\lambda^2)+AK^2r^4,
\]
proving \eqref{eq:quantitative-exact-upper}.

For the lower bound, combine the uniform law expansion with the nondegenerate part of \cref{prop:bounded-amplitude-gap}. An exact law of energy \(M\) satisfies
\[
\delta\le(1+C_Xr)M-c_{X,u}M^2/r^2.
\]
First use \(M\ge\delta/(1+C_Xr)\), then substitute that bound in the last term. This proves \eqref{eq:quantitative-exact-lower}. If \(r=o(\lambda)\), the lower error is \(o(\lambda)\) and \(r^2/\lambda=o(\lambda)\). The dispersion comparison \eqref{eq:sharp-dispersion-moment} proves the remaining two assertions.
\end{proof}

\subsection{Restoration over the full local fiber}

Fix a sufficiently small convex \(L^\infty\)-ball \(U\subset V\) on which \cref{prop:fiber-graph-l2} holds, and set \(\mathcal F_{\mathrm{loc}}=\Phi(U)\). For a coefficient law define
\begin{equation}\label{eq:law-restoration-def}
\Delta_X(\Pi)=\inf_{\nu\in\mathcal F_{\mathrm{loc}}}
\bigl(\Energy_\Pi(\nu)-\Energy_\Pi(\beta_\Pi)\bigr).
\end{equation}
The Hadamard variance inequality \cite[Proposition~4.4]{sturm2003probability} gives, for every coefficient \(\nu\),
\begin{equation}\label{eq:variance-L2-lower}
\Energy_\Pi(\nu)-\Energy_\Pi(\beta_\Pi)
\ge\distB(\nu,\beta_\Pi)^2\ge4\|\nu-\beta_\Pi\|_2^2.
\end{equation}
We also need the local Hessian of the loss.

\begin{lemma}[Local loss estimates]\label{lem:distance-quartic}\label{lem:F-C2}\label{lem:ambient-hessian}
For \(|a|+|b|\) sufficiently small,
\begin{equation}\label{eq:distance-quartic}
\distH(a,b)^2=4|a-b|^2+R_4(a,b),\qquad
|D_a^jR_4(a,b)|\le C(|a|+|b|)^{4-j},\quad 0\le j\le2.
\end{equation}
The integrated loss is twice real Fr\'echet differentiable on coefficient balls. If its samples and evaluation point are \(O(r)\) in supremum norm, then
\begin{equation}\label{eq:F-hessian}
D^2\Energy_\Pi(\nu)[h,h]=8\|h\|_2^2+O(r^2)\|h\|_2^2,
\end{equation}
uniformly in the bounded direction \(h\).
\end{lemma}
\begin{proof}
The formula \eqref{eq:poincare-distance} is real analytic after squaring near the diagonal and has quadratic term \(4|a-b|^2\), with no cubic term. Its remainder and first two derivatives have the stated bounds. On compact coefficient disks, the pointwise derivatives are bounded; dominated convergence permits differentiation under expectation and integration. Integrating the second-derivative estimate proves \eqref{eq:F-hessian}.
\end{proof}

The following estimate is relative to the exact output. Its constants are uniform over the centered laws, their energies and their concentration.

\begin{proposition}[Full local-fiber restoration at an exact output]\label{prop:uniform-exact-restoration}
For every sufficiently small \(r>0\) and every centered law in \eqref{eq:centered-law-class}, put \(y=\tau_X(\pi(\beta_\Pi))\). Then
\begin{equation}\label{eq:uniform-exact-restoration}
\left|\Delta_X(\Pi)-4\|y\|_{\WP}^2\right|
\le C_Xr^2\|y\|_{\WP}^2.
\end{equation}
Equivalently,
\(\Delta_X(\Pi)=4d_{\WP}([X],\pi(\beta_\Pi))^2(1+O_X(r^2))\).
The constants are independent of the law and its energy.
\end{proposition}

\begin{proof}
Write \(\mu=\Phi(v)\), \(\beta=\beta_\Pi\), and use the orthogonal decomposition
\[
b=P_V\beta,\qquad h=P_H\beta,\qquad \beta=b+h.
\]
Put \(M=\E\|v\|_2^2\le\operatorname{Area}(X)r^2\). Centering gives \(\E\mu=\E\Psi(v)\in H\). The estimates in the proof of \cref{prop:uniform-law-expansion} yield
\begin{equation}\label{eq:uniform-restoration-vertical}
\|b\|_2\le C_Xr^2\sqrt M\le C_Xr^3,
\qquad \|\beta\|_2^2\le C_Xr^2M\le C_Xr^4.
\end{equation}
Pointwise, \(|\beta-\E\mu|\le C_Xr\E|\mu|^2\le C_Xr^3\), while \(\|\beta-\E\mu\|_1\le C_XrM\). The finite-dimensional harmonic projection is bounded from \(L^1\) to \(H\) in \(L^\infty\). Therefore
\[
\|b\|_\infty\le C_X(r^3+rM)\le C_Xr^3,
\]
and \(b\in U\) for small \(r\). Set
\[
\nu_0=\Phi(b),\qquad n=h-\Psi(b)\in H,
\qquad \beta=\nu_0+n.
\]
The graph estimate gives \(\|\nu_0\|_2\le C_Xr^3\). On the segment \(\nu_0+tn\), the \(L^2\) norm is at most \(C_Xr^2\), and the \(L^\infty\) norm remains in the fixed Bers neighborhood.

Let \(h_y\in H\) be the harmonic lift of \(y\). Since \(\mathfrak P(\nu_0)=0\), integrate the derivative estimate \eqref{eq:bers-nonlinear-first} along this segment. The normalized Bers coordinate gives
\[
\|y-D\pi_0[n]\|_{\WP}
\le C_Xr^2\|n\|_2.
\]
The restriction \(D\pi_0|_H\) is an isometry for the Weil--Petersson norm. Absorbing the right-hand side first bounds \(\|n\|_2\le C_X\|y\|_{\WP}\), and then gives
\begin{equation}\label{eq:uniform-restoration-normal}
\|n-h_y\|_2\le C_Xr^2\|y\|_{\WP}.
\end{equation}
If \(y=0\), this forces \(n=0\), so \(\beta=\nu_0\) lies in the fiber and \(\Delta_X(\Pi)=0\). Assume henceforth that \(\delta=\|y\|_{\WP}>0\). By \cref{prop:uniform-law-expansion} and \cref{thm:l2-tame}, \(\delta\le C_XM\le C_Xr^2\).

The samples, \(\nu_0\), and \(\beta\) are \(O(r)\) in \(L^\infty\). Since \(\beta\) minimizes the ambient loss, \cref{lem:ambient-hessian} gives
\[
\Energy_\Pi(\nu_0)-\Energy_\Pi(\beta)
=4\|n\|_2^2+O_X(r^2\|n\|_2^2)
=4\delta^2+O_X(r^2\delta^2).
\]
This proves the upper bound.

Choose a \(\delta^4\)-near-minimizer \(\nu=\Phi(w)\) over the full local graph. The variance inequality \eqref{eq:variance-L2-lower} gives
\[
\|\nu-\beta\|_2\le C_X\delta,
\qquad \|w-b\|_2\le C_X\delta,
\qquad \|w\|_2\le C_X(r^3+\delta).
\]
The fixed parameter neighborhood \(U\) is convex. Thus \eqref{eq:Psi-L2-Lipschitz} applies to \(w,b\in U\), giving
\[
\|\Psi(w)-\Psi(b)\|_2
\le C_X\delta(r^3+\delta)\le C_Xr^2\delta.
\]
The harmonic component of \(\nu-\beta\) is \(-n+\Psi(w)-\Psi(b)\). Equation~\eqref{eq:uniform-restoration-normal}, orthogonality, and the variance inequality give
\[
\Energy_\Pi(\nu)-\Energy_\Pi(\beta)
\ge4\|P_H(\nu-\beta)\|_2^2
\ge4\delta^2-C_Xr^2\delta^2.
\]
The near-minimality error satisfies \(\delta^4\le C_Xr^2\delta^2\) and is absorbed in this bound. Finally, expressing \(\tau_X\), whose derivative at \([X]\) is the identity, in Weil--Petersson normal coordinates gives
\(d_{\WP}([X],\tau_X^{-1}(y))^2=\|y\|_{\WP}^2+O_X(\|y\|_{\WP}^3)\).
Since \(\|y\|_{\WP}\le C_Xr^2\), this proves the distance formulation.
\end{proof}

\subsection{Simultaneous creation and restoration}

\begin{corollary}[Exact pairs with sharp creation and restoration]\label{cor:simultaneous-sharp-costs}
Let \(u\ne0\), \(r_j\to0\), and \(\delta_j/r_j^2\to0\), with \(r_j,\delta_j>0\). For all sufficiently large \(j\), there is \(v_j\in V\) such that
\[
\|v_j\|_\infty\le r_j,
\qquad \Pi_j=\tfrac12\delta_{\Phi(v_j)}+\tfrac12\delta_{\Phi(-v_j)},
\qquad \pi(\beta_{\Pi_j})=\tau_X^{-1}(\delta_j u),
\]
and this one sequence of exact two-point laws satisfies
\begin{equation}\label{eq:simultaneous-sharp-costs}
\frac{\|v_j\|_2^2}{\delta_j}\longrightarrow\gamma_X(u),
\qquad
\frac{\cU_2(\Pi_j)^2}{4\delta_j}\longrightarrow\gamma_X(u),
\qquad
\frac{\Delta_X(\Pi_j)}{4\delta_j^2}\longrightarrow\|u\|_{\WP}^2.
\end{equation}
Every sequence of centered exact-output laws with radius at most \(r_j\) has the same restoration limit. Conversely, simultaneous sharp creation and restoration for nonzero \(u\) requires \(\delta_j/r_j^2\to0\).
\end{corollary}

\begin{proof}
Choose the exact equal-weight pairs supplied by \cref{cor:sharp-cost-scale}, with \(\|v_j\|_2^2/\delta_j\to\gamma_X(u)\). The dispersion estimate \eqref{eq:sharp-dispersion-moment} gives the second limit. Applying \cref{prop:uniform-exact-restoration} to these same laws gives the third, and gives the stated restoration limit for every centered exact-output sequence. The necessary scale is again \cref{cor:sharp-cost-scale}.
\end{proof}

\subsection{The two-cost ratio and degeneration}

For the sharply creating pairs, the exact-output restoration estimate defines a directional ratio
\begin{equation}\label{eq:two-cost-ratio}
\mathcal R_X(u):=\frac{\|u\|_{\WP}^2}{4\gamma_X(u)^2}
=\lim_{j\to\infty}\frac{\Delta_X(\Pi_j)}{\cU_2(\Pi_j)^4},
\qquad u\ne0,
\end{equation}
for the laws in \cref{cor:simultaneous-sharp-costs}. This ratio connects the fiber problem to point evaluation in the Weil--Petersson cotangent space.

\begin{corollary}[Extremal ratio and the quadratic Bergman kernel]\label{cor:bergman-cost-ratio}
Let \(\phi_1,\ldots,\phi_n\) be a Weil--Petersson orthonormal basis of \(\QD(X)\), and let
\[
\mathcal B_X(p):=\sum_{k=1}^n|\phi_k(p)|_h^2
\]
be its hyperbolic quadratic Bergman kernel. Then
\begin{equation}\label{eq:extremal-bergman-ratio}
\sup_{u\ne0}\mathcal R_X(u)
=\frac14\sup_{\phi\ne0}\frac{\|\phi\|_{\infty,h}^2}{\|\phi\|_{\WP,*}^2}
=\frac14\max_{p\in X}\mathcal B_X(p).
\end{equation}
The supremum is realized by a spatial evaluation atom at a point maximizing \(\mathcal B_X\). For closed surfaces of fixed genus, there are \(c_g,C_g,\ell_g>0\) such that
\begin{equation}\label{eq:systolic-cost-ratio}
\frac{c_g}{\ell_{\mathrm{sys}}(X)}
\le\sup_{u\ne0}\mathcal R_X(u)
\le\frac{C_g}{\ell_{\mathrm{sys}}(X)},
\qquad \ell_{\mathrm{sys}}(X)\le\ell_g.
\end{equation}
For every fixed-genus sequence with \(\ell_{\mathrm{sys}}(X)\to0\),
\begin{equation}\label{eq:sharp-systolic-cost-ratio}
\ell_{\mathrm{sys}}(X)\sup_{u\ne0}\mathcal R_X(u)
\longrightarrow\frac{1}{2\pi}.
\end{equation}
Thus no constant depending only on genus bounds restoration energy by squared dispersion along all sharply creating exact two-point sequences.
\end{corollary}

\begin{proof}
Duality and \cref{thm:cotangent-calibration} give
\[
\sup_{\gamma_X(u)\le1}\|u\|_{\WP}
=\sup_{\|\phi\|_{\WP,*}\le1}\gamma_X^*(\phi)
=\sup_{\|\phi\|_{\WP,*}\le1}\|\phi\|_{\infty,h}.
\]
For fixed \(p\), the squared norm of point evaluation on \(\QD(X)\) is \(\mathcal B_X(p)\). Taking the supremum over \(p\) proves \eqref{eq:extremal-bergman-ratio}.

Choose \(p\) maximizing \(\mathcal B_X\), and let \(\phi_p\) be a unit Weil--Petersson cotangent vector attaining the evaluation norm there. Cauchy--Schwarz in the orthonormal basis gives
\[
|\phi_p(z)|_h\le\sqrt{\mathcal B_X(z)}
\le\sqrt{\mathcal B_X(p)}=|\phi_p(p)|_h.
\]
Choose the phase of the atom \(e_X(p,a)\) from \cref{cor:geometric-response-body} so that its pairing with \(\phi_p\) is positive. That atom lies in the \(\gamma_X\)-unit ball, and \(\phi_p/\|\phi_p\|_{\infty,h}\) pairs with it to give one. Hence its cost norm is one, while its Weil--Petersson norm is \(\sqrt{\mathcal B_X(p)}\). It realizes the asserted extremal ratio.

For the lower systolic bound, put \(\ell=\ell_{\mathrm{sys}}(X)\). The collar concentration differential in \cite[Lemma~2.6]{rupflin2015uniform}, normalized to have \(\|\Omega\|_{\WP,*}=1\), satisfies
\[
\Omega=b_0dz^2+e\quad\hbox{on the collar},
\qquad \|e\|_{\infty,h}\le C_g\sqrt\ell,
\qquad |b_0|\|dz^2\|_{2,\mathrm{collar}}=1+O_g(\sqrt\ell).
\]
In conformal collar coordinates, \(\rho(s)=\ell/(2\pi\cos(\ell s/(2\pi)))\). Thus \(\|dz^2\|_{2,\mathrm{collar}}\asymp\ell^{-3/2}\), and on the central geodesic \(|dz^2|_h=(2\pi/\ell)^2\). More precisely, in the norm convention of this paper,
\[
\|dz^2\|_{2,\mathrm{collar}}^2
=2\pi\int_{-S_\ell}^{S_\ell}\rho(s)^{-2}\,ds
=\frac{8\pi^5}{\ell^3}(1+O(\ell^3)),
\]
where \(S_\ell=(2\pi/\ell)(\pi/2-\arctan\sinh(\ell/2))\). It follows that
\[
\liminf_{\ell\to0}\ell
\sup_{\phi\ne0}\frac{\|\phi\|_{\infty,h}^2}{\|\phi\|_{\WP,*}^2}
\ge\frac2\pi.
\]
For the upper bound, the sharp pointwise estimate of \cite[Theorem~3.7]{bridgeman2021uniform} gives, for every \(\eta>0\) and all surfaces with sufficiently small systole,
\[
\|\phi\|_{\infty,h}
\le(1+\eta)\sqrt{\frac{2}{\pi\ell_{\mathrm{sys}}(X)}}\,\|\phi\|_{\WP,*}.
\]
Take the limsup and then \(\eta\downarrow0\). Equation~\eqref{eq:extremal-bergman-ratio} gives \eqref{eq:sharp-systolic-cost-ratio}, and hence \eqref{eq:systolic-cost-ratio}.

Finally choose the maximizing direction on each surface in a degenerating sequence, and apply \cref{cor:simultaneous-sharp-costs} separately on each surface. Equation~\eqref{eq:two-cost-ratio} and the lower bound in \eqref{eq:systolic-cost-ratio} make the limiting ratio unbounded. The local realization scales may be chosen separately for the different surfaces.
\end{proof}

The sharp comparison grows as \((2\pi\ell_{\mathrm{sys}})^{-1}\). We next establish nonlinear estimates uniform in the systole, before studying the degeneration of the creation norm itself.

\section{Uniform nonlinear realization}\label{sec:uniform-scales}

Put \(\sigma_X=\min\{1,\ell_{\mathrm{sys}}(X)\}\). We give a quantitative realization range with constants depending only on genus and normalize \(\gamma_X(u)=1\). For a general nonzero direction, replace \(u\) by \(u/\gamma_X(u)\) and \(\delta\) by \(\delta\gamma_X(u)\). Throughout this section the restoration infimum is taken over \(\Phi(U_X)\), where \(U_X\) is the systole-dependent graph neighborhood fixed in \cref{prop:systole-nonlinear}.

\subsection{Uniform graph and law estimates}

\begin{proposition}[Systole bounds for the nonlinear errors]\label{prop:systole-nonlinear}
There are a universal \(\varepsilon_0>0\) and \(\rho_g,C_g^{\mathrm{an}}>0\) with the following properties on every closed genus-\(g\) surface. The fiber graph exists for \(\|v\|_\infty<\varepsilon_0\) and satisfies
\begin{equation}\label{eq:uniform-infinity-graph}
\|\Psi(v)\|_\infty\le C_g^{\mathrm{an}}\|v\|_\infty^2,
\qquad
\|G_X(v)-q_X(v)\|_{\gamma_X}\le C_g^{\mathrm{an}}\|v\|_\infty^4.
\end{equation}
On \(U_X=\{v\in V:\|v\|_\infty<\rho_g\sigma_X^2\}\), the quadratic graph and graph-derivative bounds \eqref{eq:Psi-L2-quadratic}--\eqref{eq:DPsi-L2} hold with constant \(C_g^{\mathrm{an}}\sigma_X^{-2}\), and the cubic bounds \eqref{eq:Psi-L2-cubic}--\eqref{eq:DPsi-L2-cubic} hold with constant \(C_g^{\mathrm{an}}\sigma_X^{-4}\).

For every centered law with vertical radius \(r\le \rho_g\sigma_X^4\), put \(M=\E\|v\|_2^2\) and \(y=\tau_X(\pi(\beta_\Pi))\). Then
\begin{align}
\gamma_X\bigl(y-\E q_X(v)\bigr)
&\le C_g^{\mathrm{an}}\sigma_X^{-4}rM,\label{eq:systole-law-error}\\
4M\le\cU_2(\Pi)^2
&\le4(1+C_g^{\mathrm{an}}\sigma_X^{-4}r^2)M,\label{eq:systole-dispersion-error}\\
\left|\Delta_X(\Pi)-4\|y\|_{\WP}^2\right|
&\le C_g^{\mathrm{an}}\sigma_X^{-4}r^2\|y\|_{\WP}^2.\label{eq:systole-restoration-error}
\end{align}
Here the restoration infimum is taken over the full graph \(\Phi(U_X)\).
\end{proposition}

\begin{proof}
See Appendix~\ref{app:uniform-proofs}.
\end{proof}

\subsection{Quantitative concentration and a disjoint response basis}

\begin{lemma}[Uniform localized realization]\label{lem:uniform-localized-realization}
There are genus-dependent constants \(\eta_g,A_g>0\) and a fixed response bound \(\kappa_g=1/2\) with the following properties. For \(\gamma_X(u)=1\) and \(0<h\le\eta_g\sigma_X\), there is a smooth vertical direction \(v_0\), supported in finitely many disjoint disks, such that
\begin{equation}\label{eq:uniform-near-optimal-direction}
q_X(v_0)=u,\qquad
\|v_0\|_2^2\le1+A_gh,\qquad
\|v_0\|_\infty\le A_g/h.
\end{equation}
There are also \(n=3g-3\) directions \(a_i\) on mutually disjoint disks outside this support, with \(\|a_i\|_2=1\), \(\|a_i\|_\infty\le A_g/h\), and
\begin{equation}\label{eq:uniform-response-basis}
q_X\!\left(\sum_i z_ia_i\right)=\sum_i z_i^2q_X(a_i),
\qquad
\gamma_X\!\left(\sum_i z_i^2q_X(a_i)\right)
\ge\kappa_g\max_i|z_i|^2.
\end{equation}
\end{lemma}

\begin{proof}
See Appendix~\ref{app:uniform-basis-proof}.
\end{proof}

\subsection{An explicit realization range}

\begin{theorem}[Uniform exact two-point realization]\label{thm:uniform-exact-realization}
There are \(c_g,C_g>0\) such that, for every closed genus-\(g\) surface, every \(u\) with \(\gamma_X(u)=1\), and parameters
\begin{equation}\label{eq:uniform-realization-range}
0<r\le c_g\sigma_X^4,
\qquad 0<h\le c_g\sigma_X,
\qquad 0<\delta\le c_g\min\{r^2h^2,h^4\},
\end{equation}
there is an equal-weight pair \(\Pi_v\) in the fiber over \(X\) with \(\|v\|_\infty\le r\) and exact output \(\tau_X^{-1}(\delta u)\). It satisfies
\begin{align}
\frac1{1+C_g\sigma_X^{-4}r}
\le\frac{\mathcal I_X(\tau_X^{-1}(\delta u);r)}\delta
&\le\frac{\|v\|_2^2}\delta
\le1+C_g\left(h+\frac\delta{h^4}\right),\label{eq:uniform-exact-creation}\\
1-C_g\sigma_X^{-4}r
\le\frac{\mathcal K_X(\tau_X^{-1}(\delta u);r)}{4\delta}
&\le\frac{\cU_2(\Pi_v)^2}{4\delta}
\le1+C_g\left(h+\frac\delta{h^4}+\sigma_X^{-4}r^2\right).\label{eq:uniform-exact-dispersion}
\end{align}
Every centered law of radius at most \(r\) with this exact output obeys
\begin{equation}\label{eq:uniform-exact-restoration-systole}
\left|\frac{\Delta_X(\Pi)}{4\delta^2\|u\|_{\WP}^2}-1\right|
\le C_g\sigma_X^{-4}r^2.
\end{equation}
\end{theorem}

\begin{proof}
We first fix the constants that enter the degree argument. Let \(A_g,\eta_g,\kappa_g\) be as in \cref{lem:uniform-localized-realization}, and let \(Q_g=C_g^{\mathrm{an}}\) bound the quartic midpoint remainder in \eqref{eq:uniform-infinity-graph}. Increase \(A_g\) to at least one and choose
\begin{equation}\label{eq:uniform-degree-constant-choice}
K_g^2=1+\frac{2Q_gA_g^4}{\kappa_g}.
\end{equation}
Only after fixing these constants, choose the realization threshold so that
\begin{equation}\label{eq:uniform-final-constant-choice}
0<c_g\le
\min\left\{\rho_g/2,\eta_g,K_g^{-2},A_g^{-2},1\right\}.
\end{equation}
The first two bounds impose the analytic and localization ranges; the third controls the correction amplitude, and the fourth imposes the prescribed vertical radius. The response lower bound \(\kappa_g\) is unchanged by this choice.

Take \(v_0,a_i\) from \cref{lem:uniform-localized-realization}. For \(w=\sum_i z_ia_i\), consider the polydisk \(\max_i|z_i|\le K_g\delta/h^2\). Since \(\delta/h^4\le c_g\le K_g^{-2}\),
\[
\|w\|_\infty\le A_gK_g\delta/h^3\le A_g\sqrt\delta/h.
\]
Disjoint support, \(\delta\le c_gr^2h^2\), and \(c_g\le A_g^{-2}\) give, throughout this polydisk,
\[
\|\sqrt\delta\,v_0+w\|_\infty
\le A_g\sqrt\delta/h\le r\le c_g\sigma_X^4
\le\tfrac12\rho_g\sigma_X^4.
\]
Thus the whole parameter region is in the analytic range of \cref{prop:systole-nonlinear} before the homotopy is formed. The midpoint expansion gives
\[
G_X(\sqrt\delta\,v_0+w)-\delta u=q_X(w)+R_\delta(w),
\qquad
\gamma_X(R_\delta(w))\le Q_gA_g^4\delta^2/h^4.
\]

On the boundary, \eqref{eq:uniform-response-basis} gives
\[
\gamma_X(q_X(w))\ge\kappa_gK_g^2\delta^2/h^4
>Q_gA_g^4\delta^2/h^4.
\]
The homotopy \(q_X(w)+sR_\delta(w)\) has no boundary zero. Its coordinate-square degree \(2^n\) gives a zero \(w_\delta\), with
\[
\|w_\delta\|_2^2=\sum_i|z_i|^2
\le nK_g^2\delta^2/h^4.
\]
Thus \(v=\sqrt\delta\,v_0+w_\delta\) has exact output \(\delta u\), and
\[
\|v\|_2^2\le\delta(1+A_gh)+nK_g^2\delta^2/h^4.
\]
The preceding polydisk bound gives \(\|v\|_\infty\le r\).

For any centered law with this exact output, \eqref{eq:systole-law-error} gives
\(\delta\le(1+C_g^{\mathrm{an}}\sigma_X^{-4}r)M\).
Combine this with \eqref{eq:systole-dispersion-error} and the constructed pair's energy to obtain both cost inequalities. Fix the final upper constant \(C_g\) to dominate the analytic, energy, and dispersion constants; this changes none of the preceding choices. Equation~\eqref{eq:systole-restoration-error} proves the restoration estimate.
\end{proof}

\begin{corollary}[Prescribed scales on degenerating surfaces]\label{cor:prescribed-degenerating-scales}
Let \(X_j\) be closed genus-\(g\) surfaces, let \(\gamma_{X_j}(u_j)=1\), and set \(\sigma_j=\sigma_{X_j}\). For positive parameters \(r_j,\delta_j\) satisfying
\begin{equation}\label{eq:prescribed-sharp-scales}
\frac{r_j}{\sigma_j^4}\longrightarrow0,
\qquad
\frac{\delta_j}{r_j^2\sigma_j^2}\longrightarrow0,
\end{equation}
there are exact equal-weight pairs, defined for all sufficiently large \(j\), of radius at most \(r_j\), with output \(\tau_{X_j}^{-1}(\delta_j u_j)\), satisfying
\begin{equation}\label{eq:uniform-joint-sharp-costs}
\frac{\|v_j\|_2^2}{\delta_j}\to1,
\qquad \frac{\cU_2(\Pi_j)^2}{4\delta_j}\to1,
\qquad
\frac{\Delta_{X_j}(\Pi_j)}{4\delta_j^2\|u_j\|_{\WP}^2}\to1.
\end{equation}
In particular, for \(0<\eta<c_g\), the explicit choice
\begin{equation}\label{eq:explicit-power-scales}
r=c_g\eta\sigma_X^4,
\qquad \delta=c_g^6\eta^6\sigma_X^{10}
\end{equation}
gives creation and normalized dispersion \(1+O_g(\eta)\), and restoration ratio \(1+O_g(\eta^2\sigma_X^4)\).
\end{corollary}

\begin{proof}
Put \(\lambda_j=\delta_j/(r_j^2\sigma_j^2)\) and choose \(h_j=c_g\sigma_j\lambda_j^{1/4}\). Then
\[
\frac{\delta_j}{r_j^2h_j^2}=c_g^{-2}\lambda_j^{1/2}\to0,
\qquad
\frac{\delta_j}{h_j^4}=c_g^{-4}\frac{r_j^2}{\sigma_j^2}\to0,
\qquad h_j\to0.
\]
Theorem~\ref{thm:uniform-exact-realization} applies for all large \(j\) and proves \eqref{eq:uniform-joint-sharp-costs}. For \eqref{eq:explicit-power-scales}, use \(h=c_g\eta\sigma_X\). Then
\[
\frac{\delta}{r^2h^2}=c_g^2\eta^2,\qquad
\frac{\delta}{h^4}=c_g^2\eta^2\sigma_X^6.
\]
Since \(c_g\le1\) and \(0<\eta<c_g\), both ratios are at most \(c_g\), so the constants already chosen satisfy \eqref{eq:uniform-realization-range}. Substitution in the three error bounds gives the stated estimates.
\end{proof}

The powers in \eqref{eq:explicit-power-scales} give a sufficient realization range on all of moduli space. The sharp-cost regime on each fixed surface remains \(\delta=o(r^2)\), as proved in \cref{cor:sharp-cost-scale}.

\section{Boundary geometry of the creation metric}\label{sec:cost-degeneration}

The cotangent calibration determines the degeneration of the entire creation norm. The limit separates the quadratic differentials that persist on the punctured components from one complex principal part at each node. The resulting tangent norm is additive across these pieces.

Let \(X_j\) be closed hyperbolic surfaces of genus \(g\), converging to a finite-area punctured surface \(\Sigma\) by pinching the disjoint curves \(\alpha_1,\ldots,\alpha_k\). Write \(\ell_{i,j}=\ell_{\alpha_i}(X_j)\to0\) and \(C_{i,j}\) for their standard collars. In collar coordinates \(z=s+i\vartheta\), with \(\vartheta\in\mathbb R/(2\pi\mathbb Z)\), set
\[
\rho_{i,j}(s)=\frac{\ell_{i,j}}{2\pi\cos(\ell_{i,j}s/(2\pi))},
\qquad
\phi=\sum_{m\in\mathbb Z}b_m^{i,j}(\phi)e^{mz}dz^2.
\]
Define the persistent subspace
\begin{equation}\label{eq:persistent-cotangent}
W_j:=\{\phi\in\QD(X_j):b_0^{i,j}(\phi)=0,\ 1\le i\le k\}.
\end{equation}
On \(\Sigma\), let \(\QD(\Sigma)\) consist of integrable holomorphic quadratic differentials, equivalently differentials with at most simple poles at the punctures. Their hyperbolic pointwise norm tends to zero at each cusp. Write \(\gamma_\Sigma\) for the norm on \(\QD(\Sigma)^*\) dual to \(\|\cdot\|_{\infty,h}\).

The collar normalization is fixed by the direct integral
\begin{equation}\label{eq:collar-L2-normalization}
\begin{aligned}
\|dz^2\|_{2,C_{i,j}}^2
&=2\pi\int_{-S_{\ell_{i,j}}}^{S_{\ell_{i,j}}}\rho_{i,j}(s)^{-2}\,ds\\
&=\frac{16\pi^4}{\ell_{i,j}^3}
\left(\frac\pi2-x_{i,j}+\frac12\sin(2x_{i,j})\right)
=\frac{8\pi^5}{\ell_{i,j}^3}\bigl(1+O(\ell_{i,j}^3)\bigr),
\end{aligned}
\end{equation}
where \(x_{i,j}=\arctan\sinh(\ell_{i,j}/2)\) and
\(S_\ell=(2\pi/\ell)(\pi/2-\arctan\sinh(\ell/2))\). Thus the leading \(L^2\)-normalized principal coefficient is of order \(\ell^{3/2}\), whereas the chosen principal coefficient \(\ell^2\) has central-circle hyperbolic norm exactly \(4\pi^2\). All collar coefficients below use this normalization.

\begin{lemma}[Normalized collar differentials]\label{lem:normalized-collar}
For all large \(j\), there are complex-linear isomorphisms \(L_j:\QD(\Sigma)\to W_j\), preserving the Weil--Petersson cotangent inner product, whose values converge smoothly locally to their arguments under the degeneration identifications. There are unique \(\chi_{i,j}\in W_j^{\perp_{\WP,*}}\) with
\begin{equation}\label{eq:collar-principal-normalization}
b_0^{m,j}(\chi_{i,j})=\delta_{im}\ell_{i,j}^2.
\end{equation}
They satisfy
\begin{equation}\label{eq:normalized-collar-error}
\|\chi_{i,j}\|_{\infty,h;X_j\setminus C_{i,j}}
+\|\chi_{i,j}-\ell_{i,j}^2dz^2\|_{\infty,h;C_{i,j}}
\le C\ell_{i,j}^2.
\end{equation}
For every fixed \(\varepsilon>0\), their norms on the \(\varepsilon\)-thick part tend to zero. Moreover,
\begin{equation}\label{eq:persistent-thin-decay}
\|L_jw\|_{\infty,h;\varepsilon\text{-thin}(X_j)}
\le C\varepsilon^{-2}e^{-\pi/\varepsilon}\|w\|_{\WP,*},
\end{equation}
for all sufficiently small \(\varepsilon\) and all large \(j\).
\end{lemma}

\begin{proof}
The persistent-space convergence and thin-part decay follow from \cite[Lemma~2.4(i), (iii)]{rupflin2013asymptotics}. Choose converging orthonormal bases and define \(L_j\) by their correspondence. Fix a sufficiently small \(\bar\eta>0\) so that, for all large \(j\), the pinching curves are exactly the simple closed geodesics of length at most \(2\bar\eta\). Then \cite[Lemma~2.1, (2.11)--(2.12)]{grosse2019holomorphic} applies with constants independent of \(j\), and the principal-part map on \(W_j^\perp\) is an isomorphism. It gives unit differentials \(\Omega_{i,j}\) with positive principal coefficient \(b_{i,j}\) on collar \(i\), vanishing principal coefficients on the other pinching collars, and
\[
\|\Omega_{i,j}\|_{\infty,h;X_j\setminus C_{i,j}}
+\|\Omega_{i,j}-b_{i,j}dz^2\|_{\infty,h;C_{i,j}}
\le C\ell_{i,j}^{3/2},
\qquad
b_{i,j}\|dz^2\|_{2,C_{i,j}}=1+O(\ell_{i,j}^3).
\]
The collar integral computed in the proof of \cref{cor:bergman-cost-ratio} gives \(b_{i,j}\asymp\ell_{i,j}^{3/2}\). Thus
\[
\chi_{i,j}=\frac{\ell_{i,j}^2}{b_{i,j}}\Omega_{i,j}
\]
has the required principal parts and satisfies \eqref{eq:normalized-collar-error}. On a fixed thick part, \(\rho_{i,j}\) has a positive lower bound, so the principal term \(\ell_{i,j}^2\rho_{i,j}^{-2}\) also tends to zero. This proves the remaining assertion.
\end{proof}

Set
\begin{equation}\label{eq:degeneration-identification}
T_j(w,a):=L_jw+\sum_{i=1}^ka_i\chi_{i,j},
\qquad (w,a)\in\QD(\Sigma)\oplus\mathbb C^k.
\end{equation}
For \(u_0\in\QD(\Sigma)^*\) and \(c\in\mathbb C^k\), define \(u_j\in T_{[X_j]}\Teich\) by
\begin{equation}\label{eq:degeneration-tangent-identification}
\langle u_j,T_j(w,a)\rangle_T
=\langle u_0,w\rangle+\sum_{i=1}^kc_ia_i.
\end{equation}

\begin{theorem}[Additive degeneration of creation cost]\label{thm:additive-cost-degeneration}
Uniformly on bounded subsets of \(\QD(\Sigma)\oplus\mathbb C^k\),
\begin{equation}\label{eq:cotangent-max-limit}
\|T_j(w,a)\|_{\infty,h}
\longrightarrow
N_0(w,a):=\max\{\|w\|_{\infty,h},\ 4\pi^2\max_i|a_i|\}.
\end{equation}
Under \eqref{eq:degeneration-tangent-identification}, uniformly on bounded sets of \((u_0,c)\),
\begin{equation}\label{eq:additive-cost-limit}
\gamma_{X_j}(u_j)
\longrightarrow
\gamma_\Sigma(u_0)+\frac1{4\pi^2}\sum_{i=1}^k|c_i|.
\end{equation}
If \(\Sigma\) has components \(\Sigma_\nu\) and \(u_0=(u_\nu)_\nu\), then \(\gamma_\Sigma(u_0)=\sum_\nu\gamma_{\Sigma_\nu}(u_\nu)\).
\end{theorem}

\begin{proof}
Fix \((w,a)\). On every compact subset of \(\Sigma\), \(T_j(w,a)\) converges to \(w\), so
\[
\liminf_j\|T_j(w,a)\|_{\infty,h}\ge\|w\|_{\infty,h}.
\]
On the central circle of \(C_{i,j}\), the mean Fourier coefficient is exactly \(a_i\ell_{i,j}^2\), and the conformal factor is \(\ell_{i,j}/(2\pi)\). Taking the mean on that circle gives the further lower bound
\[
\|T_j(w,a)\|_{\infty,h}\ge4\pi^2|a_i|.
\]

For the upper bound, fix a small \(\varepsilon>0\). On the \(\varepsilon\)-thick part, \cref{lem:normalized-collar} and smooth local convergence give a limsup bounded by \(\|w\|_{\infty,h}\). On the \(\varepsilon\)-thin part of \(C_{i,j}\), the persistent term is bounded by \eqref{eq:persistent-thin-decay}, the other collar differentials tend uniformly to zero, and
\[
|\ell_{i,j}^2dz^2|_h
=4\pi^2\cos^2(\ell_{i,j}s/(2\pi))\le4\pi^2.
\]
Consequently,
\[
\limsup_j\|T_j(w,a)\|_{\infty,h}
\le\max\{\|w\|_{\infty,h},\ 4\pi^2\max_i|a_i|\}
+C\varepsilon^{-2}e^{-\pi/\varepsilon}\|w\|_{\WP,*}.
\]
Let \(\varepsilon\downarrow0\). The same estimates uniformly bound the operator norms of \(T_j\) on the fixed finite-dimensional parameter space. A finite-net argument therefore upgrades pointwise convergence to uniform convergence on bounded sets.

Uniform convergence on the \(N_0\)-unit sphere yields numbers \(\eta_j\to0\) such that
\[
(1-\eta_j)N_0\le\|T_j(\cdot)\|_{\infty,h}\le(1+\eta_j)N_0.
\]
By \cref{thm:cotangent-calibration}, \(\gamma_{X_j}(u_j)\) is the dual of the middle norm. The dual of the maximum norm \(N_0\) is
\[
N_0^*(u_0,c)=\gamma_\Sigma(u_0)+\frac1{4\pi^2}\sum_i|c_i|:
\]
the choices of \(w\) and of the phases of \(a_i\) are independent in its product unit ball. Taking polars proves \eqref{eq:additive-cost-limit}. Applying the same argument to the components of \(\Sigma\) proves the last assertion.
\end{proof}

\subsection{Length and angular twist}

Choose a pants decomposition containing the pinching curves. Use the angular twist convention in \cref{tab:normalizations}: \(\theta_i=2\pi s_i/\ell_i\), with Dehn period \(2\pi\). Partial derivatives below hold all other Fenchel--Nielsen coordinates fixed. Let \(u_j^0\) denote the persistent dual lift specified by
\[
\langle u_j^0,T_j(w,a)\rangle_T=\langle u_0,w\rangle.
\]

\begin{corollary}[The normal coordinate scales]\label{cor:fn-creation-scales}
Suppose the remaining Fenchel--Nielsen data stay in a compact set away from further degeneration. For fixed \(s,t\in\mathbb R^k\),
\begin{equation}\label{eq:fn-cost-limit}
\gamma_{X_j}\!\left(u_j^0+\sum_i
\left(s_i\partial_{\ell_i}+\frac{t_i}{\ell_{i,j}^2}\partial_{\theta_i}\right)\right)
\longrightarrow
\gamma_\Sigma(u_0)+\sum_i
\sqrt{\left(\frac\pi2s_i\right)^2+\left(\frac{t_i}{4\pi}\right)^2}.
\end{equation}
The convergence is uniform on bounded sets of \((u_0,s,t)\). In particular,
\begin{equation}\label{eq:length-twist-sharp-scales}
\gamma_{X_j}(\partial_{\ell_i})\longrightarrow\frac\pi2,
\qquad
\frac{\gamma_{X_j}(\partial_{\theta_i})}{\ell_{i,j}^2}
\longrightarrow\frac1{4\pi}.
\end{equation}
On every closed surface and for every pants curve,
\begin{equation}\label{eq:angular-twist-global-upper}
\gamma_X(\partial_{\theta_i})\le\frac{\ell_i^2}{4\pi}.
\end{equation}
\end{corollary}

\begin{proof}
By \cref{tab:normalizations}, a horizontal variation \(\operatorname{Re}(\Upsilon)\) pairs with \(\phi\) as \(\tfrac12\langle\phi,\Upsilon\rangle_{\WP,*}\), with the cotangent inner product linear in its first argument.

Let \(\Psi_i^{\mathrm{FN}}\) be the quadratic differential representing the metric derivative in \(\theta_i\). The angular twist identity \cite[(2.35)]{grosse2019holomorphic} is stated for the Hermitian tensor inner product \(\langle\cdot,\cdot\rangle_{\mathrm{tensor}}=4\langle\cdot,\cdot\rangle_{\WP,*}\). Since the real tensor pairing is one half its real part, while the Teichm\"uller cotangent pairing with a variation \(\operatorname{Re}(\Upsilon)\) is \(\frac12\langle\phi,\Upsilon\rangle_{\WP,*}\), it gives
\begin{equation}\label{eq:angular-twist-pairing}
\langle\partial_{\theta_i},\phi\rangle_T
=i\pi b_0^i(\phi).
\end{equation}
It follows that \(\ell_{i,j}^{-2}\partial_{\theta_i}\) has exactly the dual coordinates \((0,i\pi e_i)\) in \eqref{eq:degeneration-tangent-identification}. The mean Fourier coefficient on the central circle satisfies
\[
|b_0^i(\phi)|\le\frac{\ell_i^2}{4\pi^2}\|\phi\|_{\infty,h}.
\]
This proves \eqref{eq:angular-twist-global-upper}.

For the length derivative, \cite[Proposition~1.1 and Theorem~1.2]{grosse2019holomorphic}, in particular (1.7), (1.13) and (1.14) there, give its metric differential \(\Lambda_i^{\mathrm{FN}}\) in the form
\[
\Lambda_i^{\mathrm{FN}}=-\frac{\ell_i}{2\pi^2}dz^2+R_i
\quad\hbox{on }C_i,
\qquad
\|R_i\|_{\infty,h;C_i}
+\|\Lambda_i^{\mathrm{FN}}\|_{\infty,h;X\setminus C_i}\le C\ell_i.
\]
Fourier orthogonality and the collar integral then yield, uniformly for all \(\phi\in\QD(X)\),
\begin{equation}\label{eq:length-principal-pairing}
\langle\partial_{\ell_i},\phi\rangle_T
=-\frac{2\pi^3}{\ell_i^2}b_0^i(\phi)
+O\bigl(\ell_i\|\phi\|_{\infty,h}\bigr).
\end{equation}
Indeed, the main integral is
\(-\ell_i b_0^i(\phi)\|dz^2\|_{2,C_i}^2/(4\pi^2)\), and the remaining pairing is bounded by \(C\ell_i\|\phi\|_{1,h}\). The latter is at most \(C\ell_i\|\phi\|_{\infty,h}\), since the total area is fixed. The \(O(\ell_i^3)\) relative error in the collar integral is absorbed using the central-circle bound for \(b_0^i(\phi)\).

Evaluating \eqref{eq:length-principal-pairing} on \(T_j(w,a)\) shows that \(\partial_{\ell_i}\) has limiting dual coordinates \((0,-2\pi^3e_i)\). Equations~\eqref{eq:angular-twist-pairing} and \eqref{eq:additive-cost-limit} therefore give \eqref{eq:fn-cost-limit}, since
\[
\frac{|-2\pi^3s_i+i\pi t_i|}{4\pi^2}
=\sqrt{\left(\frac\pi2s_i\right)^2+\left(\frac{t_i}{4\pi}\right)^2}.
\]
The uniform estimates prove uniform convergence on bounded sets.
\end{proof}

\subsection{Distance to a pinching stratum}

The norm \(\gamma\) defines a length metric
\[
d_\gamma(X,Y):=\inf_{c:X\to Y}\int\gamma_{c(t)}(\dot c(t))\,dt.
\]
Write \(\widehat\Teich_\gamma\) for its metric completion. Cotangent norm comparison gives
\begin{equation}\label{eq:gamma-WP-distance-domination}
\gamma_X(u)\le\sqrt{4\pi(g-1)}\,\|u\|_{\WP},
\qquad
d_\gamma\le\sqrt{4\pi(g-1)}\,d_{\WP}.
\end{equation}
Consequently the identity extends to a Lipschitz map
\[
\iota:\overline\Teich^{\WP}\longrightarrow\widehat\Teich_\gamma.
\]
Here \(\overline\Teich^{\WP}\) is augmented Teichm\"uller space. Masur established the boundary extension of the Weil--Petersson metric \cite{masur1976extension}; its completion and stratum structure are described in \cite{wolpertgeometry}. For a multicurve \(\Gamma\), let \(D_\Gamma\) be the image under \(\iota\) of its closed pinching stratum, including further degenerations.

\begin{lemma}[Uniform calibration of short length functions]\label{lem:short-length-calibration}
There are \(\ell_g,C_g>0\) such that, on every closed surface of genus \(g\), the real differential of a geodesic length \(\ell_\alpha<\ell_g\) is represented by a cotangent vector \(\kappa_\alpha\):
\[
d\ell_\alpha(u)=\operatorname{Re}\langle u,\kappa_\alpha\rangle_T.
\]
For any disjoint collection of such curves and complex numbers \(\lambda_i\),
\begin{equation}\label{eq:short-length-cometric}
\left\|\sum_i\lambda_i\kappa_{\alpha_i}\right\|_{\infty,h}
\le\frac2\pi\max_i|\lambda_i|
+C_g\sum_i|\lambda_i|\ell_{\alpha_i}^2.
\end{equation}
\end{lemma}

\begin{proof}
Use a pants decomposition containing all sufficiently short curves, with lengths bounded in terms of genus; the extension and bounds are recorded in \cite[Remark~2.11]{grosse2019holomorphic}. The length variation identity \cite[(2.8)]{grosse2019holomorphic}, together with \eqref{eq:angular-twist-pairing}, gives
\[
\kappa_{\alpha_i}=-\frac{2\pi i}{\ell_i}\Psi_i^{\mathrm{FN}}.
\]
For example, a metric variation \(\operatorname{Re}(\Upsilon)\) pairs with this differential to give \(-2\pi^2\operatorname{Re}b_0^i(\Upsilon)/\ell_i\), the derivative of \(\ell_i\).

The normalized twist estimates of \cite[Proposition~1.1 and Theorem~1.3]{grosse2019holomorphic}, with their pointwise tensor norm twice ours, imply
\[
\kappa_{\alpha_i}
=-\frac{\ell_i^2}{2\pi^3}(1+O_g(\ell_i^3))dz^2+E_i
\quad\hbox{on }C_i,
\qquad
\|E_i\|_{\infty,h;C_i}
+\|\kappa_{\alpha_i}\|_{\infty,h;X\setminus C_i}
\le C_g\ell_i^2.
\]
To check the error, the unit collar differential has decay part \(O_g(\ell_i^{3/2})\), and the twist differential has size \(O_g(\ell_i^{3/2})\) in cotangent \(L^2\). The other basis coefficients in Theorem~1.3 are \(O_g(\ell_i^{3/2}\ell_m^{3/2})\); their pointwise contribution after multiplication by \(\ell_i^{-1}\) is \(O_g(\ell_i^2\ell_m)\). The pants lengths are bounded, so these contributions have the stated bound. The principal coefficient follows from \(\|dz^2\|_{2,C_i}^2=8\pi^5\ell_i^{-3}(1+O(\ell_i^3))\).

The principal terms have disjoint supports and pointwise maximum \(2/\pi+O_g(\ell_i^3)\). Summing the errors proves \eqref{eq:short-length-cometric}.
\end{proof}

\subsection{Identification of the metric completion}

The short-length calibration controls the normal approach to a stratum. For the completion problem we also need length functions on the surviving components.

\begin{lemma}[Bounded-length calibration]\label{lem:bounded-length-calibration}
There is a universal constant \(C_0\) such that, for every simple closed curve \(\alpha\), its length differential satisfies
\begin{equation}\label{eq:bounded-length-cometric}
\|d\ell_\alpha\|_{\gamma,*}\le C_0e^L
\qquad\text{whenever }\ell_\alpha\le L.
\end{equation}
Consequently \(f_{\alpha,L}:=\min\{\ell_\alpha,L\}\) is \(C_0e^L\)-Lipschitz for \(d_\gamma\), on all of \(\Teich\).
\end{lemma}

\begin{proof}
Conjugate the uniformizing group in the upper half-plane so that the axis of \(\alpha\) is \(i\mathbb R_+\). Gardiner's variation formula \cite{gardiner1975schiffer}, in the normalization recorded in \cite[\S2.3]{wolpert2008behavior}, represents \(d\ell_\alpha\) by the real part of the quadratic differential
\[
\kappa_\alpha=\frac2\pi\sum_{B\in\langle A\rangle\backslash\Gamma}
B^*\left(\frac{dz}{z}\right)^2.
\]
The hyperbolic pointwise norm of a summand at \(p\) is \(\operatorname{sech}^2d(p,B^{-1}(i\mathbb R_+))\). The lifts of a simple geodesic have disjoint collars of width
\[
w(\ell_\alpha)=\operatorname{arsinh}\!\left(\frac1{\sinh(\ell_\alpha/2)}\right);
\]
see \cite[Ch.~4]{buser2010geometry}. Put \(a=\min\{1,w(L)/2\}\). The radius-\(a\) disks centered at the nearest points of the distinct axes are disjoint. The number of these centers within distance \(R\) of \(p\) is at most \(Ca^{-2}e^R\), by comparison with the area of a hyperbolic disk. Summation over the shells \(m\le d<m+1\) gives
\[
\sum_B\operatorname{sech}^2d(p,B^{-1}(i\mathbb R_+))
\le Ca^{-2}\sum_{m\ge0}e^{-m}\le Ce^L.
\]
The last inequality follows from the collar-width formula, with an enlarged universal constant for \(L\le1\). This proves \eqref{eq:bounded-length-cometric}. Integration, with the absolute-continuity chain rule at the truncation level, proves the Lipschitz assertion.
\end{proof}

\begin{lemma}[Short transversals and marked compactness]\label{lem:marked-transversal-compactness}
Fix a marked pants decomposition \(P\), a submulticurve \(\Gamma\subset P\), and \(0<a\le B\). Suppose that
\[
\ell_\alpha\le B\quad(\alpha\in P),
\qquad a\le\ell_\alpha\quad(\alpha\in P\setminus\Gamma).
\]
Choose seam origins for the real Fenchel--Nielsen twists \(s_\alpha\), measured in units of length, so that a positive Dehn twist adds \(\ell_\alpha\). There are reference transversals \(\beta_\alpha\), disjoint from \(P\setminus\{\alpha\}\), and constants \(L,C\) depending only on \(a,B\), with the following properties.
\begin{enumerate}
\item For each \(\alpha\in P\setminus\Gamma\), some integral twist \(T_\alpha^k\beta_\alpha\) has length at most \(L\). The bound is uniform when any adjacent boundary lengths tend to zero.
\item If \(m_\alpha=i(\alpha,\beta_\alpha)\in\{1,2\}\), then
\begin{equation}\label{eq:transversal-real-twist}
|s_\alpha-k\ell_\alpha|
\le \frac{\ell_{T_\alpha^k\beta_\alpha}}{m_\alpha}+C.
\end{equation}
\item Suppose \(X_j\) satisfy the preceding length bounds, \(\ell_\gamma(X_j)\to0\) for \(\gamma\in\Gamma\), and one fixed marked transversal \(T_\alpha^{k_\alpha}\beta_\alpha\) has bounded length for each \(\alpha\in P\setminus\Gamma\). Then the surviving marked Fenchel--Nielsen coordinates lie in a compact set, and \(X_j\) have a subsequence converging in augmented Teichm\"uller space to the \(\Gamma\)-stratum.
\end{enumerate}
\end{lemma}

\begin{proof}
Cut along all cuffs other than \(\alpha\). Its adjacent piece is a one-holed torus or a four-holed sphere, allowing cusp boundary components. A seam marking supplies a transversal with \(m_\alpha=1\) or \(2\), respectively. Write \(t=s_\alpha-k\ell_\alpha\). The elementary-move length formulas take the form
\begin{equation}\label{eq:transversal-length-formulas}
\cosh\frac{\ell_{T_\alpha^k\beta_\alpha}}2
=
\begin{cases}
U\cosh(t/2),&m_\alpha=1,\\
U\cosh t+V,&m_\alpha=2.
\end{cases}
\end{equation}
Here \(U>0\) and \(V\ge0\) depend continuously on the cuff and external boundary lengths, including zero external lengths; see \cite[Propositions~3.1 and 4.1]{alessandrini2012behaviour}. For example, in the torus case
\[
U=\frac{\sqrt{(\cosh\ell_\alpha+\cosh(\ell_0/2))/2}}
{\sinh(\ell_\alpha/2)}.
\]
The four-holed-sphere coefficient is a product of two positive square roots divided by \(\sinh^2(\ell_\alpha/2)\). Since \(\ell_\alpha\in[a,B]\) and the external lengths belong to \([0,B]\), compactness of this parameter box gives
\[
0<U_-\le U\le U_+,\qquad 0\le V\le V_+,
\]
with bounds depending only on \(a,B\). These bounds hold also when external cuffs of the cut piece are identified elsewhere in the surface.

Choose \(k\) so that \(|t|\le\ell_\alpha/2\le B/2\). Equation~\eqref{eq:transversal-length-formulas} gives a uniform upper bound \(L\) for the transversal length, proving (i). Conversely,
\[
\frac{U_-}{2}e^{m_\alpha|t|/2}
\le\cosh\frac{\ell_{T_\alpha^k\beta_\alpha}}2
\le e^{\ell_{T_\alpha^k\beta_\alpha}/2}.
\]
Taking logarithms proves (ii), with \(C=2\max\{0,\log(2/U_-)\}\). A different fixed seam convention changes \(s_\alpha\) by a bounded multiple of \(\ell_\alpha\) and enlarges \(C\).

For (iii), the integers \(k_\alpha\) are fixed. Equation~\eqref{eq:transversal-real-twist} therefore bounds each surviving real twist. The corresponding lengths lie in \([a,B]\), so these coordinates range over a compact box. Extract a convergent subsequence in that box. The other lengths tend to zero, and their twists are forgotten in the augmented topology \cite{wolpertgeometry}. The subsequence consequently converges with the marking on every surviving component retained.
\end{proof}

The use of pants curves and transversals is shared with the earthquake-metric completion argument of \cite[Lemma~10.15]{huang2024earthquake}. Here the global calibration of truncated length functions supplies the Cauchy control of those fixed marked curves.

\begin{theorem}[Completion of the creation metric]\label{thm:creation-completion}
The identity of \(\Teich\) extends to a mapping-class-equivariant homeomorphism
\begin{equation}\label{eq:creation-completion}
\iota:\overline\Teich^{\WP}\xrightarrow{\ \cong\ }\widehat\Teich_\gamma.
\end{equation}
A sequence in \(\Teich\) is \(d_\gamma\)-Cauchy if and only if it is \(d_{\WP}\)-Cauchy, with the marked identification of their limits. Under \eqref{eq:creation-completion}, the boundary points are precisely marked noded surfaces; the marking on each surviving component is retained. In particular, \(D_\Gamma\) is the closed \(\Gamma\)-pinching stratum.
\end{theorem}

\begin{proof}
The extension \(\iota\) is Lipschitz by \eqref{eq:gamma-WP-distance-domination}. Each \(f_{\alpha,L}\) extends to the \(\gamma\)-completion by \cref{lem:bounded-length-calibration}. On augmented Teichm\"uller space its value is \(\min\{\ell_\alpha,L\}\), where the extended length is zero at a node, is the ordinary length on a surviving component, and is infinite if \(\alpha\) crosses a node. This is the length-function description of the augmented topology; see \cite{wolpertgeometry,wolpert2008behavior}.

These functions separate augmented points. Different node sets are separated by a curve belonging to one set and not the other. With the node set fixed, the marked length spectra on the components determine their Fenchel--Nielsen coordinates: pants lengths determine the pieces, and lengths of transverse curves and their twists determine the gluing parameters \cite[Chs.~3 and 6]{buser2010geometry}. Thus \(\iota\) is injective.

Let \((X_j)\) be \(d_\gamma\)-Cauchy. Fix \(\varepsilon>0\) below the collar threshold, so that all curves of length less than \(2\varepsilon\) are disjoint. Fix a genus-dependent pants bound \(B_g\) for extending such a short multicurve \cite{buser2010geometry,grosse2019holomorphic}. Choose \(N\) so that the diameter of the tail \(j\ge N\) is less than
\[
\min\left\{\frac{\varepsilon}{2C_0e^{2\varepsilon}},
\frac{1}{C_0e^{B_g+1}}\right\},
\]
and let \(\Gamma\) be the curves shorter than \(\varepsilon\) on \(X_N\). The Lipschitz bounds for \(f_{\alpha,\varepsilon}\) and \(f_{\alpha,2\varepsilon}\) are independent of \(\alpha\), so this choice gives, simultaneously for every curve,
\begin{equation}\label{eq:completion-fixed-short-curves}
\ell_\alpha(X_j)\ge\varepsilon/2\quad(\alpha\notin\Gamma),
\qquad
\ell_\alpha(X_j)<2\varepsilon\quad(\alpha\in\Gamma),
\qquad j\ge N.
\end{equation}
For example, a drop from at least \(\varepsilon\) to less than \(\varepsilon/2\) would require \(d_\gamma\)-distance at least \(\varepsilon/(2C_0e^\varepsilon)\). Hence only the fixed finite multicurve \(\Gamma\) can degenerate along the tail.

Extend \(\Gamma\) to a pants decomposition \(P\) with lengths at \(X_N\) bounded by \(B_g\). The second diameter bound and \(f_{\alpha,B_g+1}\) ensure that every pants length stays below \(B_g+1\). The lengths in \(\Gamma\) converge, since they agree with their truncations at \(2\varepsilon\). Write \(\Gamma_0\subset\Gamma\) for those with zero limit. All lengths in \(P\setminus\Gamma_0\) then lie in \([a,B_g+1]\) on a late tail, for some \(a>0\).

Apply \cref{lem:marked-transversal-compactness}(i) with \(B=B_g+1\). Fix \(L=L(a,B_g+1)\), and choose a late term \(X_{N'}\) so that the subsequent Cauchy diameter is less than \((C_0e^{L+1})^{-1}\). At \(X_{N'}\), select a transversal \(\beta_\alpha'=T_\alpha^{k_\alpha}\beta_\alpha\) of length at most \(L\) for each cuff outside \(\Gamma_0\). These are now fixed marked curves. The Lipschitz bound for \(f_{\beta_\alpha',L+1}\) gives
\[
\ell_{\beta_\alpha'}(X_j)<L+1,\qquad j\ge N'.
\]
Indeed, reaching the truncation level would change that function by at least one. Part (ii) bounds the surviving real twists with the integers \(k_\alpha\) fixed. Part (iii) then gives a subsequence converging in \(d_{\WP}\) to a marked noded surface \(Z\), with precisely the curves in \(\Gamma_0\) pinched.

Every subsequence has a further subsequence of this form, using the same pants and transverse curves. All resulting augmented limits have the same values of every \(f_{\alpha,L}\), since these functions are Cauchy on the original sequence. Separation gives the same limit \(Z\). It follows that the entire sequence converges in \(d_{\WP}\). Its \(\gamma\)-limit is \(\iota(Z)\), proving surjectivity and the Cauchy equivalence. The reverse implication in that equivalence follows directly from \eqref{eq:gamma-WP-distance-domination}.

Finally, if \(\iota(Z_j)\to\iota(Z)\), choose interior points \(Y_j\) with \(d_{\WP}(Y_j,Z_j)<1/j\). They are \(d_\gamma\)-Cauchy and hence converge in \(d_{\WP}\) by the preceding argument. Injectivity identifies their limit with \(Z\), so \(Z_j\to Z\). This proves continuity of the inverse. Equivariance follows from the invariant cotangent norm defining \(\gamma\).
\end{proof}

\begin{corollary}[Completion of the dual \(L^p\) metrics]\label{cor:Lp-completion}
For \(2\le p\le\infty\), let \(\gamma_{p,X}\) be the tangent norm dual to the hyperbolic \(L^p\) cotangent norm, and let \(d_p\) be its length metric. All \(d_p\) have the same Cauchy sequences and the same marked completion, canonically homeomorphic to augmented Teichm\"uller space.
\end{corollary}

\begin{proof}
For \(A=4\pi(g-1)\), H\"older's inequality gives
\[
A^{-1/p}\gamma_X(u)\le\gamma_{p,X}(u)
\le A^{1/2-1/p}\|u\|_{\WP},
\qquad
A^{-1/p}d_\gamma\le d_p\le A^{1/2-1/p}d_{\WP}.
\]
Theorem~\ref{thm:creation-completion} identifies the Cauchy sequences at the two ends. The same inequalities give continuity in both directions on the completions.
\end{proof}

\begin{theorem}[Linear distance to pinching strata]\label{thm:linear-pinching-distance}
Let \(\Gamma=\{\alpha_1,\ldots,\alpha_k\}\) be a nonempty multicurve. In a Fenchel--Nielsen family whose remaining data stay in a compact set away from further degeneration, put \(S(X)=\sum_{i=1}^k\ell_{\alpha_i}(X)\). As all these lengths tend to zero,
\begin{equation}\label{eq:linear-pinching-distance}
d_\gamma(X,D_\Gamma)=\frac\pi2S(X)+O(S(X)^2).
\end{equation}
In the Weil--Petersson normalization of this paper, it follows that
\begin{equation}\label{eq:pinching-distance-WP-square}
\frac{d_\gamma(X,D_\Gamma)}{d_{\WP}(X,\overline\Teich_\Gamma^{\WP})^2}
\longrightarrow\frac14.
\end{equation}
At fixed lengths and other coordinates, angular twisting satisfies
\begin{equation}\label{eq:angular-twist-distance}
d_\gamma(X(\theta),X(\theta'))
\le\frac1{4\pi}\sum_i\ell_i^2|\theta_i'-\theta_i|.
\end{equation}
\end{theorem}

\begin{proof}
We obtain a lower bound for arbitrary joining paths. Fix \(0<\ell_0<\ell_g\), and define on all of Teichm\"uller space
\[
S_0(X):=\sum_i\min\{\ell_{\alpha_i}(X),\ell_0\},
\qquad
F(X):=\int_0^{S_0(X)}\frac{ds}{2/\pi+C_gs^2}.
\]
At differentiability points, only curves shorter than \(\ell_0\) contribute to \(dS_0\). They are disjoint, and \cref{lem:short-length-calibration} bounds its dual norm by
\[
\|dS_0\|_{\gamma,*}
\le\frac2\pi+C_g\sum_{\ell_{\alpha_i}<\ell_0}\ell_{\alpha_i}^2
\le\frac2\pi+C_gS_0^2.
\]
Thus \(F\) is 1-Lipschitz for \(d_\gamma\), by integration along piecewise smooth paths and the absolute-continuity chain rule at level crossings. It extends to \(\widehat\Teich_\gamma\) and vanishes on \(D_\Gamma\). For all small pinching lengths,
\begin{equation}\label{eq:pinching-distance-global-lower}
d_\gamma(X,D_\Gamma)\ge F(X)
=\frac\pi2S(X)-O_g(S(X)^3).
\end{equation}

For the upper bound, decrease \(\ell_i\) along \(\ell_i(t)=t\ell_i(X)\), keeping all other coordinates fixed, with \(0<t\le1\). This path tends in augmented Teichm\"uller space to a point of the \(\Gamma\)-stratum. The central-circle coefficient bound and \eqref{eq:length-principal-pairing} give
\[
\gamma_X\!\left(\sum_i s_i\partial_{\ell_i}\right)
\le\frac\pi2\sum_i|s_i|+C\sum_i\ell_i|s_i|.
\]
Along the pinching path the constant is uniform under the stated compactness assumption. Its \(\gamma\)-length is therefore at most
\[
\int_0^1\left(\frac\pi2\sum_i\ell_i(X)
+Ct\sum_i\ell_i(X)^2\right)dt
=\frac\pi2S(X)+O(S(X)^2).
\]
Together with \eqref{eq:pinching-distance-global-lower}, this proves \eqref{eq:linear-pinching-distance}. The Weil--Petersson stratum-distance expansion \cite[Corollary~21]{wolpertgeometry} is \(d_{\WP}(X,\overline\Teich_\Gamma^{\WP})=\sqrt{2\pi S(X)}+O(S(X)^2)\); it gives \eqref{eq:pinching-distance-WP-square}. Finally, integrate \eqref{eq:angular-twist-global-upper} along successive angular coordinate segments to obtain \eqref{eq:angular-twist-distance}.
\end{proof}

The radial and angular estimates describe finite-distance pinching in the creation metric. The lower bound uses a global Lipschitz function, so excursions outside the coordinate family cannot change its leading radial coefficient.

\subsection{Exact laws along a degenerating family}

\begin{corollary}[Realization of the limiting additive cost]\label{cor:degenerating-exact-laws}
Let \(u_j\) have the coordinates \((u_0,c)\ne0\) in \eqref{eq:degeneration-tangent-identification}, and put
\[
G_0=\gamma_\Sigma(u_0)+\frac1{4\pi^2}\sum_i|c_i|,
\qquad \sigma_j=\min\{1,\ell_{\mathrm{sys}}(X_j)\}.
\]
For any prescribed positive scales satisfying
\[
r_j/\sigma_j^4\to0,
\qquad \delta_j/(r_j^2\sigma_j^2)\to0,
\]
for all sufficiently large \(j\) there are equal-weight pairs in the fibers over \(X_j\), of vertical radius at most \(r_j\), with exact output \(\tau_{X_j}^{-1}(\delta_ju_j)\), such that
\begin{equation}\label{eq:degenerating-exact-costs}
\frac{\|v_j\|_{2,X_j}^2}{\delta_j}\to G_0,
\qquad
\frac{\cU_{2,X_j}(\Pi_j)^2}{4\delta_j}\to G_0,
\qquad
\frac{\Delta_{X_j}(\Pi_j)}{4\delta_j^2\|u_j\|_{\WP,X_j}^2}\to1.
\end{equation}
\end{corollary}

\begin{proof}
By \cref{thm:additive-cost-degeneration}, \(\gamma_{X_j}(u_j)\to G_0>0\). Apply \cref{cor:prescribed-degenerating-scales} to \(u_j/\gamma_{X_j}(u_j)\) with displacement parameter \(\delta_j\gamma_{X_j}(u_j)\). The hypotheses are unchanged up to bounded factors, and \eqref{eq:uniform-joint-sharp-costs} gives the three limits.
\end{proof}

\section{Further questions}\label{sec:discussion}

The Hessians of nondegenerate cotangent maxima give the explicit lower coefficient \eqref{eq:explicit-excess-liminf}. Determining the full quadratic excess coefficient requires minimizing the other two terms in the calibration deficit together with the vector output constraint. The exposed-face description \eqref{eq:finite-contact-face} specifies the possible contact masses, but a genericity statement for supporting cotangent differentials alone does not determine their prevalence among tangent directions.

The finite-radius limit is locally uniform on the interior of the feasible body. Boundary outputs may require a different recovery construction, since the positive-weight correction used in \cref{thm:finite-radius-limit} loses its interior slack. It also remains to determine whether the constrained single-direction value \(J_X^{(2)}\) agrees with its convexification \(J_X\).

The boundary norm has linear radial cost and angular cost of order \(\ell^2\). A distance model for simultaneous length and twist variation would require estimates for actual marked Fenchel--Nielsen velocities and lower bounds controlling all joining paths. The systole powers in the uniform realization theorem give sufficient ranges. Collar-adapted fields may enlarge those ranges for directions concentrated near a short geodesic.

\appendix

\section{Nonlinear Bers and fiber estimates}\label{app:bers-nonlinear}\label{sec:l2tame}

The exact-output estimates use two hyperbolic \(L^2\) factors in each nonlinear remainder. We prove them by moving the exterior point of Bers evaluation to infinity.

\begin{proposition}[Nonlinear \(L^2\)-tameness of the Bers projection]\label{prop:bers-nonlinear-l2}
There exist constants \(\delta_X>0\) and \(C_X<\infty\) such that, whenever
\[
  \|\mu\|_\infty\le \delta_X,
  \qquad
  h\in\Bspace(X),
\]
one has
\begin{align}
  \|\mathfrak P(\mu)-D\mathfrak P_0[\mu]\|_B
  &\le C_X\|\mu\|_2^2,\label{eq:bers-nonlinear-zero}\\
  \|D\mathfrak P_\mu[h]-D\mathfrak P_0[h]\|_B
  &\le C_X\|\mu\|_2\|h\|_2.\label{eq:bers-nonlinear-first}
\end{align}
The remainders after the quadratic term satisfy
\begin{align}
\|\mathfrak P(\mu)-D\mathfrak P_0[\mu]-\tfrac12D^2\mathfrak P_0[\mu,\mu]\|_B
&\le C_X\|\mu\|_\infty\|\mu\|_2^2,\label{eq:bers-cubic-zero}\\
\|D\mathfrak P_\mu[h]-D\mathfrak P_0[h]-D^2\mathfrak P_0[\mu,h]\|_B
&\le C_X\|\mu\|_\infty\|\mu\|_2\|h\|_2.\label{eq:bers-cubic-first}
\end{align}
Here \(\mathfrak P\) is the Bers projection associated with the fixed uniformization used above. One can take \(\delta_X=a_0>0\) universal and \(C_X=C\min\{1,\ell_{\mathrm{sys}}(X)\}^{-2}\).
\end{proposition}

\begin{lemma}[Uniform lift bound]\label{lem:uniform-lift-bound}
For \(\sigma_X=\min\{1,\ell_{\mathrm{sys}}(X)\}\), every disk uniformization and every lifted Beltrami differential satisfy
\begin{equation}\label{eq:uniform-euclidean-lift}
\|\nu^\sharp\|_{2,\mathrm e}^2
\le C\sigma_X^{-2}\|\nu\|_{2,X}^2.
\end{equation}
The constant is independent of the position of the disk origin.
\end{lemma}

\begin{proof}
Unfold the Euclidean integral over a fundamental set, using the hyperbolic area on the quotient. The resulting kernel at a lift \(z\) is
\[
K(z)=\frac14\sum_{\gamma\in\Gamma}(1-|\gamma z|^2)^2.
\]
The orbit points are at least \(\sigma_X\) apart in the hyperbolic metric. Their radius-\(\sigma_X/4\) disks are disjoint. On each disk the function \((1-|\zeta|^2)^2/4\) is comparable to its value at the center, with a universal constant. Summing the disk integrals gives
\[
K(z)\le C\sigma_X^{-2}\int_\D
\frac{(1-|\zeta|^2)^2}{4}\,dA_\D(\zeta)
=C\pi\sigma_X^{-2}.
\]
The orbit separation holds after any change of disk origin. Unfolding now proves \eqref{eq:uniform-euclidean-lift}.
\end{proof}

\begin{proof}[Proof of Proposition~\ref{prop:bers-nonlinear-l2}]
Lift the coefficients to \(\D\), extend them by zero, and use Euclidean \(L^p\) norms until the last step. Let
\[
\mathscr T f(z)=-\frac1\pi\int_{\mathbb C}\frac{f(\zeta)}{\zeta-z}\,dA_{\mathrm e}(\zeta),
\qquad
\mathscr Hf(z)=-\frac1\pi\operatorname{p.v.}\!\int_{\mathbb C}\frac{f(\zeta)}{(\zeta-z)^2}\,dA_{\mathrm e}(\zeta).
\]
These are the Cauchy and Beurling normalizations used in \cite[\S2, Theorems~2--3]{guo2003geodesic}. The principal solution is \(F_\mu=z+\mathscr T\omega\), with \(\omega\) as below. Put
\[
A_\mu=M_\mu\mathscr H,\qquad R_\mu=(I-A_\mu)^{-1},
\qquad \omega=R_\mu\mu.
\]
Choose a universal \(a_0>0\) with
\[
a_0\max\{\Lambda_2,\Lambda_4,\Lambda_{4/3}\}<\tfrac12,
\]
where \(\Lambda_p\) is the \(L^p\) Beurling norm and \(\Lambda_2=1\). The three resolvents are then bounded by \(2\) \cite[\S4.5.2]{astala2008elliptic}. Write \(a=\|\mu\|_\infty\le a_0\). The resolvent equation gives
\begin{align*}
\|\omega\|_{2,\mathrm e}&\le2\|\mu\|_{2,\mathrm e},\\
\|\omega-\mu\|_{1,\mathrm e}&\le2\|\mu\|_{2,\mathrm e}^2,
&
\|\omega-\mu\|_{2,\mathrm e}&\le2a\|\mu\|_{2,\mathrm e}.
\end{align*}
The cubic coefficient remainder is
\[
\kappa=\omega-\mu-\mu\mathscr H\mu
=\mu\mathscr H(\omega-\mu),
\qquad
\|\kappa\|_{1,\mathrm e}\le2a\|\mu\|_{2,\mathrm e}^2.
\]

Differentiating the resolvent equation yields
\[
\dot\omega-h
=R_\mu(\mu\mathscr Hh)+R_\mu(h\mathscr H\omega),
\qquad \dot\omega=D\omega_\mu[h].
\]
For the first series, place the outermost \(\mu\) in \(L^2\) and all remaining coefficient factors in \(L^\infty\):
\[
\|A_\mu^m(\mu\mathscr Hh)\|_{1,\mathrm e}
\le a^m\|\mu\|_{2,\mathrm e}\|h\|_{2,\mathrm e},
\qquad m\ge0.
\]
The zeroth term of the second series satisfies
\(\|h\mathscr H\omega\|_{1,\mathrm e}
\le2\|h\|_{2,\mathrm e}\|\mu\|_{2,\mathrm e}\).
For every higher term, H\"older's inequality with exponents \(4\) and \(4/3\) gives
\begin{align}
\|A_\mu^m(h\mathscr H\omega)\|_{1,\mathrm e}
&\le \Lambda_{4/3}\|\mu\|_{4,\mathrm e}
(a\Lambda_{4/3})^{m-1}\|h\mathscr H\omega\|_{4/3,\mathrm e}\notag\\
&\le C(a\Lambda_{4/3})^{m-1}
\|\mu\|_{4,\mathrm e}\|h\|_{2,\mathrm e}\|\omega\|_{4,\mathrm e}\notag\\
&\le Ca(a\Lambda_{4/3})^{m-1}
\|\mu\|_{2,\mathrm e}\|h\|_{2,\mathrm e},
\qquad m\ge1.\label{eq:bers-mixed-tail}
\end{align}
The last step uses \(\|\omega\|_{4,\mathrm e}\le2\|\mu\|_{4,\mathrm e}\) and
\(\|\mu\|_{4,\mathrm e}^2\le a\|\mu\|_{2,\mathrm e}\). Thus
\[
\|\dot\omega-h\|_{1,\mathrm e}
\le C\|\mu\|_{2,\mathrm e}\|h\|_{2,\mathrm e}.
\]

Subtracting the two quadratic variation terms gives exactly three tails:
\begin{equation}\label{eq:bers-cubic-derivative-tails}
D\kappa_\mu[h]
=\underbrace{\sum_{m\ge1}A_\mu^m(\mu\mathscr Hh)}_{T_1}
+\underbrace{h\mathscr H(\omega-\mu)}_{T_2}
+\underbrace{\sum_{m\ge1}A_\mu^m(h\mathscr H\omega)}_{T_3}.
\end{equation}
\begin{samepage}
Their bounds are
\begin{align*}
\|T_1\|_{1,\mathrm e}
&\le\frac{a}{1-a}\|\mu\|_{2,\mathrm e}\|h\|_{2,\mathrm e},\\
\|T_2\|_{1,\mathrm e}
&\le2a\|\mu\|_{2,\mathrm e}\|h\|_{2,\mathrm e},\\
\|T_3\|_{1,\mathrm e}
&\le\frac{Ca}{1-a\Lambda_{4/3}}
\|\mu\|_{2,\mathrm e}\|h\|_{2,\mathrm e}.
\end{align*}
\end{samepage}
The first is the \(L^2\) geometric tail; the second uses the displayed \(L^2\) bound for \(\omega-\mu\); the third follows by summing \eqref{eq:bers-mixed-tail}. In particular,
\begin{equation}\label{eq:bers-cubic-derivative-L1}
\|D\kappa_\mu[h]\|_{1,\mathrm e}
\le Ca\|\mu\|_{2,\mathrm e}\|h\|_{2,\mathrm e}.
\end{equation}

We now move an exterior evaluation point \(p\) to infinity. For finite \(|p|>1\), set \(q=1/\bar p\) and
\[
A_p(z)=\frac{z-q}{1-\bar qz},\qquad B_p=A_p^{-1}.
\]
This disk automorphism satisfies \(A_p(p)=\infty\); use the identity when \(p=\infty\). In the coordinate \(w=A_p(z)\), the coefficient and its variation become
\begin{equation}\label{eq:bers-mobius-coefficients}
\mu_p(w)=\mu(B_p(w))\frac{\overline{B_p'(w)}}{B_p'(w)},
\qquad
h_p(w)=h(B_p(w))\frac{\overline{B_p'(w)}}{B_p'(w)}.
\end{equation}
They remain supported in \(\D\), with unchanged supremum norms. After postcomposition by a M\"obius normalization, \(F_\mu\circ B_p\) is the normalized solution for \(\mu_p\). The Schwarzian chain rule gives
\begin{equation}\label{eq:bers-mobius-schwarzian}
\mathfrak P(\mu_p)(w)
=\mathfrak P(\mu)(B_p(w))B_p'(w)^2.
\end{equation}
Since \(p\) is fixed independently of \(\mu\), differentiating this identity gives the same pullback formula for \(D\mathfrak P_\mu[h]\) and for \(D^2\mathfrak P_0[\mu,h]\). It therefore applies to both remainders and their first variations. The Bers weight obeys the exact identity
\begin{equation}\label{eq:bers-mobius-weight}
(|w|^2-1)^2|\phi(B_p(w))B_p'(w)^2|
=(|B_p(w)|^2-1)^2|\phi(B_p(w))|.
\end{equation}
Thus the weighted value at \(p\) is the absolute value of the \(w^{-4}\) coefficient at infinity in the new uniformization.

For its principal solution the Cauchy representation gives
\[
F_{\mu_p}(w)=w+\frac{c_p}{w}+O(w^{-2}),
\qquad c_p=\frac1\pi\int_\D\omega_p(\zeta)\,d\xi\,d\eta,
\qquad \omega_p=R_{\mu_p}\mu_p.
\]
The product term in the Schwarzian is \(O(w^{-6})\), so
\begin{equation}\label{eq:bers-at-infinity}
\lim_{w\to\infty}w^4\mathfrak P(\mu_p)(w)
=-\frac6\pi\int_\D\omega_p(\zeta)\,d\xi\,d\eta.
\end{equation}
The linear and quadratic coefficients replace \(\omega_p\) by \(\mu_p\) and \(\mu_p\mathscr H\mu_p\), respectively. The four \(L^1\) estimates above consequently control the quadratic remainder, its derivative, the cubic remainder, and its derivative at \(p\).

The conjugate disk uniformization represents the same quotient, hence
\(\|\mu_p\|_{2,X}=\|\mu\|_{2,X}\) and
\(\|h_p\|_{2,X}=\|h\|_{2,X}\). Apply \cref{lem:uniform-lift-bound} to both factors in each new uniformization, and take the supremum over \(p\). All four estimates in \cref{prop:bers-nonlinear-l2} follow with
\begin{equation}\label{eq:explicit-Bers-constant}
C_X=C\sigma_X^{-2},\qquad \delta_X=a_0.
\end{equation}
\end{proof}

\subsection{Nonlinear control of the local fiber}

The Bers estimate also controls the graph \(\Psi\) throughout a fixed small \(L^\infty\)-neighborhood. This is the input needed to compare all local-fiber competitors in the quartic restoration problem.

\begin{proposition}[Nonlinear \(L^2\)-tame fiber graph]\label{prop:fiber-graph-l2}
After shrinking the graph neighborhood in Proposition~\ref{prop:fiber-graph}, there is a constant \(C_X\) such that
\begin{align}
  \|\Psi(v)\|_2
  &\le C_X\|v\|_2^2,\label{eq:Psi-L2-quadratic}\\
  \|D\Psi_v[h]\|_2
  &\le C_X\|v\|_2\|h\|_2,
  \qquad h\in V.\label{eq:DPsi-L2}
\end{align}
For the quadratic graph term \(B=D^2\Psi_0\),
\begin{align}
\|\Psi(v)-\tfrac12B(v,v)\|_2
&\le C_X\|v\|_\infty\|v\|_2^2,\label{eq:Psi-L2-cubic}\\
\|D\Psi_v[h]-B(v,h)\|_2
&\le C_X\|v\|_\infty\|v\|_2\|h\|_2.\label{eq:DPsi-L2-cubic}
\end{align}
Consequently, whenever the line segment from \(u\) to \(v\) remains in this neighborhood,
\begin{equation}\label{eq:Psi-L2-Lipschitz}
  \|\Psi(v)-\Psi(u)\|_2
  \le C_X(\|u\|_2+\|v\|_2)\|v-u\|_2.
\end{equation}
\end{proposition}

\begin{proof}
Set \(\mu=\Phi(v)=v+\Psi(v)\).  Since \(\mathfrak P(\mu)=0\), \(D\mathfrak P_0[v]=0\), and \(D\mathfrak P_0|_H\) is an isomorphism onto the tangent space of the Bers image, \eqref{eq:bers-nonlinear-zero} gives
\[
  \|\Psi(v)\|_2
  \le C_X\|\mu\|_2^2
  =C_X\bigl(\|v\|_2^2+\|\Psi(v)\|_2^2\bigr).
\]
The last equality uses \(V\perp_{L^2}H\).  For a sufficiently small graph neighborhood, the quadratic term in \(\|\Psi(v)\|_2\) is absorbed, proving \eqref{eq:Psi-L2-quadratic}.

Differentiate the identity \(\mathfrak P(\Phi(v))=0\) in the direction \(h\in V\):
\begin{equation}\label{eq:fiber-diff-bers}
  D\mathfrak P_\mu[h+D\Psi_v[h]]=0.
\end{equation}
Because \(D\mathfrak P_0[h]=0\), Proposition~\ref{prop:bers-nonlinear-l2} and \eqref{eq:Psi-L2-quadratic} give
\[
  \|D\mathfrak P_\mu[h]\|_B
  \le C_X\|\mu\|_2\|h\|_2
  \le C_X\|v\|_2\|h\|_2.
\]
The finite-dimensional operator \(D\mathfrak P_\mu|_H\) remains uniformly invertible for \(\mu\) in a sufficiently small \(L^\infty\)-neighborhood of the origin.  Since \(D\Psi_v[h]\in H\), equation~\eqref{eq:fiber-diff-bers} proves \eqref{eq:DPsi-L2}.  Integrating \eqref{eq:DPsi-L2} along the segment \(u+t(v-u)\) gives \eqref{eq:Psi-L2-Lipschitz}.

The first-variation estimate at \(\mu=ta\), divided by \(|t|\) and passed to the limit, gives
\[
\|D^2\mathfrak P_0[a,b]\|_B\le C_X\|a\|_2\|b\|_2
\]
for all bounded Beltrami differentials \(a,b\). Put \(\psi=\Psi(v)\). The fiber identity and \(D\mathfrak P_0[B(v,v)]=-D^2\mathfrak P_0[v,v]\) imply
\[
D\mathfrak P_0[\psi-\tfrac12B(v,v)]
=-D^2\mathfrak P_0[v,\psi]-\tfrac12D^2\mathfrak P_0[\psi,\psi]-R_3(\mu),
\]
where \(\|R_3(\mu)\|_B\le C_X\|\mu\|_\infty\|\mu\|_2^2\). Since \(\|\mu\|_2\le C_X\|v\|_2\), \(\|\mu\|_\infty\le C_X\|v\|_\infty\), and \(\|v\|_2\le C_X\|v\|_\infty\), inversion of \(D\mathfrak P_0|_H\) proves \eqref{eq:Psi-L2-cubic}.

For the derivative, substitute
\[
D\mathfrak P_\mu[k]
=D\mathfrak P_0[k]+D^2\mathfrak P_0[\mu,k]+R'_3(\mu)[k]
\]
into \eqref{eq:fiber-diff-bers}, with \(k=h+D\Psi_v[h]\). The identity \(D\mathfrak P_0[B(v,h)]=-D^2\mathfrak P_0[v,h]\) leaves
\[
D\mathfrak P_0[D\Psi_v[h]-B(v,h)]
=-D^2\mathfrak P_0[\psi,h]-D^2\mathfrak P_0[\mu,D\Psi_v[h]]-R'_3(\mu)[k].
\]
The bilinear estimate, \eqref{eq:DPsi-L2}, and \eqref{eq:bers-cubic-first} bound the right side by \(C_X\|v\|_\infty\|v\|_2\|h\|_2\). Inverting on \(H\) proves \eqref{eq:DPsi-L2-cubic}.
\end{proof}

\section{Uniform graph and localization estimates}\label{app:uniform-proofs}

\subsection{Uniform graph bounds}

We prove \cref{prop:systole-nonlinear}.

\begin{proof}
Write \(L_X=D\mathfrak P_0\). In disk coordinates the Ahlfors--Weill reflected inverse \cite{ahlfors1962uniqueness,nag1987conformally} is
\[
(S_X\phi)(z)=-\frac{(1-|z|^2)^2}{2\bar z^4}\phi(1/\bar z),
\qquad z\in\D,
\]
with its value at zero defined by continuation. The automorphy of \(\phi\) makes \(S_X\phi\) a harmonic Beltrami differential on \(X\). The Bers reproducing formula gives \(L_XS_X=\id\) \cite{nag1987conformally}. Since
\((|1/\bar z|^2-1)^2=(1-|z|^2)^2/|z|^4\),
\[
\|S_X\phi\|_\infty=\tfrac12\|\phi\|_B,
\qquad P_H=S_XL_X,
\qquad \|P_H\|_{\infty\to\infty}\le3.
\]
The last estimate follows from \(\|L_X\|_{\infty\to B}\le6\), obtained by differentiating the universal Bers bound \(\|\mathfrak P(\mu)\|_B\le6\) at zero.

Fix \(\eta=1/2\) and put \(N_X(\mu)=\mathfrak P(\mu)-L_X\mu\). Holomorphic Cauchy estimates on a larger universal coefficient ball give a universal \(C_0\ge1\) such that, for \(\|\mu\|_\infty\le\eta\),
\[
\|N_X(\mu)\|_B\le C_0\|\mu\|_\infty^2,
\qquad
\|DN_X(\mu)[k]\|_B\le C_0\|\mu\|_\infty\|k\|_\infty.
\]
For \(v\in V\), put \(a=\|v\|_\infty\). Since \(L_Xv=0\), the fiber equation is equivalent to the fixed-point equation
\[
\psi=T_{X,v}(\psi):=-S_XN_X(v+\psi),
\qquad \psi\in H.
\]
Choose \(0<\varepsilon_0\le\min\{\eta/2,(4C_0)^{-1}\}\). If \(a\le\varepsilon_0\) and \(\|\psi\|_\infty\le2C_0a^2\), then \(\|v+\psi\|_\infty\le2a\le\eta\), and
\[
\|T_{X,v}(\psi)\|_\infty\le2C_0a^2,
\qquad \operatorname{Lip}(T_{X,v})\le C_0a\le\tfrac14.
\]
Contraction gives a unique fixed point \(\Psi(v)\) in this harmonic ball, holomorphic in \(v\), with \(\|\Psi(v)\|_\infty\le2C_0\|v\|_\infty^2\). This proves the uniform graph radius and the first bound in \eqref{eq:uniform-infinity-graph}.

The disk midpoint has the explicit formula
\begin{equation}\label{eq:uniform-midpoint-formula}
w=\frac{b-a}{1-\bar a b},\qquad
t=\frac{w}{1+\sqrt{1-|w|^2}},\qquad
\mathfrak m(a,b)=\frac{a+t}{1+\bar a t}.
\end{equation}
For \(|a|,|b|\le1/4\), its real derivatives through order four have universal bounds: the denominators and the square-root argument stay uniformly away from zero. The same bounds hold for the pointwise map on \(L^\infty\). On a smaller universal graph ball, Cauchy estimates on the contraction ball bound \(D^j\Psi\) for \(1\le j\le4\); the universal Bers bound on a larger coefficient ball bounds \(D^j\mathfrak P\) there. The chain rule therefore gives a universal \(C_4\) such that the harmonic lift
\[
\widehat G_X(v)
=S_X\mathfrak P\bigl(\mathfrak m(\Phi(v),\Phi(-v))\bigr),
\qquad D\pi_0[\widehat G_X(v)]=G_X(v),
\]
satisfies \(\|D^4\widehat G_X(v)\|_{\infty\to\infty}\le C_4\) on that ball. Its evenness removes the odd Taylor terms, and \cref{prop:B-S} identifies the quadratic term with the harmonic lift of \(q_X\). Taylor's integral remainder and \(\gamma_X(D\pi_0[h])\le4\pi(g-1)\|h\|_\infty\) prove the second bound in \eqref{eq:uniform-infinity-graph}. Fix the resulting universal \(\varepsilon_0\), reducing it also so that \(2\varepsilon_0\le a_0\), the coefficient radius of Appendix~\ref{app:bers-nonlinear}.

Appendix~\ref{app:bers-nonlinear} proves all four Bers estimates in \cref{prop:bers-nonlinear-l2} with constant \(C\sigma_X^{-2}\) on a universal coefficient ball. The reflected inverse has norm at most \(\sqrt{4\pi(g-1)}/2\) from the Bers norm to harmonic \(L^2\). Fix \(D_g\ge1\) to dominate the reflected quadratic and first-variation constants and the area bound coming from the graph supremum estimate. Before absorption,
\[
\|\Psi(v)\|_2\le D_g\sigma_X^{-2}
\bigl(\|v\|_2^2+\|\Psi(v)\|_2^2\bigr),
\qquad
\|\Psi(v)\|_2\le D_g\|v\|_\infty^2.
\]
Set \(E_g=\sqrt{4\pi(g-1)}+D_g\) and choose a preliminary graph threshold
\[
0<\rho_g^{(0)}\le
\min\{\varepsilon_0/8,1,(4D_g)^{-1},(4D_gE_g)^{-1}\}.
\]
For \(\|v\|_\infty\le\rho_g^{(0)}\sigma_X^2\),
\(D_g\sigma_X^{-2}\|\Psi(v)\|_2
\le D_g^2(\rho_g^{(0)})^2\le1/16\).
The quadratic graph term is therefore absorbed, giving
\[
\|\Psi(v)\|_2\le2D_g\sigma_X^{-2}\|v\|_2^2,
\qquad \|\mu\|_2\le2\|v\|_2,
\qquad \|\mu\|_2\le E_g\rho_g^{(0)}\sigma_X^2.
\]
The harmonic inverse is controlled directly in \(L^2\). For
\(A_{X,\mu}=S_X(D\mathfrak P_\mu-L_X)|_H\), the reflected first-variation estimate gives
\begin{equation}\label{eq:uniform-harmonic-inverse-absorption}
\|A_{X,\mu}\|_{2\to2}
\le D_g\sigma_X^{-2}\|\mu\|_2
\le D_gE_g\rho_g^{(0)}\le\tfrac14,
\qquad
\|(I_H+A_{X,\mu})^{-1}\|_{2\to2}\le\tfrac43.
\end{equation}
For \(J=D\Psi_v[k]\), differentiate the fiber equation and use \(L_Xk=0\):
\[
(I_H+A_{X,\mu})J=-S_X(D\mathfrak P_\mu-L_X)k,
\qquad
\|J\|_2\le\tfrac83D_g\sigma_X^{-2}\|v\|_2\|k\|_2.
\]
This proves the quadratic derivative bound.

Put \(\psi=\Psi(v)\) and let \(R_3\) be the cubic Bers remainder. The two graph identities are
\begin{align*}
\psi-\tfrac12B(v,v)
&=-S_X\bigl(D^2\mathfrak P_0[v,\psi]
+\tfrac12D^2\mathfrak P_0[\psi,\psi]+R_3(\mu)\bigr),\\
J-B(v,k)
&=-S_X\bigl(D^2\mathfrak P_0[\psi,k]
+D^2\mathfrak P_0[\mu,J]+DR_3(\mu)[k+J]\bigr).
\end{align*}
The reflected bilinear bound contributes \(D_g\sigma_X^{-2}\) and the bounds on \(\psi,J\) contribute another \(C_g\sigma_X^{-2}\). Using \(\|v\|_2\le\sqrt{4\pi(g-1)}\|v\|_\infty\) bounds the mixed terms by \(C_g\sigma_X^{-4}\|v\|_\infty\|v\|_2^2\) and \(C_g\sigma_X^{-4}\|v\|_\infty\|v\|_2\|k\|_2\), respectively. The pure cubic remainders have the smaller factor \(C_g\sigma_X^{-2}\). These prove the asserted cubic bounds.

An embedded disk of radius \(\sigma_X/4\), the mean-value estimate, and harmonic duality give
\begin{equation}\label{eq:systole-linear-L1}
\|P_H a\|_2\le C_g^{\mathrm{an}}\sigma_X^{-1}\|a\|_1.
\end{equation}
In fact a cotangent vector of unit \(L^2\)-norm has pointwise norm at most \(C\sigma_X^{-1}\) on that disk. For \(\mu=\Phi(v)\), orthogonality and the graph bounds give
\[
\E\|\mu\|_2^2\le(1+C_g^{\mathrm{an}}\sigma_X^{-4}r^2)M,
\qquad \|\E\mu\|_2\le C_g^{\mathrm{an}}\sigma_X^{-2}M.
\]
The pointwise barycenter estimate has a universal constant and yields
\[
\|\beta-\E\mu\|_1\le C_g^{\mathrm{an}}rM,
\qquad
\|\beta-\E\mu\|_2^2\le C_g^{\mathrm{an}}r^4M,
\qquad
\|\beta\|_2^2\le C_g^{\mathrm{an}}\sigma_X^{-4}r^2M.
\]
The mean graph remainder is bounded by \(C_g^{\mathrm{an}}\sigma_X^{-4}rM\) in \(L^2\). The linear barycenter correction is bounded by \(C_g^{\mathrm{an}}\sigma_X^{-1}rM\) after projection, using \eqref{eq:systole-linear-L1}. For the nonlinear Bers remainder, the two factors are explicitly
\[
\|\mathfrak P(\beta)-L_X\beta\|_B
\le C\sigma_X^{-2}\|\beta\|_2^2
\le C_g^{\mathrm{an}}\sigma_X^{-6}r^2M.
\]
The first factor \(\sigma_X^{-2}\) comes from the lift estimate in \cref{lem:uniform-lift-bound} applied to the two \(L^2\) inputs, and the remaining \(\sigma_X^{-4}\) comes from the preceding bound for \(\|\beta\|_2^2\). As \(r\le\rho_g\sigma_X^4\) and \(0<\sigma_X\le1\), this remainder is at most \(C_g^{\mathrm{an}}\sigma_X^{-4}rM\). The other two errors obey the same final bound. The comparison \(\gamma_X\le\sqrt{4\pi(g-1)}\|\cdot\|_{\WP}\) proves \eqref{eq:systole-law-error}. The Euclidean variance and the disk-metric comparison give \eqref{eq:systole-dispersion-error}.

For restoration, use \(b=P_V\beta\), \(h=P_H\beta\), \(\nu_0=\Phi(b)\), and \(n=h-\Psi(b)\), as in \cref{prop:uniform-exact-restoration}. Since \(P_V\E\mu=0\) and \(\|P_V\|_{\infty\to\infty}\le4\),
\[
\|b\|_\infty+\|b\|_2\le C_g^{\mathrm{an}}r^3,
\qquad \|\beta\|_2+\|n\|_2\le C_g^{\mathrm{an}}\sigma_X^{-2}r^2.
\]
Thus \(b\in U_X\). Integrating the Bers derivative along \(\nu_0+tn\) gives
\[
\|y-D\pi_0[n]\|_{\WP}
\le C_g^{\mathrm{an}}\sigma_X^{-4}r^2\|n\|_2.
\]
Absorption identifies \(n\) with the harmonic lift of \(y\) to relative error \(C_g^{\mathrm{an}}\sigma_X^{-4}r^2\). The loss at \(\nu_0\) gives the upper bound in \eqref{eq:systole-restoration-error}. A near-minimizer \(\Phi(w)\) satisfies, by the variance inequality,
\[
\|w-b\|_2\le C_g^{\mathrm{an}}\|y\|_{\WP},
\qquad
\|w\|_2\le C_g^{\mathrm{an}}(r^3+\|y\|_{\WP}),
\qquad \|y\|_{\WP}\le C_g^{\mathrm{an}}\sigma_X^{-2}r^2.
\]
The graph Lipschitz bound therefore yields
\[
\|\Psi(w)-\Psi(b)\|_2
\le C_g^{\mathrm{an}}\sigma_X^{-4}r^2\|y\|_{\WP}.
\]
Orthogonality and the variance inequality give the matching lower bound throughout the full graph \(\Phi(U_X)\).

All graph and law constants obtained above are bounded by a genus-dependent \(L_g^{\mathrm{an}}\), independent of \(X\). The absorptions in the law and restoration estimates require only that \(L_g^{\mathrm{an}}\sigma_X^{-4}r^2\le1/4\), together with the preliminary graph range. Fix
\[
0<\rho_g\le
\min\{\rho_g^{(0)},(4L_g^{\mathrm{an}})^{-1/2},1\},
\qquad C_g^{\mathrm{an}}\ge\max\{D_g,L_g^{\mathrm{an}}\}.
\]
Since \(r\le\rho_g\sigma_X^4\), these choices satisfy all the absorptions. The constants are now fixed for the rest of the section.
\end{proof}

\subsection{Disjoint localized directions}\label{app:uniform-basis-proof}

We prove \cref{lem:uniform-localized-realization}.

\begin{proof}
Set \(n=3g-3\). The response body in \cref{cor:geometric-response-body} represents \(u\) as a convex combination of at most \(2n+1\) evaluation atoms. Maximize the absolute determinant of \(n\) point-evaluation columns, using unit Beltrami frames at their centers. The maximum is positive because point evaluations span \(\QD(X)^*\). If \(e_i\) are maximizing columns and \(\phi_i\) are their dual quadratic differentials, replacement of the \(i\)-th column by any evaluation atom gives
\[
\|\phi_i\|_{\infty,h}\le1,\qquad
|\langle z,\phi_i\rangle|\le\gamma_X(z).
\]
Thus, for \(E(c)=\sum_i c_ie_i\),
\begin{equation}\label{eq:uniform-atom-basis-inverse}
\|E^{-1}z\|_\infty\le\gamma_X(z).
\end{equation}

On the universal cover, a quadratic differential of unit hyperbolic supremum norm has uniformly bounded first derivatives on a fixed-radius disk. Normalize a real bump displacement of radius \(h\) so that its vertical derivative \(a\) has \(\|a\|_2=1\). The local isotopy formula \eqref{eq:local-response}, followed by a phase rotation, gives
\begin{equation}\label{eq:uniform-bump-atom}
\|a\|_\infty\le C/h,\qquad
\gamma_X(q_X(a)-e_X(p,b))\le Ch.
\end{equation}
Its response is the weighted average of the local cotangent coefficient against \(|W_z|^2\); its energy uses \(\rho^2|W_{\bar z}|^2\). These derivative magnitudes agree for a real bump. Variation of the coefficient and density on the disk gives the \(O(h)\) error, uniformly over cotangent vectors of unit supremum norm.

We need at most \(N_g=4n+1\) disks: \(2n+1\) for the convex representation and two copies of the atom basis. Fix \(R_g>1\) so that hyperbolic area comparison allows \(N_g\) centers to be chosen successively, each within \(R_gh\) of its prescribed center and at distance at least \(2h\) from earlier centers. This choice depends only on \(N_g\). Require \((R_g+1)h<\sigma_X/4\), so the disks used in the construction are embedded. Moving centers and applying \eqref{eq:uniform-bump-atom} changes every response column by at most \(L_gh\) in \(\gamma_X\), for a now fixed genus-dependent \(L_g\).

For either resulting basis \(F(c)=\sum_i c_iq_X(b_i)\), \eqref{eq:uniform-atom-basis-inverse} gives
\[
\|E^{-1}(F-E)\|_{\ell^\infty\to\ell^\infty}\le nL_gh.
\]
Choose once and for all
\begin{equation}\label{eq:localization-constant-choice}
0<\eta_g\le
\min\left\{1,\frac1{8(R_g+1)},\frac1{2nL_g}\right\}.
\end{equation}
Both basis matrices are then invertible, with
\[
\|F^{-1}z\|_\infty\le2\gamma_X(z),
\qquad \gamma_X(F(c))\ge\tfrac12\|c\|_\infty.
\]
This fixes \(\kappa_g=1/2\), independently of every later smallness threshold.

Let \(v_A\) be the sum of the convex-representation bumps with coefficients equal to square roots of their weights. Disjoint support gives
\[
\|v_A\|_2^2=1,\qquad
\gamma_X(q_X(v_A)-u)\le L_gh.
\]
Use the first perturbed basis to write
\(u-q_X(v_A)=\sum_i c_iq_X(b_i)\), where
\(\sum_i|c_i|\le2nL_gh\le1\).
Choose \(z_i^2=c_i\) and set \(v_0=v_A+\sum_i z_ib_i\). Disjoint support makes \(q_X(v_0)=u\), adds at most \(2nL_gh\) to its energy, and bounds its supremum norm by \(C/h\). Fix \(A_g\) to dominate these energy and supremum constants. The second perturbed basis supplies the \(a_i\) and \eqref{eq:uniform-response-basis}.
\end{proof}

\begingroup
\small
\bibliographystyle{amsplain}
\bibliography{reference}
\endgroup
\end{document}